\documentclass{siamart190516}
\usepackage{amssymb}
\usepackage{url}
\usepackage{hyperref}
\usepackage{srcltx}
\usepackage{xcolor}
\usepackage{amsmath}
\usepackage{epsfig}

\def\Id{{\mathrm{Id}}}
\def\n{\mathrm{nom}}
\def\bT{{\mathbf{T}}}

\def\Id{{\mathrm{Id}}}

\def\three?{3}
\def\four?{4}
\def\clsr{{\mathrm{cl}}}
\def\ten?{10}

\def\Image{{\mathrm{Im}}}

\newtheorem{remark}{Remark}[section]

\def\norm2to2{{\|\cdot\|_{2,2}}}
\def\Prob{\hbox{\rm Prob}}

\def\bE{{\mathbf{E}}}

\def\inter{\hbox{\rm  int}}

\def\Diag{\hbox{\rm  Diag}}
\def\Prob{\hbox{\rm  Prob}}

\def\Opt{\hbox{\rm Opt}}

\def\Conv{\hbox{\rm  Conv}}

\def\Tr{{\mathop{\hbox{\rm  Tr}}}}
\def\cA{{\cal A}}
\def\cB{{\cal B}}
\def\cC{{\cal C}}

\def\cE{{\cal E}}
\def\cF{{\cal F}}

\def\cL{{\cal L}}

\def\cN{{\cal N}}

\def\cQ{{\cal Q}}
\def\cR{{\cal R}}

\def\cT{{\cal T}}
\def\cU{{\cal U}}
\def\cV{{\cal V}}
\def\cW{{\cal W}}
\def\cX{{\cal X}}

\def\cZ{{\cal Z}}

\def\B{{\cal B}}

\def\R{{\cal R}}
\def\S{{\cal S}}
\def\T{{\cal T}}

\def\V{{\cal V}}
\def\W{{\cal W}}
\def\X{{\cal X}}

\def\rank{{\mathop{\hbox{\rm  Rank}}}}
\def\Ker{{\hbox{\rm  Ker}\,}}

\def\bS{{\mathbf{S}}}

\def\e{{\hbox{\rm e}}}

\def\qed{\ \hfill$\square$\par\smallskip}

\def\bR{{\mathbf{R}}}

\def\argmin{\mathop{\hbox{\rm argmin}}}

\def\Risk{{\hbox{\rm Risk}}}

\newcommand{\ov}[1]{{\overline{#1}}}

\newcommand{\rf}[1]{~(\ref{#1})}
\newcommand{\be}{\begin{eqnarray}}
\newcommand{\ee}[1]{\label{#1}\end{eqnarray}}
\newcommand{\nn}{\nonumber \\}
\newcommand{\ese}{\end{eqnarray*}}
\newcommand{\bse}{\begin{eqnarray*}}
\newcommand{\third}{ \mbox{\small$\frac{1}{3}$}}
\newcommand{\half}{ \mbox{\small$\frac{1}{2}$}}
\newcommand{\hide}[1]{{}}
\newcommand{\aic}[2]{{\color{blue}~#2}}
\newcommand{\anc}[2]{{\color{violet}~#2}}
\title{Tight Computationally Efficient Approximation of Matrix Norms with Applications}
\author{
Anatoli Juditsky
\footnotemark[2]\ \footnotemark[5]
\and Georgios Kotsalis\footnotemark[3]\ \footnotemark[6]
\and Arkadi Nemirovski\footnotemark[4]\ \footnotemark[5]}

\begin{document}\maketitle
\renewcommand{\thefootnote}{\fnsymbol{footnote}}
\footnotetext[2]{LJK, Universit\'e Grenoble Alpes, 700 Avenue Centrale,  38401 Domaine Universitaire de Saint-Martin-d'Hères, France,
{\tt anatoli.juditsky@univ-grenoble-alpes.fr}}
\footnotetext[3]{Georgia Institute of Technology, Atlanta, Georgia 30332, USA, {\tt gkotsalis3@gatech.edu}}
\footnotetext[4]{Georgia Institute of Technology, Atlanta, Georgia 30332, USA, {\tt nemirovs@isye.gatech.edu}}
\footnotetext[5]{Research of this author was supported by MIAI {@} Grenoble Alpes (ANR-19-P3IA-0003).}
\footnotetext[6]{Research of this author was supported by NIFA grant 2020-67021-31526.}
\renewcommand{\thefootnote}{\arabic{footnote}}
\date{}

\begin{abstract}
We address the problems of computing operator norms of matrices induced by given norms on the argument  and the image space. It is known that aside of a fistful of ``solvable cases,'' most notably, the case when both given norms are Euclidean, computing operator norm of a matrix is NP-hard. We specify rather general families of
 norms on the argument and the images space (``ellitopic'' and ``co-ellitopic,'' respectively) allowing for reasonably tight computationally efficient upper-bounding of the associated operator norms. We extend these results to bounding ``robust operator norm of uncertain matrix with box uncertainty,'' that is, the maximum of operator norms of matrices representable as a linear combination, with coefficients of magnitude $\leq1$, of a collection of given matrices. Finally, we consider some applications of norm bounding, in particular, (1) computationally
  efficient synthesis of affine non-anticipative finite-horizon control of discrete time linear dynamical systems under bounds on the peak-to-peak gains, (2) signal recovery with uncertainties in sensing matrix, and (3) identification of parameters of time invariant discrete time linear dynamical systems via noisy observations of states and inputs on a given time horizon, in the case of ``uncertain-but-bounded'' noise varying in a box.

\end{abstract}
\section{Introduction}\label{sintro}
In this paper, our theoretical focus is on two problems as follows:
\par
{\bf A.} [{approximating} operator norms] {\sl Given norms $\|\cdot\|_\cX$ and
$\|\cdot\|_{\cB}$ with unit balls $\cX\subset\bR^n$ and  ${\cB}\subset \bR^m$, estimate the induced norm
$
\|A\|_{\cB,\cX}:=\max_{x:\|x\|_\cX\leq1}\|Ax\|_\cB$
of an $m\times n$ matrix $A$};
\par
{\bf B.} [{approximating} robust norm of uncertain matrix with box uncertainty] {\sl With $\|\cdot\|_\cX$, $\|\cdot\|_{\cB}$ as in {\bf A}, given an ``uncertain $m\times  n$ matrix with box uncertainty''---set of the form
$
\cA=\Big\{A_\n+\sum_{s=1}^S\epsilon_sA_s:\|\epsilon\|_\infty\leq 1\Big\},
$ ($A_\n,A_1,...,A_S\in\bR^{m\times n}$),
estimate the robust norm
$
\|\cA\|_{\cB,\cX}=\max_{A\in\cA}\|A\|_{\cB,\cX}.
$
of the uncertain matrix $\cA$.}
\par
Applications motivating our interest in these problems will be discussed later; we start with outlining the research status of these problems as ``academic entities'' and our related results.
\par$\bullet$
Aside of {few}  special cases, e.g., the case of the spectral norm ($\cX$ and $\cB$ are unit Euclidean balls in the respective spaces), {\bf A} is NP-hard; this is so, e.g.,
when $\|\cdot\|_\cX=\|\cdot\|_p$, $\|\cdot\|_\cB=\|\cdot\|_r$, and $p\geq 2\geq r\geq1$ with $p\neq r$ \cite{Daureen}. {\bf B} is NP-hard already when $\cB$, $\cX$
are unit Euclidean balls, $A_\n=0$, and $A_s$ are restricted to be symmetric matrices of rank 2 \cite{BTNMC}.
Hardness of {\bf A}, {\bf B} makes it natural to look for efficiently computable reasonably tight upper bounds on the norms in question. Below we build these bounds for the case
where $\cX$ and {\sl the polar $\cB_*$} of $\cB$ are {\sl ellitopes.} \begin{quote}
Sufficient for our current purposes example of an ellitope in $\bR^k$ is a bounded set $\cZ$ cut of $\bR^k$ by convex constraint on the vector  $[z^TP_1z;...;z^TP_Jz]$ of values of convex homogeneous quadratic forms of $z$:
    $
    \cZ=\{z\in\bR^k: \exists t\in\cT: z^TP_jz\leq t_j,j\leq J\},
    $
    where $P_j\succeq0$, $\sum_jP_j\succ0$, and $\cT$ is a convex compact subset of $\bR^J_+$ with a nonempty interior which is monotone, i.e., $0\leq t'\leq t\in\cT$ implies that $t'\in\cT$. A simple example x is the intersection of finitely many ellipsoids/elliptic cylinders centered at the origin.
    \end{quote}
We demonstrate that in the ellitopic case one can build efficiently computable upper bounds $\Phi(A)$ on $\|A\|_{\cB,\cX}$ and $\Psi(A_1,...,A_N)$ on $\|\cA\|_{\cB,\cX}$ which are convex in $A$, resp., in $(A_1,...,A_N)$, such that
\begin{subequations}
\begin{align}\label{eq0a}
 \|A\|_{\cB,\cX}&\leq \Phi(A)\leq O(1)\sqrt{\ln(2K)\ln(2L)}\|A\|_{\cB,\cX}, \\
   \|\cA\|_{\cB,\cX}&\leq \Psi(A_1,...,A_N)\leq O(1)\sqrt{\ln(2K)\ln(2L)}\vartheta(\kappa)\|\cA\|_{\cB,\cX}
   \label{eq0b}
\end{align}
\end{subequations}
where $K$ and $L$ are {\em ellitopic sizes} (numbers of quadratic forms in the description) of $\cX$ and $\cB_*$, $\kappa$ is the maximum of ranks of $A_i$, and $\vartheta(\cdot)$ is a certain universal function of $\kappa$.
    \par
$\bullet$ {\sl Relation to existing literature, problem {\bf A}.} {\bf  A} is the problem of maximizing a quadratic (specifically, bilinear) form on $\cB_*\times\cX$, and there exists significant literature on tractable relaxations, semidefinite and alike, of these problems. To the best of our knowledge, the most advanced existing results are those in the seminal papers \cite{NesSDP,YuLp} of Yu. Nesterov. As applied to {\bf A}, those results, in our present language, state that when the positive semidefinite matrices participating in description of $\cX$ and $\cB_*$ are diagonal, the appropriate  efficiently computable relaxation bound on $\|A\|_{\cB,\cX}$ (which in fact is nothing but the bound $\Phi$ participating in \rf{eq0a}) is tight within {\sl absolute constant factor}
(for details, see Remark \ref{nesterov}). It should be stressed that ``tightness within an absolute constant''  heavily exploits diagonality of the matrices describing $\cX$ and $\cB_*$; in the case of general ellitopes, logarithmic tightness factors in \rf{eq0a} seem to be unavoidable.\footnote{For instance, it was shown in \cite{NemRoTer} that when $\|\cdot\|_\cB=\|\cdot\|_2$ and $\cX$ is the intersection of $K$ ''stripes'' centered at the origin (i.e., the corresponding positive semidefinite matrices are of rank 1), the relaxation bounds in question can indeed be larger than the true quantity by factor $O\big(\sqrt{\ln K}\big)$.}
    \par
    The results on tight computationally tractable upper-bounding of maxima of quadratic forms over general-type ellitopes (same as the notion of an ellitope itself)
    originate from \cite{JudNem2018} and are further developed in \cite{STOPT}. As compared to those results, dealing with bilinear rather than with general quadratic forms allows
    us below to refine the analysis, and, as a result, to reduce the tightness factor in \rf{eq0a} to $O\sqrt{\ln(2K)\ln(2L)}$ instead of $O(1)\ln(K+L)$ guaranteed by \cite{STOPT}.
    \par
$\bullet$ {\sl Relation to existing literature, problem {\bf B}.} The only known to us preceding results on bounding robust norms of uncertain matrices deal
with the spectral norm ($\cX$ and $\cB$ are unit Euclidean balls), in which case the tightness factor in \rf{eq0b} boils down to
$\vartheta(2\kappa)$; these results can be easily derived from  the ``Matrix Cube Theorem'' in \cite{BTNMC}.
\paragraph{Applications} While {\bf A} and {\bf B} look legitimate academic problems, and the outlined results---legitimate academic results, the actual motivation for what follows stems from specific applications of problems {\bf A} and {\bf B} we are about to consider.
\par
Our principal motivation for problem {\bf A} comes from control and is the necessity to handle {\sl peak-to-peak} design specifications in synthesis of linear controllers. Specifically, given a linear dynamical system\\
\centerline{$
x_{t+1}=A_tx_t+B_tu_t+D_td_t,\;x_0=z,\,\,y_t=C_tx_t+E_td_t
$}
with states $x_t$, controls $u_t$, observed outputs $y_t$, and external disturbances $d_t$, we want to build an affine non-anticipating controller $
u_t=g_t+\sum_{\tau=0}^tG^t_\tau y_\tau
$
in such a way that the trajectory $w^N=\{x_t,1\leq t\leq N;\,y_t,u_t,0\leq t<N\}$ of the closed loop system on a given time horizon satisfies a given set of design specifications. With smart nonlinear reparameterization of affine non-anticipating controllers (passing from affine output-based control to the control which is affine in {\em purified outputs}, see \cite{KLN} and references therein), the system trajectory becomes affine function of the initial state $z$ and the sequence $d^N=[d_0;...;d_{N-1}]$ of external disturbances, with the matrices and constant terms in these affine functions  {\sl affine in the vector $\chi$ of controller's parameters} varying in certain $\bR^\nu$.  Bi-affinity of $w^N$ in $(d^N,z)$ and in $\chi$ is the key to computationally efficient processing of design specifications of appropriate structure3.
In this paper, we address an important (and considered as difficult in control) specification, namely, {\sl peak-to-peak gain} defined as follows.\footnote{For the sale of definiteness, we focus on ``disturbance-to-state'' peak-to-peak gain; peak-to-peak gains from disturbance to controls, or to outputs, or from initial state to states, etc., are defined similarly and can be processed in the same way.} Let us fix some a norm $\|\cdot\|_{(d)}$ on the space where the disturbances $d_t$ live, and norm $\|\cdot\|_{(x)}$ on the space where the states $x_t$ live. We equip the space $D^N$ of disturbance sequences $d^N=[d_0;...;d_{N-1}]$ with the norm $\|d^N\|_{d,\infty}=\max_t\|d_t\|_{(d)}$, and the space $X^N$ of state trajectories $x^N=[x_1;...;x_N]$ with the norm $\|x^N\|_{x,\infty}=\max_t\|x_t\|_{(x)}$. With affine in purified outputs controller $\chi$,
$x^N$ is an affine function of $d^N$ and $z$; let $X[\chi]$ be the matrix of coefficients at $d^N$ in this affine dependence. Peak-to-peak disturbance-to-state gain stemming from $\|\cdot\|_{(d)}$ and $\|\cdot\|_{(x)}$ is, by definition, the norm of $X[\chi]$ induced by the norms $\|d^N\|_{d,\infty}$ and $\|x^N\|_{x,\infty}$, and the corresponding design specification is just an upper bound on this gain. Since $X[\chi]$, as was already mentioned, is affine in $\chi$, this specification is a convex constraint on $\chi$.
However, this constraint can be difficult to handle because the operator norm in question is typically difficult to compute (this is so already when $\|\cdot\|_{(d)}$ and $\|\cdot\|_{(x)}$ are $\|\cdot\|_2$-norms). In such case, we can utilize our results on problem {\bf A} to safely approximate the design specification in question by replacing difficult-to-compute induced norm of $X=X[\chi]$ by its efficiently computable convex in $X$ and reasonably tight upper bound, as explained in details in Section \ref{peak-to-peak}.
\par
Our main motivating application for problem {\bf B} is identification of parameters $A$ of discrete time linear time invariant dynamical system\\
\centerline{$
x_{t+1}=A[x_t;r_t],
$}
from corrupted by noise observations of states $x_0,...,x_N$ and inputs $r_0,...,r_{N-1}$  on a given time horizon. We focus on the case of {\sl
uncertain-but-bounded noise}, in which deviations of entries in observations from the actual values of the corresponding entries in $x_t$ and $r_t$ are bounded in magnitude. We discuss an approach (to the best of our knowledge, new),heavily utilizes our results on problem {\bf B}, to computationally efficient identification of $A$ and to generating on-line upper bounds on recovery errors.
\par
Note that there is some literature on the first, and huge literature on the second of the just outlined applications. Instead of positioning our results with respect to this literature in the introduction,
we find it more productive to postpone this positioning till appropriate parts of the main body of the paper.
\par
{\sl Structure} of the paper is as follows. Section \ref{prelim} presents background on ellitopes. Section \ref{mainresellnorms} is devoted to problem {\bf A}, and Section \ref{secrobnorms}---to problem {\bf B}. Technical proofs are relegated to the appendix, where we present additional results on system identification, same as describe how our results can be extended from ellitopes to an essentially wider family of sets---{\sl spectratopes}.

\section{Preliminaries: ellitopes and spectratopes}\label{prelim}
Ellitopes and their extensions, spectratopes, introduced in \cite{STOPT}, are convex compact sets well-suited for tight upper-bounding maxima of quadratic forms over the sets. To make the paper more readable,
in its main body we focus on ellitopes; (always straightforward) extensions to spectratopes are relegated to Appendix.

\subsection{Ellitopes: definition and basic examples}\label{ss21} A {\sl basic ellitope} is a set $\cW$ represented as
\begin{equation}\label{2020ell1}
\cW=\{w\in\bR^p:\exists t\in\cT: w^TT_kw\leq t_k,\,1\leq k\leq K\}
\end{equation}
where $T_k\succeq0$, $k\leq K$, $\sum_kT_k\succ0$, and $\cT$ is a convex computationally tractable compact {\sl monotone} subset of $\bR^K_+$ with $\inter\cT\neq\emptyset$, monotonicity meaning that when $0\leq t\leq t'$ and $t'\in\cT$, we have $t\in\cT$ as well.
\par
An {\em  ellitope} $\cX$ is a linear image of a basic ellitope:
\begin{equation}\label{2020ell2}
\cX=P\cW=\{x\in\bR^n:\exists w\in\cW:x=Pw\} \hbox{\  with $\cW$ given by (\ref{2020ell1})}
\end{equation}
 We call $K$ {\sl ellitopic size} of ellitopes (\ref{2020ell1}) and (\ref{2020ell2}).
\par
Clearly, every ellitope is a convex compact set symmetric w.r.t. the origin;  a basic ellitope, in addition, has a nonempty interior.
\paragraph{Examples} {\bf A.} Bounded intersection $\cX$ of $K$ centered at the origin ellipsoids/elli\-p\-tic cylinders $\{x\in \bR^n:x^TT_kx\leq1\}$ [$T_k\succeq0$] is a basic ellitope:\\
\centerline{$
\cX=\{x\in \bR^n:\exists t\in\cT:=[0,1]^K: x^TT_kx\leq t_k,\,k\leq K\}
$}
In particular, the unit box $\{x\in \bR^n:\|x\|_\infty\leq1\}$ is a basic ellitope.\\
{\bf B.} A $\|\cdot\|_p$-ball in $\bR^n$ with $p\in[2,\infty]$ is a basic ellitope:\\
\centerline{$
\{x\in\bR^n:\|x\|_p\leq1\} =\{x:\exists t\in\cT=\{t\in\bR^n_+,\|t\|_{p/2}\leq 1\}:\underbrace{ x_k^2}_{x^TT_k x}\leq t_k,\,k\leq K\}.
$}
Ellitopes admit fully algorithmic "calculus:" this family is closed with respect to basic operations preserving convexity
and symmetry w.r.t. the origin, e.g., taking finite intersections, linear images, inverse images under linear embedding, direct products, arithmetic summation
 (for details, see \cite[Section 4.6]{STOPT}); what is missing, is taking convex hulls of finite unions.
 \subsection{Bounding maximum of  quadratic form over an ellitope}\label{bounding_ellitope}
 The starting point of what follows is the problem
 \begin{equation}\label{quadprob}
\Opt_*(C)=\max_{x\in\cX}x^TCx,\,\,C\in\bS^n
\end{equation}
of maximizing a homogeneous quadratic form  over a convex compact set $\cX\subset\bR^n$. It is well known that basically the only generic case when the problem is easy is the one where $\cX$ is an ellipsoid. It is shown in \cite{STOPT}
that when $\cX$ is an ellitope, (\ref{quadprob}) admits reasonably tight efficiently computable upper bound. Specifically, when $\cX$ is given by (\ref{2020ell2}),
$\lambda\in\bR^k_+$ is such that $P^TCP\preceq\sum_k\lambda_kT_k$ and $x\in \cX$, one has for some $t\in\cT$\\
\centerline{$
x^TCx=w^TP^TCPw\leq w^T[\sum_k\lambda_k T_k]w \leq \sum_k\lambda_kt_k,
$}
implying the validity of the implication\\
\centerline{$
\lambda\geq0,\;P^TCP\preceq\sum_k\lambda_kT_k\;\Rightarrow \;\Opt_*(C)\leq\phi_{\cT}(\lambda):=\max_{t\in\cT}\lambda^Tt,
$}
and thus---the first claim of the following
{\begin{theorem}\label{2020Prop4.6} {\rm \cite[Proposition 4.6]{STOPT}} Given ellitope {\rm (\ref{2020ell2})} and
a matrix $C\in\bS^n$, consider the quadratic maximization problem {\rm (\ref{quadprob})}
along with its relaxation
\begin{equation}\label{2020eq10}
\Opt(C)=\min_\lambda\left\{\phi_{\cT}(\lambda): \lambda\geq0,P^TCP\preceq{\sum}_k\lambda_kT_k\right\}
\end{equation} The problem  is computationally tractable and solvable, and $\Opt(C)$ is an efficiently computable upper bound on $\Opt_*(C)$.
This upper bound is reasonably tight:
\[
\Opt_*(C)\leq \Opt(C)\leq3\ln(\sqrt{3}K)\Opt_*(C).
\]
\end{theorem}}\noindent
To the best of our knowledge, the first result of this type was established in \cite{NemRoTer} for $\cX$ which is an intersection of $K$ concentric elliptic cylinders/ellipsoids; in this case, (\ref{quadprob}) becomes a special case of quadratic quadratically constrained optimization problem, and (\ref{2020eq10}) is the standard Shor's semidefinite relaxation (see, e.g., \cite[Section 4.3]{LMCO}) of this problem. In \cite{NemRoTer} it is shown
that the ratio $\Opt(C)/\Opt_*(C)$ indeed can be as large as $O(\ln(K))$, even when all $T_k=a_ka_k^T$ are of rank 1 and  $\cX$ is the polytope $\{x:|a^T_kx|\leq1,k\leq K\}$.

\section{Bounding operator norms}\label{mainresellnorms}
As stated in Introduction, one of the subjects of this paper is tight efficiently computable upper-bounding of the operator norm\\
\centerline{$
\|A\|_{\cB,\cX}=\max\limits_x\left\{\|Ax\|_\cB: \|x\|_\cX\leq 1\right\}
$}
of a linear mapping $x{\to} Ax:\bR^n\to\bR^m$ induced by norms $\|\cdot\|_\cX$ and $\|\cdot\|_\cB$ on the argument and the destination spaces,
with $\|\cdot\|_\cU$ standing for the norm with unit ball $\cU$. Our approach works for the case when $\cX$ and {\sl the polar} $\cB_*$
of $\cB$ are ellitopes with nonempty interiors:
{\small\begin{equation}\label{2020ell2*}
\begin{array}{rcl}
\cX&=&P\cW=\{x\in\bR^n:\exists w\in\cW:x=Pw\},\,\inter \cX\neq\emptyset,\\
&&\cW=\{w\in\bR^p: \exists t\in\cT:w^TT_kw\leq t_k,k\leq K\}:\\
\cB_*&:=&\{v\in\bR^m:v^Ty\leq1\,\forall y\in\cB\}\\
&=&Q\cZ=\{y\in\bR^m:\exists z\in\cZ:y=Qz\},\,\inter \cB_*\neq\emptyset,\\
&&\cZ=\{z\in\bR^q: \exists r\in\cR: z^TR_\ell z\leq r_\ell,\ell\leq L\}
\end{array}
\end{equation}}
with $T_k$, $\cT$, $R_\ell$, $\cR$ as required in the definition of a basic ellitope.
\par
Under the assumptions just introduced, $\|A\|_{\cB,\cX}$  is the maximum of a quadratic form over a basic ellitope $\cZ\times\cW$:
{\small\bse
\|A\|_{\cB,\cX}=\max\limits_{x\in\cX}\|Ax\|_\cB=\max\limits_{y\in\cB_*,x\in\cX}y^TAx=\max\limits_{w\in\cW,z\in\cZ}z^TQ^TAPw\\
={{1\over 2}}\max\limits_{[z;w]\in\cZ\times\cW}[z;w]^T\left[\begin{array}{c|c}&Q^TAP\cr\hline
P^TA^TQ\cr\end{array}\right][z;w].
\ese}\noindent
In this case relaxation (\ref{2020eq10}) provides efficiently computable upper bound on $\|A\|_{\cB,\cX}$. Immediate computation taking into account the direct product structure of the ellitope $\cZ\times\cW$ and bilinearity of the quadratic form we are maximizing over this ellitope shows that this bound is
\begin{equation}\label{normbound}
\Opt(A)=\min\limits_{\lambda,\upsilon}\left\{\phi_\cT(\lambda)+\phi_\cR(\upsilon):\lambda\geq0,\upsilon\geq0,\left[\begin{array}{c|c}
\sum_\ell\upsilon_\ell R_\ell&{1\over 2}Q^TAP\cr\hline{1\over 2}P^TA^TQ&\sum_k\lambda_kT_k\cr\end{array}\right]\succeq0\right\}.
\end{equation}
Note that $\Opt(A)$ clearly is a convex function of $A$, and Theorem \ref{2020Prop4.6} implies that\\
\centerline{$
\|A\|_{\cB,\cX}\leq \Opt(A)\leq  3\ln(\sqrt{3}[K+L])\|A\|_{\cB,\cX}.
$}
Our main goal is to demonstrate that the latter bound can be refined{.}
{\begin{theorem}\label{verynewRelaxationTheorem} In the case of {\rm \rf{2020ell2*}} one has
{\small\begin{equation}\label{betterbound}
\|A\|_{\cB,\cX}\leq\Opt(A)\leq
{\varsigma(K,L)}\|A\|_{\cB,\cX},\;\;\varsigma(K,L)=\left\{\begin{array}{ll}3\sqrt{\ln(4K)\ln(4L)},&\max[K,L]>1\\
1,&K=L=1.
\end{array}\right.
\end{equation}}\noindent
\end{theorem}}\noindent
\begin{remark}~\label{nesterov}{\rm
Results of \cite{NesSDP,YuLp} imply that in some cases the tightness factor $\kappa$ in (\ref{betterbound}) can be improved to an {\sl absolute} constant. Specifically,
\par
1) In the case of  \rf{2020ell2*} {\sl with diagonal matrices $T_k$ and $R_\ell$}, it follows from \cite[Theorem 13.2.1]{NesSDP} that one can take $\varkappa={\pi\over 4-\pi}\approx3.660$
\par
2) When $\|\cdot\|_\cX=\|\cdot\|_p$, $\|\cdot\|_\cB=\|\cdot\|_r$ with $\infty\geq p\geq2$, $1\leq r\leq2$ (this is a special case of 1)), Nesterov \cite{NesSDP,YuLp} proved that the upper bound
\begin{equation}\label{pq}
{1\over 2}\min\limits_{\lambda,\mu}\left\{\|\lambda\|_{{p\over p-2}}+\|\mu\|_{{r\over 2-r}}:\left[\begin{array}{c|c}\Diag\{\mu\}&A\cr\hline A^T&\Diag\{\lambda\}\cr\end{array}\right]\succeq0\right\}
\end{equation}
on $\|A\|_{p\to r}:=\max_{\|x\|_p\leq1}\|Ax\|_r$ (this bound coincides with $\Opt(A)$ when $\cX$ is the ellitope $\{x:\|x\|_p\leq1\}$, and $\cB_*$ is the ellitope $\{v:\|v\|_{{r\over r-1}}\leq 1\}$) is tight within (even better than in 1))
\begin{itemize}
\item factor ${\pi\over 2\sqrt{3}-2\pi/3}\approx2.2936$ in the entire range $p\in[2,\infty],r\in[1,2]$,
\item factor $\sqrt{\pi/2}\approx1.2533$ when $p=2$ and $r\in[1,2]$.\footnote{Using the identity $\|A\|_{\cB,\cX}=\|A^T\|_{\cX_*,\cB_*}$, where $\cX_*$ is the polar of $\cX$ (as is immediately seen, this identity is respected by our bounding scheme), we see that $\Opt(A)$ is within $\sqrt{\pi/2}$ from $\|A\|_{p\to r}$ when $p\geq2$ and $r=2$.}
\end{itemize}
\par\noindent
Needless to say, when $p=r=2$, the tightness factor is 1. In addition, it is shown in \cite{Daureen} that in the range $\infty\geq p\geq2$, $1\leq r\leq 2$ bound (\ref{pq})
is {\sl exactly equal} to the corresponding norm of $A$ {\sl for entrywise nonnegative matrices}.
}\end{remark}
%%%%%%%%%%%%%%%%%%%%
\par
Note that there is a simple case when $\Opt(A)=\|A\|_{\cB,\cX}$---the one where $A$ is a row vector,  $\cB=[-1,1]\subset\bR$, and, therefore,
$$
\|A\|_{\cB,\cX}=\max_{x\in\cX} Ax.
$$
Our bounding is intelligent enough to recognize this situation. Indeed, in the case in question (\ref{normbound}) reads
$$
\Opt(A)=\min\limits_{\lambda,\upsilon}\left\{\phi_\cT(\lambda)+\upsilon:\lambda\geq0,
\left[\begin{array}{c|c}
\upsilon&{1\over 2}AP\cr\hline{1\over 2}P^TA^T&\sum_k\lambda_kT_k\cr\end{array}\right]\succeq0\right\}
$$
while, by Lagrange duality,
{\small$$
\begin{array}{rl}
\max\limits_{x\in\cX} Ax&=\max\limits_{w,t}\left\{APw:w^TT_kw\leq t_k,k\leq K, t\in\cT\right\}\\
&=\min\limits_{\lambda\geq0}\max\limits_{t\in \cT,w}\left\{APw+\sum_k\lambda_kt_k-w^T\left[\sum_k\lambda_kT_k\right]w\right\}\\
&=\min\limits_{\lambda\geq0}\max\limits_{w}\left\{\phi_\cT(\lambda)+APw-w^T\left[\sum_k\lambda_kT_k\right]w\right\}\\
&=
\min\limits_{\lambda\geq0,\upsilon}\left\{\phi_\cT(\lambda)+\upsilon: \upsilon-APw+w^T\left[\sum_k\lambda_kT_k\right]w\geq0\,\forall w\right\}\\
&=\min\limits_{\lambda,\upsilon}\left\{\phi_\cT(\lambda)+\upsilon:\lambda\geq0,
\left[\begin{array}{c|c}
\upsilon&{1\over 2}AP\cr\hline{1\over 2}P^TA^T&\sum_k\lambda_kT_k\cr\end{array}\right]\succeq0\right\}=\Opt(A).
\end{array}
$$}\noindent
To put this immediate observation into a proper perspective, see Section \ref{extension}.
\par
The just outlined results are stronger than what in the case in question is stated by Theorem \ref{verynewRelaxationTheorem}. This being said, it can be proved that in the full scope of the latter theorem, logarithmic growth of the tightness factor with $K,L$ is unavoidable.
\subsection{On the scope of Theorem \ref{verynewRelaxationTheorem}}\label{ssscope}
The scope of Theorem \ref{verynewRelaxationTheorem}{---the set of the matrix norms to which the theorem applies---}is restricted to the case when the norm in the argument space is {\sl simple ellitopic norm}, meaning that its unit ball is an ellitope, and the norm on the image space
is a {\sl simple co-ellitopic} norm, meaning that the polar of its unit ball is an ellitope. Clearly, simple co-ellitopic norms (s.co-e.n's) are exactly the conjugates of simple ellitopic norms (s.e.n.'s). These classes of norms allow for certain ``calculus'' stating that some standard operations with norms preserve their ellitopic/co-ellitopic type.
\par
{\sl Basic calculus of simple ellitopic norms} is as follows.
\begin{enumerate}
\item[E.1.] {[raw materials]} When $p\in[2,\infty]$, $\|\cdot\|_p$ is a s.e.n. on $\bR^n$,
\item[E.2.] {[taking finite maxima]} When $\|\cdot\|_{(k)}$, $k\leq K$, are s.e.n.'s on $\bR^n$, so is their maximum.
\item[E.3.] {[restriction to a linear subspace]} When $\|\cdot\|$ is a s.e.n. on $\bR^n$ and $y\to Ax:\bR^{n'}\to\bR^n$ is a linear embedding, $\|y\|':=\|Ay\|$ is a s.e.n. on $\bR^{n'}$
\item[E.4.] {[passing to factor-norm]} When $\|\cdot\|$ is a s.e.n. on $\bR^n$ and $x\mapsto Ax:\bR^n\to\bR^{n'}$ is an onto mapping, the factor-norm $\|y\|'=\min_x\{\|x\|:Ax=y\}$ is s.e.n. on $\bR^{n'}$
\item[E.5.]{[``aggregation'']} Let $\|\cdot\|_{(k)}$ be s.e.n. on $\bR^{n_k}$, $k\leq K$, and let $\cA$ be a monotone convex compact set with a  nonempty interior  in $\bR^K_+$. Then the norm on $\bR^{n_1}\times...\times\bR^{n_K}$ with the unit ball
    $$
    \cX=\{[x_1;...;x_K]\in\bR^{n_1}\times...\times\bR^{n_K}:\exists {\alpha\in \cA,}{\|x_k\|_{(k)}\leq \sqrt{\alpha_k},k\leq K}\}
    $$
    is s.e.n.
    For instance, when $p_k\in[2,\infty]$ and $p\in[2,\infty]$, the norm on $\bR^{n_1}\times...\times\bR^{n_K}$ given by $
    \|[x_1;...;x_K]\|=\|[\|x_i\|_{p_1};...;\|x_K\|_{p_K}]\|_p
    $
    is s.e.n.
\end{enumerate}
All these rules are immediate consequences of ``calculus of ellitopes'' \cite[Section 4.6]{STOPT}.\par
{\sl Basic calculus of simple co-ellitopic norms}  is as follows.
\begin{enumerate}
\item[cE.1.] {[raw materials]} When $r\in[1,2]$, $\|\cdot\|_r$ is a s.co-e.n. on $\bR^n$ (cf. E.1)
\item[cE.2.] {[taking sums]} When $\|\cdot\|_{(k)}$, $k\leq K$, are s.co-e.n.'s on $\bR^n$, so is their sum.
\end{enumerate}
{\small Indeed, the unit ball $B$ of the sum of norms with polars $B_i^*$ of the unit balls is
$$B=\left\{x:{\sum}_i{\max}_{y_i\in B_i^*}y_i^Tx\leq1\right\}=\left\{x:{\max}_{y=[y_1;...;y_K]\in B_1^*\times...\times B_K^*} x^T[y_1+...+y_K]\leq1\right\},
$$
that is, the polar $B^*$ of $B$ is the image of $B_1^*\times...B_K^*$ under a linear mapping. When all $B_k^*$ are ellitopes, so is their direct product, and therefore---its linear image $B^*$. Thus, the polar of $B$ is an } ellitope, as claimed.
\begin{enumerate}
\item[cE.3.] {[restriction to a linear subspace]} When $\|\cdot\|$ is a s.co-e.n. on $\bR^n$ and $y\to Ay:\bR^{n'}\to\bR^n$ is a linear embedding, $\|y\|':=\|Ay\|$ is a s.co-e.n. on $\bR^{n'}$
\end{enumerate}
{\small
Indeed, assuming that the polar $B^*$ of the unit ball of $\|\cdot\|$ is an ellitope, we have  $\|y\|'=\max_{z\in B^*}z^TAy$. That is, the polar of the unit ball of $\|\cdot\|'$ is the linear image $A^TB^*$ of $B^*$, which is an ellitope along with $B^*$.}
\begin{enumerate}
\item[cE.4.] {[passing to factor-norm]} When $\|\cdot\|$ is a s.co-e.n. on $\bR^n$ and $x\mapsto Ax:\bR^n\to\bR^{n'}$ is an onto mapping, the factor-norm $\|y\|'=\min_x\{\|x\|:Ax=y\}$ is s.co-e.n. on $\bR^{n'}$.
\end{enumerate}
{\small Indeed, assuming the polar $B^*$ of the unit ball of $\|\cdot\|$ to be an ellitope and denoting by $A^\dag$ the pseudoinverse of the onto mapping $A$, one has
\[\|y\|^\prime={\min}_{\delta\in\Ker A}\|A^\dag y+\delta\|
={\min}_{\delta\in\Ker A}{\max}_{z\in B^*}z^T[A^\dag y+\delta]={\max}_{z\in B^*\cap{\Image A^T}}[[A^\dag]^Tz]^Ty.
\]
Thus, the polar of the unit ball of $\|\cdot\|^\prime$ is a linear image of the intersection of ellitope $B^*$ with a linear subspace, and as such is an ellitope.}
\begin{enumerate}
\item[cE.5.] {[``aggregation'']} Let $\|\cdot\|_{(k)}$ be s.co-e.n. on $\bR^{n_k}$, $k\leq K$, and let $\cA$ be a monotone convex compact set with a  nonempty interior  in $\bR^K_+$. Then the norm on $\bR^{n_1}\times...\times\bR^{n_K}$ given by
    $$
    \|[x_1;...;x_K]\|=\max_{\beta\in\cA}{\sum}_k\sqrt{\beta_k}\|x_k\|_{(k)}
    $$
    is s.co-e.n. For instance, when $r_k\in[1,2]$ and $r\in[1,2]$, the norm on $\bR^{n_1}\times...\times\bR^{n_K}$ given by
    $
    \|[x_1;...;x_K]\|=\|[\|x_1\|_{r_1};...;\|x_K\|_{r_K}]\|_r
    $
    is s.co-e.n.
\end{enumerate}
    {\small Indeed, let $\|\cdot\|_{(k)}^*$ be the s.e.n.'s  conjugate to $\|\cdot\|_{(k)}$. Setting $\cA^{1/2}=\{[\alpha_1;...;\alpha_K]\geq0:[\alpha_1^2;\alpha_2^2;...;\alpha_K^2]\in\cA\},$  we get a convex compact monotone subset of $\bR^K_+$ such that
    the unit ball $\cB$ of $\|\cdot\|$ is
    $
    \cB=\{[x_1;...;x_K]:\phi_{\cA^{1/2}}([\|x_1\|_{(1)};...;\|x_K\|_{(K)}])\leq 1\}.
    $
Hence, as is immediately seen, the polar $\cB_*$ of $\cB$
    is
    $$
    \begin{array}{rcl}
    \cB_*&=&\{[y_1;...;y_K]: \sum_k \zeta_k\|y_k\|_{(k)}^*\leq 1\,\forall (\zeta\geq0:\sum_k\alpha_k\zeta_k\leq 1\,\forall \alpha\in\cA^{1/2})\}\\
    &=&
    \{[y_1;...;y_K]: \exists \alpha\in\cA^{1/2}:\|y_k\|_{(k)}^*\leq \alpha_k,\,k\leq K\},\\
    \end{array}
    $$
   that is, $\|\cdot\|_*$ is s.e.n. by E.5.}

\subsection{An extension}\label{extension}
The above results can be straightforwardly extended from the case when $\cB_*$ and $\cX$ are ellitopes onto a more general case. Specifically, assume that
\begin{itemize}
\item[\textbf{A}.] $\cX\subset\bR^n$ is a set with nonempty interior represented  as the convex hull of a finite union of ellitopes, or, which is the same,
\end{itemize}
\begin{equation}\label{A}
\begin{array}{rcl}
\cX&=&\Conv\{{\bigcup}_{i=1}^IP_i\cX_i\}=\left\{x=\sum_{i=1}^I \lambda_iP_ix_i:x_i\in\cX_i,\lambda_i\geq0,\sum_i\lambda_i=1\right\}\\
&=&\left\{x=\sum\limits_i P_ix_i:\sum\limits_i\|x_i\|_{\cX_i}\leq1\right\},
\end{array}
\end{equation}
\begin{itemize}
\item[]
where $\cX_i\subset\bR^{n_i}$ are basic ellitopes and $\|\cdot\|_{\cX_i}$ are s.e.n. on $\bR^{n_i}$ with unit balls $\cX_i$.
\end{itemize}
Under Assumption \textbf{A}, $\cX$ is a convex compact symmetric w.r.t. the origin subset of $\bR^n$ with $0\in\inter \cX$; as such, $\cX$ is the unit ball of a norm $\|\cdot\|_\cX$. In the sequel we refer to the norms of this structure as to {\sl ellitopic norms.} Clearly, every simple ellitopic norm is ellitopic, e.g., the block $\ell_\infty$ norm
$$
\|[x_1;...;x_I]\|=\max\limits_{i\leq I}\|x_I\|_{p_i}\eqno{[p_i\in[2,\infty]\,\forall i]}
$$
on the space $\bR^{n_1}\times...\times\bR^{n_I}$ is s.e.n. (by E.1 and E.5).
In fact, the family of ellitopic norms is much wider that the family of s.e.n.'s. For example,
\begin{itemize}
\item[{\bf E.}1.] When $\|x_i\|_{(i)}$ are ellitopic norms on $\bR^{n_i}$, $i\leq I$, the associated block $\ell_1/\|\cdot\|_{(\cdot)}$ norm
\begin{equation}\label{ell1}
\|[x_1;...;x_I]\|=\sum_{i=1}^I\|x_i\|_{(i)}
\end{equation}
on $\bR^{n_1}\times...\times\bR^{n_I}$ is ellitopic.
\end{itemize}
{\small Indeed, the unit ball $\cX_i$ of $\|\cdot\|_{(i)}$ is a convex subset of $\bR^{n_i}$ of the form $\Conv\Big\{{\bigcup}_{\nu=1}^{I_i}P_{i\nu}\cX_{i\nu}\Big\}$ with basic ellitopes $\cX_{i\nu}$. Specifying linear mappings $P_i$ from $\bR^{n_i}$ to  $\bR^{n_1}\times...\times\bR^{n_I}$ as the natural embeddings
$$
[P_ix_i]_s=\left\{\begin{array}{ll}0\in\bR^{n_s},&s\neq i\\
x_i,&s=i\\
\end{array}\right\},
$$
the unit ball $\cX$ of norm (\ref{ell1}) clearly is
$
\Conv\left\{{\bigcup}_{i\leq I, \nu_i\leq I_i}P_iP_{i\nu_i}\cX_{i\nu_i}\right\}.
$
Because, in addition, this set has a nonempty interior, (\ref{ell1}) is an ellitopic norm.}
\par
Note that the property to be ellitopic is inherited when passing to factor-norms (cf. E.4):
\begin{enumerate}
\item[{\bf E.}2] When $\|\cdot\|$ is an ellitopic norm and $y\mapsto Ay:\bR^n\to\bR^{n'}$ is an onto mapping, the factor-norm $\|x\|^\prime=\min_y\{\|y\|:Ay=x\}$ on $\bR^{n'}$ induced by $\|\cdot\|$ and $A$ is ellitopic.
\end{enumerate}
{\small
Indeed, if the unit ball $\cX$ of $\|\cdot\|$ is given by \rf{A} then the unit ball $\cX'$ of $\|\cdot\|^\prime$ is the convex compact set with a nonempty interior given by
$
\cX'=A\cX=\Conv\Big\{{\bigcup}_{i=1}^I[AP_i]\cX_i\Big\}
$}
\begin{enumerate}
\item[{\bf E.}3] Let $\|\cdot\|_{(\chi)}$ be an ellitopic norm on $\bR^{n_\chi}$, $\chi=1,2$. Then the norm $\|[x_1;x_2]\|=\max[\|x_1\|_{(1)},\|x_2\|_{(2)}]$ on $\bR^{n_1}\times\bR^{n_2}$ is ellitopic.
\end{enumerate}
{\small Indeed, if the unit ball of $\|\cdot\|_{(\chi)}$ is $\Conv\Big\{{\bigcup}_{i=1}^{I_\chi}P_{i,\chi}\cX_{i,\chi}\Big\}$, then the unit ball of $\|\cdot\|$ is $\Conv\Big\{{\bigcup}_{i_1\leq I_1,i_2\leq I_2}
\Diag\{P_{i_1,1},P_{i_2,2}\}[\cX_{i_1,1}\times\cX_{i_2,2}]\Big\}$, and $\cX_{i,1}\times\cX_{i,2}$ are basic ellitopes along with $\cX_{i,1}$, $\cX_{i,2}$.}\\
By {\bf E}.3,  if $\|\cdot\|_{(i)}$ are ellitopic norms on $\bR^{n_i}$, $i\leq I$, then the norm $\|[x_1;...;x_I]\|=\max_{i\leq I}\|x_i\|_{(i)}$
on $\R^{n_1+...+n_I}$ is ellitopic as well. Note, however, that the number of ellitopes involved in the description of this norm is the product,
over $i\leq I$, of the numbers of ellitopes in the description of norms  $\|\cdot\|_{(i)}$ and thus may explode exponentially fast as $I$ grows.
\par
Assume, next, that
\begin{itemize}
\item[\textbf{B.}] $\cB\subset\bR^m$ is a set with nonempty interior which is {\sl the polar} of a set of the structure described in \textbf{A}:
\end{itemize}
{\small\begin{equation}\label{B}
\cB=\{v\in\bR^m:\max\limits_{y\in\cB_*} v^Ty\leq 1\}, \;\cB_*=\Big\{y=\sum\limits_{j=1}^J\mu_jQ_jz_j,z_j\in\cZ_j,\mu_j\geq0,{\sum}_j\mu_j=1\Big\}
\end{equation}}
\begin{itemize}
\item[]
where $\cZ_j\subset\bR^{m_j}$ are basic ellitopes and $Q_j$, $\cZ_j$ are such that $\cB_*$ has a nonempty interior.
\end{itemize}
Under Assumption \textbf{B}, $\cB$ is a convex compact symmetric w.r.t. the origin subset of $\bR^n$ with $0\in\inter \cB$; as such, $\cB$ is the unit ball of a norm $\|\cdot\|_\cB$. In the sequel, we refer to norms of this structure as {\sl co-ellitopic.} Clearly, the conjugate of an ellitopic norm is co-ellitopic, and vice versa.
\par
Note that in the case of (\ref{B}) we have
{\small \begin{equation}\label{wehaveB}
\begin{array}{rl}
\|u\|_{\cB}&=\max\limits_{y\in\cB_*} u^Ty=\max\limits_{\{z_j,\mu_j\}}\left\{{\sum}_j\mu_ju^TQ_j:z_j\in\cZ_j,\mu_j\geq0\,\forall j,{\sum}_j\mu_j=1\right\}\\
&=\max\limits_j\max\limits_{z_j\in\cZ_j}u^TQ_jz_j=\max\limits_{j\leq J}\|Q_j^Tu\|_{\cZ^*_j}\\
\end{array}
\end{equation}}\noindent
where $\cZ^*_j$ is the polar of $\cZ_j$.
\par
Of course, every simple co-ellitopic norm is co-ellitopic. In fact, the family of co-ellitopic norms is much wider than the family of simple co-ellitopic norms due to the following  observations:
\begin{enumerate}
\item[{\bf cE.}1.] Maximum of finitely many co-ellitopic norms is co-ellitopic.
\end{enumerate}
{\small Indeed, if $\|\cdot\|_{(k)}$, $k\leq K$, are co-ellitopic norms on $\bR^n$, their conjugates $\|\cdot\|_{(k)}^*$ are ellitopic, implying by {\bf E}.1  that the  norm
$\|[y_1;...;y_K]\|=\sum_k\|y_k\|_{(k)}^*$ on $\bR^{Kn}$ is ellitopic, which by {\bf E.}2 implies that the factor-norm
$$
\|z\|_*=\min_{\{y_k\}}\left\{{\sum}_k\|y_k\|_{(k)}^*:{\sum}_ky_k=z\right\}
$$
is ellitopic. The unit ball of the latter norm is the convex compact set
\[\cB_*=\{z={\sum}_ky_k:{\sum}_k\|y_k\|_{(k)}^*\leq 1\},\] and the polar of this set
is
\begin{align*}
\cB&=\left\{x:\max_y\left\{[\sum\limits_ky_k]^Tx:\sum\limits_k\|y_k\|_{(k)}^*\leq1\right\}\leq 1\right\}\\
&=\left\{x:\max_{\lambda,y}\left\{[\sum\limits_ky_k]^Tx:\|y_k\|_{(k)}^*\leq\lambda_k,\sum\limits_k\lambda_k\leq 1\right\}\leq 1\right\}\\
&=\left\{x:\max_{\lambda}\left\{\sum\limits_k\left[\max_{y_k}\{x^Ty_k:\|y_k\|_{(k)}^*\leq\lambda_k\}\right],\sum\limits_k\lambda_k\leq 1\right\}\leq 1\right\}\\
&=\left\{x:\max_{\lambda\geq0}\left\{\sum\limits_k\lambda_k\|x\|_{(k)}:\sum\limits_k\lambda_k\leq 1\right\}\leq 1\right\}
=\{x:\max_k\|x\|_{(k)}\leq1\}.
\end{align*}\noindent
Thus, the norm $\max_k\|x\|_k$ is conjugate to the ellitopic norm $\|\cdot\|_*$ and  as such is co-ellitopic.}
\par
A closely related statement is
\begin{enumerate}
\item[{\bf cE.}2.] $\ell_\infty$-aggregation
\begin{equation}\label{ellinf}
\|[x_1;...;x_K]\|={\max}_{k\leq K} \|x_k\|_{(k)}
\end{equation}
of co-ellitopic norms $\|\cdot\|_{(k)}$ on $\bR^{n_k}$ is co-ellitopic.
\end{enumerate}
{\small
Indeed, as we have seen when justifying {\bf cE.}1, if $\|\cdot\|_{(k)}^*$ are ellitopic norms conjugate to $\|\cdot\|_{(k)}$, the norm
$$
\|[y_1;...;y_K]\|_*={\sum}_k\|y_k\|_{(k)}^*
$$
is ellitopic; clearly, norm (\ref{ellinf}) is conjugate to this ellitopic norm.}
\par
{The second observation is as follows.}
\begin{enumerate}
\item[{\bf cE.}3.] The restriction of a co-ellitopic norm onto a linear subspace is co-ellitopic.
\end{enumerate}
{\small Indeed, we should verify that if $x\mapsto Ax$ is an embedding of $\bR^{n^\prime}$ into $\bR^n$ and $\|\cdot\|$ is a co-ellitopic norm on $\bR^n$ then the norm $\|x\|'=\|Ax\|$ is co-ellitopic. This is immediate---by the standard properties of norms, under the circumstances, the norm conjugate to $\|x\|'$ is the factor-norm $\min_y\left\{\|y\|_*:A^Ty=x\right\}$ induced by the conjugate to $\|\cdot\|$ norm $\|\cdot\|_*$ on $\bR^n$. This conjugate is an ellitopic norm on $\bR^n$, and it remains to use {\bf E.}2.}
\begin{enumerate}
\item[{\bf cE.}4.] The sum of two co-ellitopic norms on $\bR^n$ is co-ellitopic.
\end{enumerate}
{\small Indeed, if $\|\cdot\|_{(\chi)}$, $\chi=1,2$, are co-ellitopic norms on $\bR^n$, and $\|\cdot\|_{(\chi)}^*$ are their conjugates, then the norm
$\|[x_1;x_2]\|_*=\max[\|x_1\|_{(1)}^*,\|x_2\|_{(2)}^*]$ is ellitopic norm on $\bR^{2n}$ by {\bf E}.3, so that its conjugate, which is $\|[x_1;x_2]\|_+=\|x_1\|_{(1)}+\|x_2\|_{(2)}$, is co-ellitopic. By {\bf cE}.3, the restriction of the latter norm on the subspace $\{[x_1;x_2]:x_1=x_2\}=[I_n;I_n]\bR^n$ also is co-ellitopic, and this restriction is nothing but the norm $\|x\|=\|x\|_{(1)}+\|x\|_{(2)}$.}
\paragraph{Simple observation} {Let $\|\cdot\|_\cX$ and $\|\cdot\|_\cB$ be norms with $\cX$ given by (\ref{A}) and $\cB$ given by (\ref{B}). Then}  the operator norm
of $A{\in\bR^{m\times n}}$ induced by the norms $\|\cdot\|_\cX$ and $\|\cdot\|_\cB$ on the  argument and image spaces can be computed as follows:
{\small\begin{equation}\label{asfollowsa}
\begin{array}{rcl}
\|A\|_{\cB,\cX}&=&\max_x\{\|Ax\|_\cB:\|x\|_\cX\leq1\}
=\max\limits_x\left\{\max\limits_j\|Q_j^TAx\|_{\cZ^*_j}:\|x\|_\cX\leq1\right\}\hbox{\ [see (\ref{wehaveB})]}\\
&=&\max\limits_x\left\{\max\limits_j\max\limits_{z_j}\{z_j^TQ_j^TAx:z_j\in\cZ_j\}:\|x\|_\cX\leq1\right\}\\
&=&
\max\limits_j\left\{\max\limits_{x,z_j}z_j^TQ_j^TAx:z_j\in\cZ_j,x\in \cX\right\}\\
&=&\max\limits_j\left\{\max\limits_{x,z_j}z_j^TQ_j^TAx:z_j\in\cZ_j,x\in\Conv\{\cup_iP_i\cX_i\}\right\}\\
&=&%\max\limits_j\left\{\max\limits_{z_j,x_1,...,x_I}z_j^TQ_j^TAP_ix_i:z_j\in\cZ_j,x_i\in\cX_i\right\}=
\max\limits_j\left\{\max\limits_{i}\left[\max\limits_{z_j\in\cZ_j,x_i\in \cX_i} z_j^TQ_j^TAP_ix_i\right]\right\}
=\max\limits_{i,j}\|Q_j^TAP_i\|_{ij}
\end{array}
\end{equation}}
where
\begin{equation}\label{asfollowsb}
\|Q_j^TAP_i\|_{ij}=
\max\limits_{z_j\in\cZ_j,x_i\in \cX_i} z_j^T[Q_j^TAP_i]x_i=\|Q_j^TAP_i\|_{\cZ_j^*,\cX_i}.
\end{equation}
Note that by the same token $\max_i\|Q_j^TAP_i\|_{ij}=\|Q_j^TA\|_{\cZ_j^*,\cX}$ and $\max_j\|Q_j^TAP_i\|_{ij}=\|AP_i\|_{\cB,\cX_i}$, so that in the case of (\ref{A}), (\ref{B}) it holds
$$
\|A\|_{\cB,\cX}={\max}_j\|Q_j^TA\|_{\cZ_j^*,\cX}={\max}_i\|AP_i\|_{\cB,\cX_i}.
$$
As we know from Theorem \ref{verynewRelaxationTheorem}, we can upper-bound $\|Q_j^TAP_i\|_{ij}$ by  $\Phi_{ij}(Q_j^TAP_i)$ with convex and efficiently computable function $\Phi_{ij}(\cdot)$, the bound being tight within the factor $\varsigma(K_i,L_j)\leq3\sqrt{\ln(4K_i)\ln(4L_j)}$, where $K_i$ and $L_j$ are the ellitopic sizes of $\cX_i$ and $\cZ_j$. As a result, the efficiently computable convex function
$$
\Phi(A)={\max}_{i,j}\Phi_{ij}(Q_j^TAP_i)
$$
is an upper bound on $\|A\|_{\cB,\cX}$ tight within the factor
$
3\sqrt{\ln(4\max_iK_i)\ln(4\max_jL_j)}.
$

In some simple situations the above tightness factor can be improved. For example, when $\cX_i=\{x_i:\|x_i\|_{p_i}\leq 1\}$, $\cZ_j=\{z_j:\|z_j\|_{q_i}\leq 1\}$ with $p_i\geq2$, $q_i\geq 2$, by  Nesterov's results of
(cf. Remark \ref{nesterov}) the tightness factor is an absolute constant (e.g., 1 in the trivial case where $p_i=q_j=2$ for all $i,j$).
%%%%%%%%%%%%%%%%%%%%%%%%%%%%%%%%%%%%%%%%%%%%%%%%%%%%%%%%%%%%%%%%%%%%%%%%%%%
\subsection{Applications}\label{ssappl}
\subsubsection{Least norm projector synthesis}
Consider the {\sl projection problem} as follows: we are given a linear {subspace} $\cF$ of linear space $\cE=\bR^n$ and a norm $\theta(\cdot)$ on $\cE$; our goal is to find a linear projector $H$ of $\cE$ onto $\cF$---a linear map $x\mapsto H x:\cE\to\cF$ with $Hx=x$ for all $x\in\cF$---which deviates the least from the identity mapping $\Id$ in the norm\\
\centerline{$
\|\cdot\|_{\theta\to\theta}: \|A\|_{\theta\to\theta}={\max}_{x\in\cE} \{\theta(Ax):\theta(x)\leq 1\}.
$}
Consider the case when the norm in question is the block $\ell_\infty/\ell_2$ norm\\
\centerline{$
\qquad\qquad\qquad\qquad\theta(x)=\max_{k\leq K} \|G_kx\|_2\hfill{[x\mapsto G_kx:{E_k}\to\bR^{\nu_k},\,\bigcap_k\Ker E_k=\{0\}]}
$}
What makes the projection problem potentially difficult is the {\sl block} $\ell_\infty$ structure of $\theta$; were $\nu_k=1$ for all $k$,  $\|\cdot\|_{\theta\to\theta}$ would have polyhedral epigraph, and minimization of $\|\Id-H\|_{\theta\to\theta}$ would be a Linear Programming problem (note that property of {$H$} to project onto $\cF$  {reduces to} a system of linear equalities on {$H$}).\footnote{{Allowing for a slight abuse of notation, we denote with $H$ the matrix of the linear mapping $H$.}} In contrast, in the general $\ell_\infty/\ell_2$ case as described above, the problem is NP hard. At the same time, the problem is within the scope of our machinery: the unit ball of $\theta$ is the ellitope
$$
\cX=\{x\in\bR^n: x^TG_k^TG_kx\leq1, k\leq K\},
$$
and therefore {$\theta$} is a simple ellitopic {norm}. At the same time, we have
$$
\theta(x)=\|Gx\|_{\infty/2},\;Gx=[G_1x;G_2x;...;G_Kx],\;\|[y_1;...;y_k]\|_{\infty/2}=\max_k\|y_k\|_2
$$
As we know, $\|\cdot\|_{\infty/2}$ is co-ellitopic (see {\bf cE.2} in Section \ref{extension}) and this property is preserved under restriction of a norm on a linear subspace ({\bf \cE.}3), and it remains to recall that $G$ is an embedding. The bottom line is that we can process the projection problem as explained in Section \ref{extension}. It is immediately seen that the {corresponding} recipe, under the circumstances, boils down to the following:
\begin{quote}
We select a linear basis $\{g_i:i\leq n\}$ in $\cE$ in such a way that the first $m=\dim\cF$ of these vectors form a basis of $\cF$; in the sequel, we identify vectors from $\cE$ with collections of their coordinates in this basis, and linear mappings from $\cE$ to $\cE$ with their matrices in this basis. Note that the (matrices of) projectors of $\cE$ onto $\cF$ are exactly block-matrices  $\left[\begin{array}{c|c}I_m&P\cr\hline
    &\cr\end{array}\right]$ with $m\times (n-m)$ blocks $P$.  Applying Theorem \ref{verynewRelaxationTheorem}, we arrive at the
efficiently solvable convex optimization problem
{\small\begin{equation}\label{proj123}
\begin{array}{l}
    \Opt=\min\limits_{P,\{\mu_k,\lambda^k:k\leq K\}}\Bigg\{
    \max_k\left[\mu_k+\sum_{j=1}^K\lambda^k_j\right]:\;
    P\in\bR^{m\times (n-m)},\\
    \lambda^k\geq0,
    \mbox{\scriptsize $\left[\begin{array}{c|c}\mu_kI_{\nu_k}&{1\over 2}G_k\left[\begin{array}{c|c}&P\cr\hline
    -I_{n-m}\cr\end{array}\right]\cr\hline
    {1\over 2}\left[\begin{array}{c|c}&\cr\hline
    P^T&-I_{n-m}\cr\end{array}\right]G_k^T&\sum_j\lambda_j^kG_j^TG_j\cr\end{array}\right]$}\succeq0,\,k\leq K\Bigg\}
    \end{array}
\end{equation}}\noindent   which is a safe tractable approximation of the problem of interest---the $P$-component of a feasible solution to the problem specifies projector of $\cE$ onto $\cF$ with the value of $\|\cdot\|_{\theta\to\theta}$
    not exceeding the value  of the objective at this solution. This approximation is tight within the factor $O(1)\sqrt{\ln(4K)}$, meaning that $\Opt$ is at most by this factor greater than the actual optimal value in the projection problem. In addition, when $\nu_k=1$\,\,$\forall k$, the tightness factor is exactly 1.
    \end{quote}
\def\lrg{{\hbox{\tiny L}}}
\def\sml{{\hbox{\tiny S}}}
\subsubsection{Illustration: projecting splines}
Consider a partition of $[0,1]$ into $M$ ``large'' segments, which are further partitioned into total of $N$ ``small'' segments.
Let also $\Gamma$ be equidistant grid on $[0,1]$ with $L$ points.
{Given nonnegative integers $\mu_\lrg\geq \nu_\lrg$, $\mu_\sml\geq \mu_\lrg$, and $\nu_\sml\leq \nu_\lrg$, let} us define ${\cF}$ as the linear space of restrictions on $\Gamma$ of splines which are polynomials of degree at most ${\mu_\lrg}$ in every large segment, with all derivatives of order $\leq {\nu_\lrg}$ continuous on the entire $[0,1]$. We define $\cE$ as the linear space of restrictions on $\Gamma$ of splines which are polynomials of degree of order $\leq \mu_\sml$ in every small segment and have continuous on $[0,1]$ derivatives of order $\leq \nu_\sml$. With the above inequalities between $\mu$'s and $\nu$'s, $\cF$ is a subspace in $\cE$. Now let $\Delta_1,...,\Delta_K$ be partitioning of $\Gamma$ into $K$ consecutive segments, and let $\theta$ be the $\ell_\infty/\ell_2$ norm on $\cE$ given by
$$
\theta(x)={\max}_{k\leq K}\sqrt{{\sum}_{i\in\Delta_k}x_{i}^2},
$$
$x_i$ being the value of spline $x\in\cE$ at the $i$-th point of $\Gamma$.
\par
\begin{figure}[htb!]
\begin{center}{\small
\begin{tabular}{cc}
\includegraphics[width=0.35\textwidth]{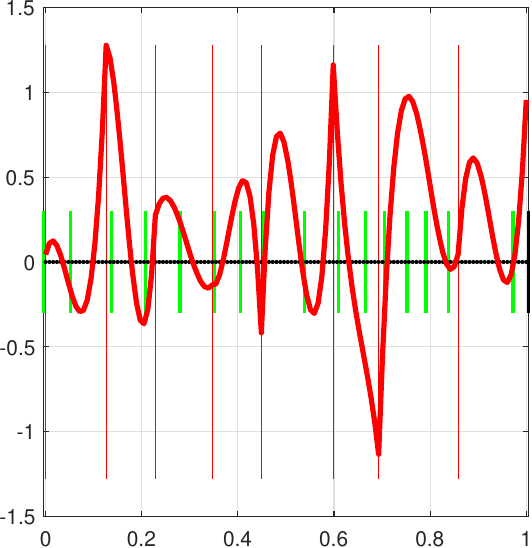}&\includegraphics[width=0.35\textwidth]{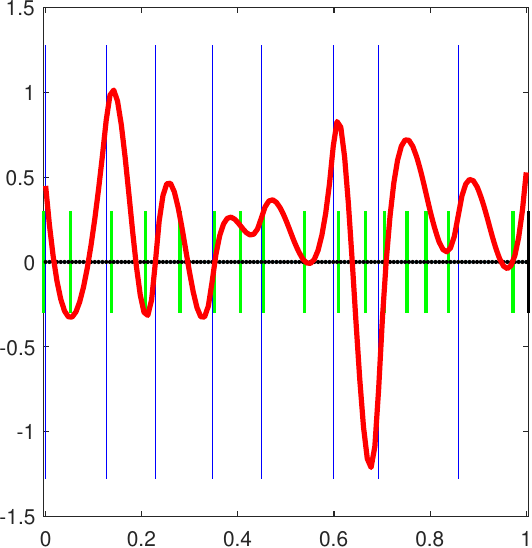}\\
\end{tabular}
\caption{\label{fig:1}  Spline $x$ from $\cE$ (left plot) and its projection $Hx$ on $\cF$ (right plot).
}}\end{center}
\end{figure}

 In Figure \ref{fig:1} we present a sample pair of a spline from $\cE$ and its projection onto $\cF$.
In this experiment, $|\Gamma|=128$, {there are eight identical large and small segments (separated by red/blue
vertical lines on the plots),} and $K=16$  (on the plots, 16 segments $\Delta_k$ are separated from each other by green
vertical lines). Splines from $\cE$ are continuous on $[0,1]$ and are polynomials of degree 3 on large/small segments, and $\cF$ is cut off $\cE$ by additional requirement for the spline to be continuously differentiable on $[0,1]$. {Solving \rf{proj123}} yields  $H$ with $\|\Id-H\|_{\theta\to\theta}\leq \Opt\approx1.255$  and $H$ is ``essentially different'' from the $\|\cdot\|_2$-orthogonal
projection\footnote{Recall that vectors from $\cE$ are restrictions of functions on $[0,1]$ onto equidistant grid in this segment and as such {$\cE$} is equipped with ``canonical'' Euclidean structure.} $\overline{H}$ of $\cE$ onto $\cF$---the spectral norm of $H-\overline{H}$ is $\approx 0.69$, and the upper bound on $\|\Id-\overline{H}\|_{\theta\to\theta}$, as given by our machinery, is $\approx1.527$. In fact both upper bounds $\approx 1.255$ on $\|\Id-H\|_{\theta\to\theta}$ and
$\approx1.527$ on $\|\Id-\overline{H}\|_{\theta\to\theta}$ happen to coincide within {four significant digits} with the quantities themselves.\footnote{One can easily build a numerical lower bound on $\|A\|_{\cB,\cX}$ by alternating maximization of the bilinear function $u^TAx$ over $u\in\cB_*$ and $x\in\cX$; in the reported experiment, these lower bounds were {within the indicated accuracy with} the upper bounds yielded by our machinery.}
\subsubsection{Synthesis of linear controller with peak-to-peak design specifications}\label{peak-to-peak} The situation we are about to address is as follows. We control a discrete time linear system
\bse
x_0=z,\,\,
x_{t+1}=A_tx_t+B_tu_t+D_td_t,0\leq t<N,\,\,
y_t=C_tx_t+E_td_t
\ese
where $x_t\in\bR^{n_x}$, $u_t\in\bR^{n_u}$, $d_t\in\bR^{n_d}$, and $y_t\in\bR^{n_y}$ are, respectively, states, controls, external disturbances, and observable outputs. When augmented with non-anticipating affine controller\\
\centerline{$
u_t=g_t +{\sum}_{\tau=0}^tG^t_\tau y_{t-\tau}
$}
the closed loop system specifies affine mappings
$$
\begin{array}{c}
(d^N:=[d_0;d_1;...;d_{N-1}],z)\mapsto  x^N:=[x_1;...;x_N]=X^Nd^N+\overline{X}^Nz+\widehat{X}^N\\
(d^N,z)\mapsto  u^N:=[u_0;...;u_{N-1}]=U^Nd^N+\overline{U}^Nz+\widehat{U}^N,\\
(d^N,z)\mapsto  y^N:=[y_0;...;y_{N-1}]=Y^Nd^N+\overline{Y}^Nz+\widehat{Y}^N\\
\end{array}
$$
{With} ``smart parameterizations'' of the controller---
passing from $\{g_t, G^t_\tau,0\leq t<N,0\leq\tau\leq t\}$ to the parameters of the {\em affine purified-output-based controller,} matrices ${X}^N$,... ,$\widehat{Y}^N$  become {\sl affine} functions of the vector $\chi$ of controller's  parameters; this vector runs through certain finite-dimensional linear space $\cC$ equipped with filtration $\cC_0\subset\cC_1\subset...\subset\cC_{N-1}=\cC$ by linear subspaces, with $\cC_d$ comprised of ``controllers with memory $d$.'' {We refer the reader to \cite{KLN} for details of the controller construction.}
\par
When designing a controller, one of natural design specifications (traditionally considered as not so easy to handle, cf., e.g., \cite{abedor1996linear,aubrecht2001minimization,balakrishnan1992computing,boyd1987comparison,diaz1993minimization} and reference therein) are bounds on ``peak-to-peak'' gains. The {\sl disturbance-to-state} gain is nothing but the norm of the matrix ${X}^N$ induced by the norm \[
\|d^N\|_{\infty/p}={\max}_{0\leq t<N}\|d_t\|_p
\] on the space of sequences $d^N$ of disturbances and the norm
\[\|x^N\|_{\infty/r}={\max}_{1\leq t\leq N} \|x_t\|_r
\] on the  space of sequences of states; disturbance-to-control and disturbance-to-output peak-to-peak gains are defined similarly.
When $\infty\geq p\geq2$ and $1\leq r\leq2$, we can {enforce} the desired bound on the peak-to-peak gain (which can be difficult to handle,
 since the corresponding norm of ${X}^N$ is, in general, difficult to compute) by bounding from above the efficiently computable upper bound, yielded by our machinery,
 on the gain. As a result, we get an efficiently tractable convex constraint on the parameters  of the controller which safely (and tightly within the factor $\sqrt{\pi/2}$, see the concluding comments in Section \ref{extension})
 approximates the design specification in question.
 \par
 {Note that our machinery remains applicable when $\|\cdot\|_p$ and $\|\cdot\|_r$ are replaced with, respectively, an s.e.n. $\|\cdot\|_{(d)}$ and a s.co-e.n. $\|\cdot\|_{(x)}$,
 same as when the design specifications impose bound on the ``restricted'' peak-to-peak gains, e.g., on the peak-to-peak disturbance-to-state gain when the disturbances
 $d^N$ are restricted to reside in a given linear subspace of the ```complete disturbance space'' $\bR^{n_dN}$.}
 \par
\begin{figure}[htb!]
\begin{center}{\small
\begin{tabular}{cc}
\includegraphics[width=0.45\textwidth]{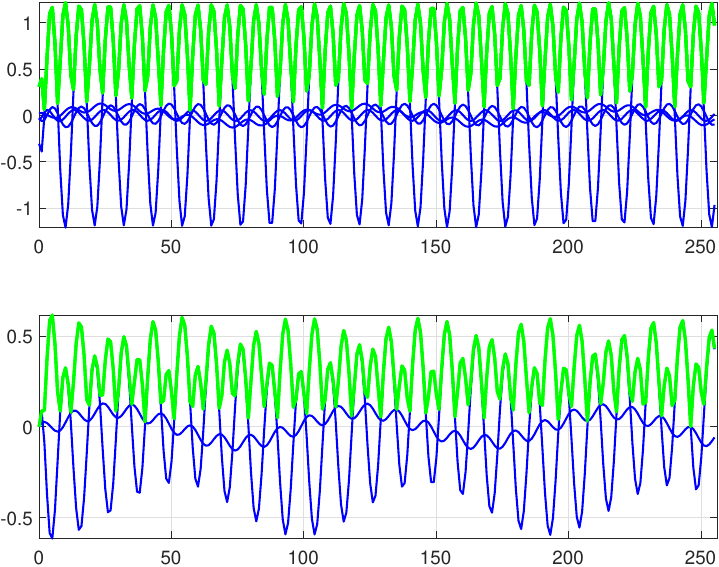}&\includegraphics[width=0.45\textwidth]{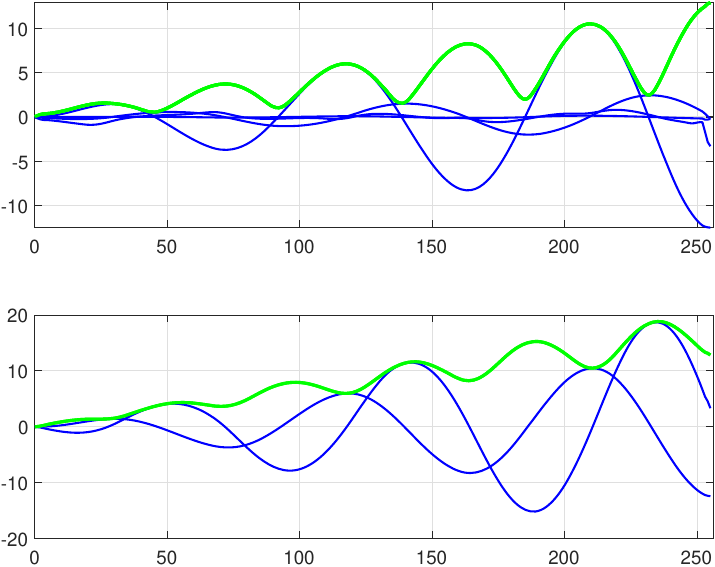}\\
\multicolumn{2}{c}{Zero control}\\
\includegraphics[width=0.45\textwidth]{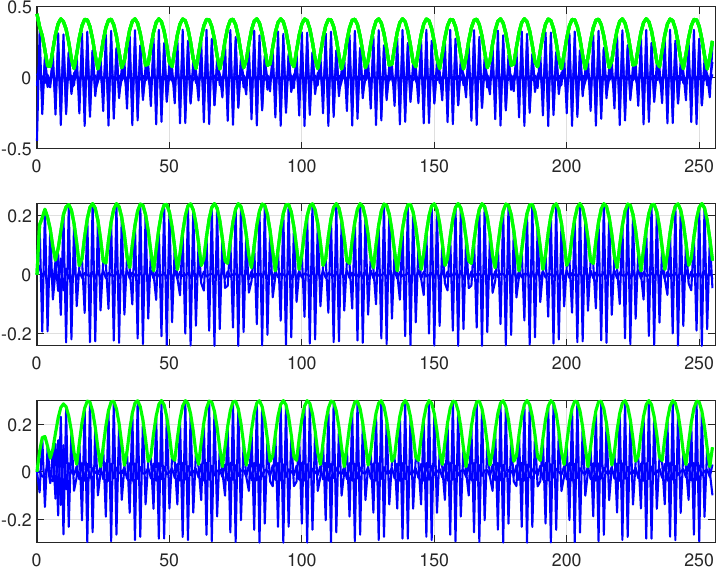}&\includegraphics[width=0.45\textwidth]{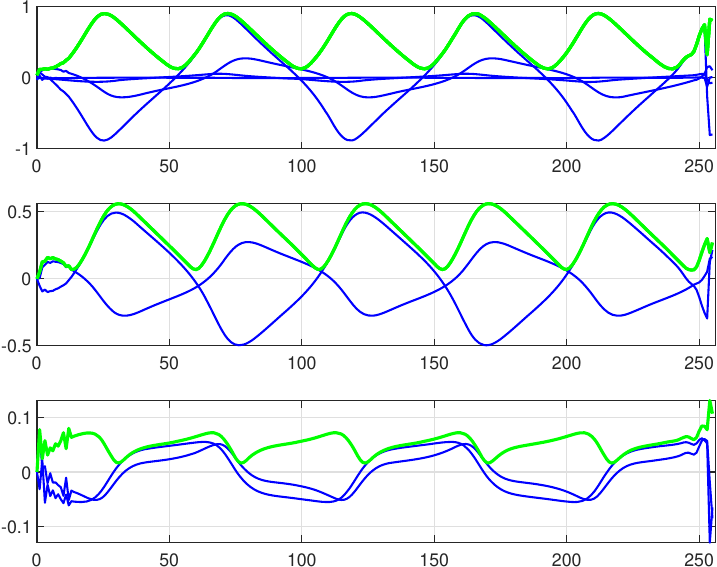}\\
\multicolumn{2}{c}{Synthesized control}
\end{tabular}
\caption{\label{fig:2}  In blue,  from top to bottom: state ($n_x=4$), output ($n_y=2$) and control ($n_u=2$, on the synthesized control plots) trajectories of the controlled plant. In the left pane: random harmonic oscillation disturbance, in the right pane: ``bad disturbance.'' In green: $\|\cdot\|_2$-norms of states, outputs and controls, respectively.}
}\end{center}
\end{figure}

{\sl Numerical illustration} we are about to present deals with minimizing disturbance-to-state $\infty/2$ peak-to-peak gain (i.e., $p=r=2$) when controlling
linearized and discretized in time motion of Boeing 747; the model we use originates from \cite{boyd_lecture}, see also Section \ref{roblinrec} below. We omit irrelevant for our purposes details (which can be found in \cite{KLN}), here it suffices to mention that the model is time-invariant (matrices $A_t\equiv A$,...,$E_t\equiv E$) with $n_x=4$ and $n_u=n_d=n_y=2$. Applying our machinery on time horizon $N=256$ to build a purified-output-based linear controller with memory depth (whatever it means) 16, we end up with controller with disturbance-to-state peak-to-peak gain $\approx 1.02$. To put this result into proper perspective, note that the matrix $A$ of the model in question is {only marginally stable} (the corresponding spectral radius is 0.9995). As a result, although trivial---identically zero---control results in uniformly bounded in $N$ peak-to-peak gain, this gain ($\approx 12$) is more than 10 times larger than the gain of the computed controller. Sample trajectories of the system with and without control are presented in Figure \ref{fig:2}.
{In the reported experiments, $\|d_t\|_2\equiv 1$ for all $t$.} ``Bad'' disturbance is selected to result in large peak-to-peak gain with vanishing control; in this case $\max_t\|x_t\|_2$ turns to be $\approx 12$, while with the control yielded by our synthesis, the same disturbances result in  $\max_t\|x_t\|_2\approx 0.9$, which is close to the upper bound on the gain ($\approx 1.02$) guaranteed  by our synthesis.
\section{Bounding robust norms of uncertain matrices}\label{secrobnorms}
\subsection{Motivation}
Consider the following problem which arises, e.g., in Robust Control:
\begin{quote}
{\sl Given box-type uncertainty set
$$
{\cal A}[\rho]=\{A={\sum}_{s=1}^Sz_s A_s: \|z\|_\infty\leq \rho\}
$$
in  the space of $m\times n$ matrices, upper-bound the quantity
$$
\Opt_*(\rho)=\max_{A\in{\cal A}[\rho]}|A|,
$$
where $|\cdot|$ stands for the spectral norm of a matrix.}
\end{quote}
This problem can be immediately reduced to the Matrix Cube problem (cf. \cite{BTNMC}, see also \cite[Section 3.4.3.1]{LMCO}): associating with $m\times n$ matrix $A$ symmetric $(m+n)\times(m+n)$ matrix
$$
\cL[A]=\left[\begin{array}{c|c}&{\half} A\cr\hline
{\half} A^T&\cr\end{array}\right],
$$
we have $|A|\leq R$ if and only if $RI_{m+n}-{2}\cL[A]\succeq0$. Therefore, the inequality
\begin{equation}\label{May2021eq1}
\Opt_*(\rho)\leq R
\end{equation}
is equivalent to
$$
RI_{m+n}+{2}{\sum}_{s=1}^Sz_s\cL[A_s]\succeq0\;\;\forall (z:\|z\|_\infty\leq\rho).
$$
According to the results of \cite{BTNMC}, reproduced in \cite[Theorem 3.4.7]{LMCO}, an  efficiently verifiable sufficient condition for the validity of the latter semi-infinite Linear Matrix Inequality (LMI)  is the solvability of the parametric system of LMIs
$$
RI_{m+n}-\rho{\sum}_{s=1}^SU_s{\succeq0},\;U_s\succeq \pm {2}\cL[A_s],\;1\leq i\leq N,\leqno{\cR[R,\rho]}
$$
in matrix variables $U_s$, and this sufficient condition is  tight within factor $\vartheta(\mu)$ depending solely of the maximum of ranks $2\rank(A_s)$ of the ``edge matrices'' $\cL[{A_s}]$.  Specifically,
\begin{quote} {\sl setting $\mu=\max\limits_{1\leq s\leq S} \rank(A_s)$, we obtain:
\begin{itemize}
\item {\rm (\ref{May2021eq1})} does take place when $\cR[R,\rho]$ is feasible, and
\item when  $\cR[R,\rho]$ is infeasible, one has $\Opt(\vartheta(2\mu)\rho)>R$,
where  $\vartheta(\cdot)$ is a universal function specified in {\rm \cite{BTNMC}} (cf. {\rm\cite[display (3.4.39)]{LMCO}} and (\ref{theta}) below) such that
\end{itemize}}
\end{quote}
\vspace{-0.5cm}
\begin{equation}\label{4.qwer8}
\vartheta(1)=1,\vartheta(2)=\pi/2,\vartheta(3)=1.7348...,\vartheta(4)=2\ \&\
\vartheta(k)\leq\half {\pi\sqrt{k}},\,\,k\geq1.
\end{equation}
The goal of this section is to extend this result onto more general matrix norms considered in Section \ref{mainresellnorms}.
\subsection{Problem setting and main result}\label{uncerftainmatr}
Let ellitopes  $\X\subset\bR^n$, $\B_*\subset\bR^m$ with nonempty interior and  basic ellitopes $\W$, $\cZ$ be given by (\ref{2020ell2*}), let $\cB$ be the polar of $\B_*$,
and let $A_s\in\bR^{m\times n}$, $1\leq s\leq S$. These data define {\sl the uncertain matrix with box uncertainty}\\
\centerline{$
\cA=\left\{A=\sum_s\epsilon_s A_s:\|\epsilon\|_\infty\leq1\right\}\subset\bR^{m\times n}
$}
and the quantity\\
\centerline{$
\|\cA\|_{\cB,\cX}={\max}_{A\in\cA}\|A\|_{\cB,\cX}
$}
which we refer to as {\sl robust $\|\cdot\|_{\cB,\cX}$-norm} of uncertain  matrix $\cA$. Note that this norm is difficult to compute already in the case of ``general position'' symmetric matrices $A_s$ of rank 2. Our goal is to conceive  a computationally efficient upper-bounding of the robust norm.
\par Let us consider the quantity
\be
\varkappa(J)=\left\{\begin{array}{ll}1,&J=1,\\
{5\over 2}\sqrt{\ln(2J)},&J>1,
\end{array}\right.
\ee{May2021upsilon-1}
and function $\vartheta$ of the positive integer argument
{\small\be
\vartheta(k)=\left[\min_{\alpha}\left\{(2\pi)^{-k/2}\int
|\alpha_1u_1^2+...+\alpha_ku_k^2|\e^{-{u^Tu/2}}du,\;\alpha\in{\bR}^k,\|\alpha\|_1=1\right\}\right]^{-1};
\ee{theta}}\noindent
note that $\vartheta(k)$ satisfies \rf{4.qwer8} \cite{BTNMC}.
Let also
{\small \be
\Opt=\min\limits_{\lambda\geq0,\upsilon\geq0,\atop\{G_s,H_s\}}\left\{\phi_{\T}(\lambda)+\phi_{\R}(\upsilon):\!\!\begin{array}{l}\left[\begin{array}{c|c}G_s&{1\over2}Q^TA_sP\cr\hline{1\over2} P^TA_s^TQ&H_s\cr\end{array}\right]\succeq0,s\leq S\\
\sum_sG_s\preceq\sum_\ell\upsilon_\ell R_\ell,
\sum_sH_s\preceq \sum_k \lambda_k T_k.
\end{array}\right\}.
\ee{May2021eq2-1}}\noindent
{\begin{proposition}\label{propmatrcube}
In the situation {of this section,} assuming that ranks of all $A_s$ are $\leq \kappa$, the efficiently computable quantity $\Opt$ as given by {\rm\rf{May2021eq2-1}} is a reasonably tight upper bound
on the robust norm $\|\cA\|_{\cB,\cX}$ of uncertain matrix $\cA$, specifically,
\be
\|A\|_{\cB,\cX}\leq\Opt\leq \varkappa(K)\varkappa(L)\vartheta(2\kappa)\|A\|_{\cB,\cX}
\ee{wehavethat}
{where $K$ and $L$ are given by {\rm \rf{2020ell2*}}.}
\end{proposition}}\noindent
%%%%%%%%%%%%%%%%%%%%%%%%%%%%%%%%%%%%%%%%%%%%%%
%%%%%%%%%%%%%%%%%%%%%%%%%%%%%%%%%%%%%%%%%%%%%%%%%%%%%%%%%%%%%%
\begin{remark}{\rm Assume that matrices $A_s=A_s[\chi]$ are affine in some vector $\chi$ of control parameters. In this case, the quantities $\|\cA\|_{\cB,\cX}$ and its efficiently computable upper bound $\Opt$ become functions  $\Opt_*(\chi)=\|\cA\|_{\cB,\cX}$ and $\Opt(\chi)$ of $\chi$, and it is immediately seen that both functions are convex. As a result, we can handle, to some extent, the problem of minimizing over $\chi$  the robust $\|\cdot\|$-norm of uncertain matrix
$$
\cA[\chi]=\left\{A={\sum}_s\epsilon_sA_s[\chi]:\|\epsilon\|_\infty\leq1\right\}.
$$
More precisely, we can minimize over $\chi$ efficiently computable convex upper bound $\Opt(\chi)$ on the robust norm $\Opt_*(\chi)$ of $\cA[\chi]$, the bound being reasonably tight provided that the ranks of matrices
$A_s[\chi]$ are small for all $\chi$ in question.}
\end{remark}
\begin{remark}\label{rem32}{\rm Note that the quantity
$$
\Opt=\min\limits_{\lambda\geq0,\upsilon\geq0,\atop\{G_s,H_s\}}\left\{\phi_{\T}(\lambda)+\phi_{\R}(\upsilon):\begin{array}{l}
\sum_sG_s\preceq\sum_\ell\upsilon_\ell R_\ell,\,
\sum_sH_s\preceq \sum_k \lambda_k T_k\\\left[\begin{array}{c|c}G_s&{1\over2}Q^TA_sP\cr\hline{1\over2} P^TA_s^TQ&H_s\cr\end{array}\right]\succeq0,s\leq S
\end{array}\right\}
$$ as given by \rf{May2021eq2-1} admits another  representation which may sometimes be more convenient. %  (cf. Section \ref{systident}).
Specifically, excluding trivial case $\Opt=0$ which takes place if and only if $Q^TA_sP=0$ for all $s$, one has
\begin{equation}\label{onehas}
{1\over \Opt}=  \max\limits_{\rho,\{G_s,H_s\},\lambda,\upsilon}\left\{\rho
:\begin{array}{l}
\lambda\geq0,\upsilon\geq 0,\phi_\cT(\lambda)\leq 1,\phi_\cR(\upsilon)\leq 1\\
\sum_sG_s\preceq\sum_\ell\upsilon_\ell R_\ell,
\sum_sH_s\preceq \sum_k \lambda_k T_k\\
\left[\begin{array}{c|c}G_s&\rho Q^TA_sP\cr\hline\rho P^TA_s^TQ&H_s\cr\end{array}\right]\succeq0,s\leq S
\end{array}
\right\}.
\end{equation}
}
\end{remark}
{\small \noindent Indeed, the optimization problem specifying $\Opt$ clearly is solvable; let $\lambda_,\upsilon,\{G_s,H_s\}$ be its optimal solution. Looking at the problem, we see, first, that
$\Opt>0$ implies $\lambda\neq0$ and $\upsilon\neq0$, and thus $\phi_\cR(\upsilon)>0$ and $\phi_\T(\lambda)>0$. Furthermore, whenever $\theta>0$, the collection
$\theta^{-1}\lambda,\theta\upsilon,\{\theta G_s,\theta^{-1}H_s\}$ is a feasible solution with the value of the objective $\theta\phi_\cR(\upsilon)+\theta^{-1}\phi_\T(\lambda)$. Since the solution we have started with is optimal, we have
$$\theta\phi_\cR(\upsilon)+\theta^{-1}\phi_\T(\lambda)\geq \phi_\cR(\upsilon)+\phi_\T(\lambda)=\Opt.
$$
This inequality holds true for all $\theta>0$, which with positive $\phi_\cR(\lambda)$ and $\phi_\T(\lambda)$ is possible if and only if $\phi_\cR(\upsilon)=\phi_\T(\lambda)=\Opt/2$. It follows that setting
$$
\overline{\lambda}=2\lambda/\Opt,\; \overline{\upsilon}=2\upsilon/\Opt,\; \overline{G}_s=2G_s/\Opt,\; \overline{H}_s=2H_s/\Opt,\; \rho=1/\Opt,
$$
we get a feasible solution to (\ref{onehas}) with the value of the objective $1/\Opt$, implying that the left hand side in (\ref{onehas}) is $\leq$ the right hand side.
On the other hand, the optimization problem in (\ref{onehas}) clearly is solvable. If $\rho,\lambda,\upsilon,\{G_s,H_s\}$ is an optimal solution {to \rf{onehas}}
then $\overline{G}_s=G_s/(2\rho)$, $\overline{H}_s=H_s/(2\rho)$, $\overline{\lambda}=\lambda/(2\rho)$, $\overline{\upsilon}=\upsilon/(2\rho)$ clearly form a feasible solution to the problem specifying $\Opt$, and the value of the objective of the latter problem
at this solution is $\leq 1/\rho$. Thus, $\Opt\leq 1/\rho$, $\rho$ being the optimal value of the optimization problem in (\ref{onehas}), so that the left hand side in (\ref{onehas}) is $\geq$ the right hand side.}

\subsubsection{An extension}\label{moreextension} Similarly to what was done in Section \ref{extension}, the above results can be straightforwardly extended to the case when $\|\cdot\|_\cX$ is ellitopic, and $\|\cdot\|_\cB$ is co-ellitopic norm.
Specifically,  for an uncertain matrix
$$\cA=\left\{{\sum}_s\epsilon_s A_s:\|\epsilon\|_\infty\leq1\right\}$$
the robust norm of $\cA$ in the case of (\ref{A}), (\ref{B}) is
$$
\max\limits_{i,j}\Big\|\Big\{{\sum}_s \epsilon_s Q_j^TA_s P_i:\|\epsilon\|_\infty\leq 1\Big\}\Big\|_{\cZ^*_j,\cX_i}
$$
where $\cX_i$ and polars $\cB_j$ of $\cZ^*_j$ are ellitopes, and we know how to efficiently upper-bound the robust norms
$\|\{{\sum}_s \epsilon_s Q_j^TA_s P_j:\|\epsilon\|_\infty\leq 1\}\|_{\cZ^*_j,\cX_i}$ and how tight such bounds are.
%%%%%
\subsubsection{Putting things together} So far, we have considered separately computationally efficient bounding of operator norms of matrices and  robust norms of uncertain matrices with box uncertainty. In applications to follow, we will be interested in a ``mixed'' setting, where we want to upper-bound the robust norm
$$
\|\cU\|_{\cB,\cX}={\max}_{A\in\cU}\|A\|_{\cB,\cX}
$$
of uncertain matrix
\begin{equation}\label{cU}
\cU=A_\n+\cA,\;\;\cA=\left\{{\sum}_{s=1}^S\epsilon_sA_s:\|\epsilon\|_\infty\leq1\right\}.
\end{equation}
The corresponding blend of our preceding results is as follows:
{\begin{proposition}\label{propmixture} Let $\cX\subset\bR^n$, $\cB,\cB_*\subset\bR^m$ be given by {\rm (\ref{A}), (\ref{B})}, with basic ellitopes
$$
\begin{array}{rcl}
\cX_i&=&\left\{x_i\in\bR^{\nu_i}:\exists t^i\in\cT^i: x_i^TT_{ki}x_i\leq t^i_k,1\leq k\leq K_i\right\},\,i\leq I\\
\cZ_j&=&\left\{z_j\in\bR^{\mu_j}:\exists s^j\in\cR^j: z_j^TR_{\ell j}z_j\leq s^j_\ell,1\leq \ell\leq L_,\right\},\,j\leq J\\
\end{array}
$$
Then
the efficiently computable quantity
\[
\Opt[\cU]=
{\max}_{i\leq I,j\leq J}\Opt_{ij}[\cU],
\]
where
{\small\be
\Opt_{ij}[\cU]&=&\min\limits_{\lambda^{ij},\upsilon^{ij},G^{ijs},H^{ijs}\atop
{\overline{G}^{ij},\overline{H}^{ij}\atop
1\leq i\leq I,1\leq j\leq J,1\leq s\leq S}}
\bigg\{\phi_{\cT^i}(\lambda^{ij})+\anc{\psi}{\phi}_{\cR^j}(\upsilon^{ij}):\;\;\lambda^{ij}\geq0,\,\upsilon^{ij}\geq0\nn
&&\qquad\quad
\left.\begin{array}{l}
{\sum}_{s=1}^SG^{ijs}+\overline{G}^{ij}\preceq {\sum}_{\ell=1}^{L_j}\upsilon^{ij}_\ell R_{\ell j},{\sum}_{s=1}^SH^{ijs}+\overline{H}^{ij}\preceq {\sum}_{k=1}^{K_i}\lambda^{ij}_kT_{ki}\\
\left[\begin{array}{c|c}G^{ijs}&{1\over 2}[Q_j^TA_sP_i]\cr\hline
{1\over 2}[Q_j^TA_sP_i]^T&H^{ijs}\cr\end{array}\right]\succeq0,s\leq S\\ \left[\begin{array}{c|c}\overline{G}^{ij}&{1\over 2}[Q_j^TA_\n P_i]\cr\hline
{1\over 2}[Q_j^TA_\n P_i]^T&\overline{H}^{ij}\cr\end{array}\right]\succeq0\end{array}
\right\}
\ee{eqmixture1}}\noindent
is an efficiently computable {\em convex} in $(A_\n,A_1,...,A_S)$ upper bound on $\|\cU\|_{\cB,\cX}$. This upper bound is reasonably tight,
specifically, setting
$$
\cU_{ij}=Q_j^TA_\n P_i+\Big\{{\sum}_{s=1}^S\epsilon_s[Q_j^TA_sP_i]:\|\epsilon\|_\infty\leq1\Big\},
$$
we have
\[\|\cU_{ij}\|_{\cZ_j^*,\cX_i}\leq\Opt_{ij}[\cU]\leq [\varsigma(K_i,L_j)+{\varkappa}(K_i){\varkappa}(L_j)\vartheta(2\kappa)] \|\cU_{ij}\|_{\cZ_j^*,\cX_i},
\]
and
\be
\|\cU\|_{\cB,\cX}&=&\max\limits_{i\leq I, j\leq J}\|\cU_{ij}\|_{\cZ_j^*,\cX_i}\leq \Opt[\cU]=\max\limits_{i\leq I, j\leq J}\Opt_{ij}[\cU]\nn
&\leq&
\left[\max\limits_{i\leq I, j\leq J}[\varsigma(K_i,L_j)+{\varkappa}(K_i){\varkappa}(L_j)\vartheta(2\kappa)]\right]\|\cU\|_{\cB,\cX}
\ee{eqmixture2}
where $\kappa$ is the maximum of ranks of $A_s$, $1\leq s\leq S$, and $\varsigma(K,L)$ and $\varkappa(\cdot),\vartheta(\cdot)$ are as defined in
Theorem \ref{verynewRelaxationTheorem} and Proposition \ref{propmatrcube}.
\end{proposition}}\noindent
{Note that ``extreme cases'' ($A_s=0$ for all $s$, on one hand, and $A_\n=0$, on the other) of Proposition \ref{propmixture} recover Theorem \ref{verynewRelaxationTheorem} and Proposition \ref{propmatrcube}, and even  their ``advanced'' versions with simple ellitopic/co-ellitopic norms extended to ellitopic/co-ellitopic ones.}
%%%%%%%%%%%%%%
\subsection{Application to robust signal recovery}
Consider the standard Signal Processing problem as follows. {\sl Given noisy observations
\begin{equation}\label{observations}
\omega = Ax+\xi,\quad\xi\sim\cN(0,I_m)
\end{equation}
of unknown signal $x$ known to belong to a given signal set $\cX\subset\bR^n$, we want to recover  $Bx\in\bR^\nu$.} Here $A\in\bR^{m\times n}$ and $B\in\bR^{\nu\times n}$ are
given matrices.
We consider {\sl linear recovery} $\widehat{x}=\widehat{x}_H(\omega):=H^T\omega$, $H\in \bR^{\nu\times m}$ and quantify the performance of a candidate estimate $\widehat{x}_H$  by its {\sl worst-case risk}
\def\Risk{{\mathrm{Risk}}}
$$
\Risk_{\|\cdot\|_\cB}[\widehat{x}_H|\cX]={\sup}_{x\in\cX}\bE_{\xi\sim\cN(0,I_m)}\left\{\|Bx-\widehat{x}_H(Ax+\xi)\|_\cB\right\},
$$
where $\|\cdot\|_\cB$ is a given norm on $\bR^\nu$. There is an extensive literature dealing with the design and performance analysis of linear estimates. In particular, it is known \cite[Proposition 4.\aic{8.1}{16}]{STOPT} that when $\cX$ is an ellitope  of ellitopic size $K$ and the polar $\cB_*$ of the unit ball $\cB$ of $\|\cdot\|$ is an ellitope of ellitopic size  $L$, the linear estimate $\widehat{x}_{H_*}$ yielded by the optimal solution to an explicit efficiently solvable convex optimization problem is optimal within logarithmic in $K$, $L$ factor:
$$
\Risk_{\|\cdot\|_\cB}[\widehat{x}_{H_*}|\cX]\leq O(1)\sqrt{\ln(2K)\ln(2L)} \Risk\Opt_{\|\cdot\|_\cB}[\cX];
$$
here  $\Risk\Opt_{\|\cdot\|_\cB}[\cX]$ is the minimax risk---the infimum of risks $\Risk_{\|\cdot\|_\cB}[\widehat{x}|\cX]$ over all estimates $\widehat{x}$, linear and nonlinear alike.
\par
The result we have just cited, as well as most of known to us results on performance of linear estimates, deals with the case when the sensing matrix $A$ is known in advance. Here we want to address the case when $A$ is subject to ``uncertain-but-bounded'' perturbations, specifically, is selected (by nature or by an adversary) from the uncertainty set
$${\cal U}=A_\n+\cA,\;\;\cA=\left\{{\sum}_{s=1}^S\epsilon_sA_s:\|\epsilon\|_\infty\leq 1\right\}.$$
{This problem can be seen as a ``non-interval'' extension of the problem of solving systems of equations affected by interval uncertainty which has received significant attention in the literature, cf., e.g., \cite{cope1979bounds,higham2002accuracy,kreinovich1993optimal,nazin2005interval,neumaier1990interval,oettli1964compatibility,polyak2003robust} and references therein.}
Assuming that given perturbation in $A$ and ``true'' signal $x$, the observation noise $\xi$  is $\cN(0,I_m)$, the worst-case risk of a linear estimate $\widehat{x}_H$ becomes
{\small\begin{align*}
\Risk^+_{\|\cdot\|_\cB}[\widehat{x}_H|\cX]&:=\sup\limits_{x\in \cX,\epsilon:\atop\|\epsilon\|_\infty\leq1}\bE_{\xi\sim\cN(0,I_m)}\left\{\Big\|[B-H^TA_\n]x-\Big[\sum\limits_s\epsilon_sH^TA_s\Big]x-H^T\xi\Big\|_\cB\right\}\\
&\leq
\|{\cal V}[H]\|_{\cB,\cX}+\bE_{\xi\sim\cN(0,I_m)}\{\|H^T\xi\|_\cB\}
\end{align*}}\noindent
where
\[
\cV[H]=\left\{[B-H^TA_\n]+{\sum}_{s=1}^S\epsilon_sH^TA_s:\|\epsilon\|_\infty\leq1\right\}
\]
and $\|{\cal V}[H]\|_{\cB,\cX}={\max}_{V\in{\cal V}[H]}\|V\|_{\cB,\cX}$.
The simplest way to build a ``presumably good'' linear estimate is to minimize over $H$ the sum of the (efficiently computable upper bound on the) robust norm
of ${\cal V}[H]$ and an efficiently computable upper bound on $\Psi(H):=\bE_{\xi\sim\cN(0,I_m)}\{\|H^T\xi\|_\cB\}$. Combining the results of Proposition \ref{propmixture} with the upper bound on $\Psi(H)$ from \cite[Lemma 4.\aic{8.6}{11}]{STOPT}, {in the case of $\cX=P\cX_1$, $\cB_*=Q\cZ_1$, $I=J=1$,} we  arrive at the efficiently solvable convex optimization problem
\be\begin{array}{rl}
~~~\Opt&=\min\limits_{H,\lambda,\upsilon,\mu,\atop
\Theta,G^s,H^s,G,H}\bigg\{\phi_\cT(\lambda)+\phi_\cR(\upsilon)+\phi_\cR(\mu)+\Tr(\Theta):\\
&\hbox{\small$\left.\begin{array}{l}\lambda\geq0,\,\upsilon\geq0,\mu\geq0,\,G+\sum_sG^s\preceq\sum_\ell \upsilon_\ell R_\ell,\,H+\sum_s H^s\preceq \sum_k \lambda_kT_k\\
\left[\begin{array}{c|c}\Theta&{1\over 2}HQ\cr\hline {1\over 2}Q^TH^T&\sum_\ell \mu_\ell R_\ell\cr\end{array}\right]\succeq0,\;\;
\left[\begin{array}{c|c}G^s&{1\over 2}Q^TH^TA_sP\cr\hline {1\over 2}P^TA_s^THQ& H^s\cr\end{array}\right]\succeq0,\,s\leq S\\
\left[\begin{array}{c|c}G&{1\over 2}Q^T[B-H^TA_\n]P\cr\hline {1\over 2}P^T[B^T-A_\n^TH]Q& H\cr\end{array}\right]\succeq0
\end{array}\right\}$}\end{array}
\ee{(!)}
(we use the notation from Proposition \ref{propmixture}\aic{, where $\cX=P\cX_1$, $\cB_*=Q\cZ_1$, $I=J=1$, allowing to skip indices $i,j$,}{} with $\nu$ in the role of $m$).
For every feasible solution to this problem, the value of the objective at the solution is an upper bound on $\Risk^+_{\|\cdot\|_\cB}[\widehat{x}_H|\cX]$, $H$ being the $H$-component of the solution in question.
Moreover, from Proposition \ref{propmixture} combined with \cite[Lemma 4.\aic{8.6}{11}]{STOPT} it follows that the function $\Opt[H]$ obtained by partial minimization of the objective in \rf{(!)} over all decision variables except $H$ is a tight, within factor $O(1)\sqrt{\ln(2K)\ln(2L)}\vartheta(2\kappa)$,
upper bound on $\Risk^+_{\|\cdot\|_\cB}[\widehat{x}_H|\cX]$; here $\kappa=\min[m,\nu,\max_s\rank A_s]$. In particular, linear estimate $\widehat{x}_{H_*}$  yielded by an optimal solution to \rf{(!)} is optimal
within the above factor, in terms of its risk $\Risk^+_{\|\cdot\|_\cB}[\cdot|\cX]$, among all linear estimates. Finally, when there is no uncertainty ($A_s=0$ for all $s$), $\widehat{x}_{H_*}$ is exactly the near-minimax-optimal estimate from \cite[Proposition 4.\aic{8.1}{16}]{STOPT}.

%\subsection{Application to  system  identification}\label{systident}
\subsubsection{The problem}
Consider situation as follows: a linear time-invariant dynamical system with states $u_t\in\bR^d$ {and inputs $r_t\in\bR^h$} evolves according to
\begin{equation}\label{system}
u_{t+1}={X[u_t;r_t]}.
\end{equation}
We are given noisy observations $\ov{u}_t$  of the states on time horizon $0\leq t\leq N$ and of the inputs on time horizon $0\leq t<N$:
\begin{equation}\label{obsstates}
\ov{u}_{ti}=u_{ti}-\xi_{ti},\;0\leq t\leq N,1\leq i\leq d;\quad
{\ov{u}_{ti}=r_{t,i-d}-\xi_{ti},\,0\leq t<N,d<i\leq d+h,}
\end{equation}
We have at our disposal upper bounds on the magnitudes of observation errors:
$$
{|\xi_{tj}|\leq\overline{\xi}_{tj}}
$$
with known $\overline{\xi}$'s.
In addition, we have partial a priori knowledge of $X$ expressed by a system of linear equations on the entries of $X$. Our goal is to recover the image $X^+$ of $X$ under a given linear mapping.
\par
Observe that the considered setting is rather different from the ``classical'' setting of linear system identification problem, cf.
\cite{aastrom1971system,eykhoff,mehra1977system,Ljung1997,sos1994}, in which
it is assumed that the states of the system are observed without errors, and the errors in observations of inputs are corrupted by random zero mean noise. The situation in which perturbations in the observation of the state of the system are {\em uncertain-but-bounded} (e.g., belong to an ellipsoid) is the subject of the significant literature (see, e.g., \cite{bertsekas1971recursive,casini2014feasible,cerone1993feasible,jaulin2001interval,kurzhanskii1977control,kurzhanskii1991identification,kurzhansky1997ellipsoidal,matasovestimators,milanese2013bounding,nazin2007ellipsoid,schweppe1973uncertain,walter1990special}
and references therein). The ``generic'' approach to the problem we develop below, to the best of our knowledge,  differs significantly from those proposed so far, and, we believe, can be considered as a meaningful contribution to the this line of research.
\par
Assigning the entries of $X$ serial indices, denoting  by $\iota(i,j)$ the index of $X_{ij}$ and setting \anc{$x^*_0=1$,}{} $x^*_{\iota(i,j)}=X_{ij}$, we get $n$-dimensional vector $x^*$,
$n=d(d+h)$, known to satisfy
the system of linear equations
\begin{equation}\label{eq(a)}
\anc{x_0=1,\,Px=0}{Px=p}
\end{equation}
{$Px=p$ ($P\in\bR^{\nu\times n}$ has linearly independent rows)} expressing our a priori information on the actual entries of $X$.
\anc{\footnote{The set of corresponding constraints can always can be written down as a homogeneous system in the extended by $x^*_0=1$ vector of entries of $X$}}{}
 Dynamic equations read
{\small\begin{align*}
\overline{u}_{t+1,i}&=
{\sum}_{j=1}^{d+h}\overline{u}_{tj}x^*_{\iota(i,j)}-\xi_{t+1,i}\anc{x^*_0}{}+{\sum}_{j=1}^{d+h}\xi_{tj}x^*_{\iota(i,j)},\;{1\leq i\leq d,\;0\leq t\leq N-1,}&\qquad{(!_{ti})}
\end{align*}}\noindent
which we rewrite as a system of linear equations on $x^*$ of the form\\
\[
Qx-{\sum}_{s=1}^S\zeta_s Q_s x= \anc{q}{q-{\sum}_{s=1}^S\zeta_sq_s}
\]
where $S=2N(d+h)+d$ is the total count of observation errors {$\xi_{tj}$}, $\zeta_1,...,\zeta_S$ are these errors written
down in certain order, $Q$ and $Q_s$ are observable $m\times n$ matrices, $m=dN$,  and \anc{$q$ is}{$q$, $q_s$ are} observable $m$-dimensional
\anc{vector.}{vectors.} Note that each matrix $Q_s$ has at most $d+1$ nonzero rows. Indeed,
{observation error $\xi_{tj}$ with $j\leq d$ participates only in equation $(!_{t-1,j})$ (this happens when $t\geq1$) and $d$ equations $(!_{ti})$, $1\leq i\leq d$, and observation error $\xi_{tj}$
with $j>d$ participates only in $d$ equations $(!_{ti})$, $1\leq i\leq d$.}
Setting
$$
\cL=\{x\in\bR^n:Px=0\},\,\bar{x}=(PP^T)^{-1}P^Tp,\, \Pi=I_n-P^T(PP^T)^{-1}P,
$$
so that $\Pi$ is an orthoprojector of $\bR^n$ onto $\cL$ and $\bar{x}$ is the orthoprojection of $x^*$ onto the orthogonal complement of $\cL$, we have
$$
x^*=\bar{x}+\Delta^*
$$
with $\Delta^*$ satisfying the relations
$$
\begin{array}{c}
\Delta^*\in\cL,\\
\exists (\epsilon^*\in\bR^S,\|\epsilon^*\|_\infty\leq1):
 Q[\bar{x}+\Delta^*]-\big[\sum_{s=1}^S\epsilon^*_s\underbrace{\overline{\zeta}_sQ_s}_{\overline{Q}_s}\big][\bar{x}+\Delta^*]
=q-\epsilon^*_s \overline{\zeta}_sq_s]\\
\end{array}
$$
Thus, $x^*=\bar{x}+\Delta^*$, where $\Delta^*$ solves, for properly selected vector $\epsilon=\epsilon^*\in\bR^S$, $\|\epsilon^*\|_\infty\leq1$, the system of linear equations
\begin{equation}\label{solveslinsyst}
\begin{array}{c}
[Q-\sum_s\epsilon_s\overline{Q}_s]\Delta=[\overline{q}+\sum_s\epsilon_s\overline{q}_s]\ \&\ \Pi\Delta=\Delta
\\
\left[\overline{q}=q-Q\bar{x},\,\overline{q}_s=\overline{Q}_s\bar{x}-\overline{\zeta}_sq_s\right]\\
\end{array}\end{equation}
in variables $\Delta\in\bR^n$.
\par
Recall that out goal is to recover from observation the image of $X$ under a given linear mapping; this is the same as to recover
$$y^*=Bx^*=\underbrace{B\bar{x}}_{\bar{y}}+\underbrace{B\Delta^*}_{\delta^*}$$ for a given $\nu\times n$ matrix $B$. Let us quantify the recovery error
by the norm $\|\cdot\|_\cB$ on $\bR^\nu$.
\subsubsection{Robust linear recovery}\label{roblinrec}
Given $m\times n$ matrix $E$ and $m\times \nu$ matrix $H$, let us recover\begin{itemize}
\item  $\Delta^*$ by the vector
$$
\begin{array}{rcl}
\widehat{\Delta}&:=&\Pi E^T\overline{q}
=\left[\Pi E^TQ-\sum_s\epsilon_s^*\Pi E^T\overline{Q}_s\right]\Delta^*-\sum_s\epsilon^*_s\Pi E^T\overline{q}_s\\
\end{array}
$$
and $x^*$---by the vector $\bar{x}+\widehat{\Delta}$,
\item $\delta^*$ by the vector
$$
\begin{array}{rcl}
\widehat{\delta}&:=&H^T\overline{q}
=\left[H^TQ-\sum_s\epsilon_s^*H^T\overline{Q}_s\right]\Delta^*-\sum_s\epsilon^*_sH^T\overline{q}_s\\
\end{array}
$$
and $y^*$---by the vector $\bar{y}+\widehat{\delta}$.
\end{itemize}\par
{\bf Performance analysis.} By (\ref{solveslinsyst}) we have
$$
\overline{q}=[Q-{\sum}_s\epsilon^*_s\overline{Q}_s]\Delta^* -{\sum}_s\epsilon^*_s\overline{q}_s.
$$
Thus, $\widehat{\Delta}\in\cL$, $\Delta^*\in\cL$ and
\begin{equation}\label{nov7eq1}
\begin{array}{c}
\begin{array}{rcl}
\widehat{\Delta}-\Delta^*&=&\left[\Pi E^TQ-I_n-\sum_s\epsilon_s^*\Pi E^T\overline{Q}_s\right]\Delta^*-[\sum_s\epsilon_s^*\Pi E^T\overline{q}_s]\\
&=&\left[\Pi[E^TQ-I_n]\Pi-\sum_s\epsilon_s^*\Pi E^T\overline{Q}_s\Pi\right]\Delta^*-[\sum_s\epsilon_s^*\Pi E^T\overline{q}_s]
\end{array}
\\
\end{array}
\end{equation}
where the concluding equality is due to $\Delta^*=\Pi\Delta^*$ and $\Pi^2=\Pi$. Besides this,
\begin{equation}\label{dec7eq1}
\begin{array}{c}
\begin{array}{rcl}
\widehat{\delta}-\delta^*&=&\left[H^TQ-B-\sum_s\epsilon_s^*H^T\overline{Q}_s\right]\Delta^*-[\sum_s\epsilon_s^*H^T\overline{q}_s]\\
&=&\left[[H^TQ-B]\Pi-\sum_s\epsilon_s^*H^T\overline{Q}_s\Pi\right]\Delta^*-[\sum_s\epsilon_s^*H^T\overline{q}_s]
\end{array}
\end{array}
\end{equation}
\par
Now let $\cX$ be the unit ball of a norm $\|\cdot\|_\cX$ on $\bR^n$; assume that this norm is both ellitopic and co-ellitopic. Let
$$
\begin{array}{rcl}
\cW_0[E]&=&\{\sum_s\epsilon_s\Pi E^T\overline{q}_s:\|\epsilon\|_\infty\leq1\}\subset\bR^n,\\
\cW[E]&=&\{\Pi[E^TQ-I_n]\Pi-\sum_s\epsilon_s\Pi E^T\overline{Q}_s\Pi:\|\epsilon\|_\infty\leq1\},\\
\end{array}
$$
and let $\Upsilon_0[E]$ and $\Upsilon[E]$ be the efficiently computable convex in $E$ upper bounds, given by our machinery, on the robust norms
{\footnotesize$$
\|\cW_0[E]\|_{\cX,[-1,1]}=\max\limits_w\{\|w\|_\cX:w\in\cW_0[E]\},\,\,\|\cW[E]\|_{\cX,\cX}=\max\limits_{W}\{\|W\|_{\cX,\cX}:W\in\cW[E]\}
$$}\noindent
of the uncertain $n\times 1$ matrix $\cW_0[E]$ and uncertain $n\times n$ matrix $\cW[E]$. By (\ref{nov7eq1}) we have
\begin{equation}\label{nov7eq2}
\|\widehat{\Delta}-\Delta^*\|_\cX\leq \Upsilon[E] \|\Delta^*\|_\cX+\Upsilon_0[E].
\end{equation}
\par
Assume from now on that $\|\cdot\|_\cB$ is a co-ellitopic norm, let
$$
\begin{array}{rcl}
\cV_0[H]&=&\{\sum_s\epsilon_sH^T\overline{q}_s:\|\epsilon\|_\infty\leq1\}\subset\bR^n,\\
\cV[H]&=&\{[H^TQ-B]\Pi-\sum_s\epsilon_sH^T\overline{Q}_s\Pi:\|\epsilon\|_\infty\leq1\}\\
\end{array}
$$
and let $\overline{\Upsilon}_0[H]$, $\overline{\Upsilon}[H]$ be the efficiently computable convex in $H$ upper bounds, given by our machinery, on the robust norms
{\footnotesize$$
\|\cV_0[H]\|_{\cB,[-1,1]}=\max\limits_w\{\|w\|_\cB:w\in\cV_0[H]\},\,\,\|\cV[H]\|_{\cB,\cX}=\max\limits_{W}\{\|W\|_{\cB,\cX}:W\in\cW[H]\}
$$}\noindent
of the uncertain $\nu\times 1$ matrix $\cV_0[H]$ and uncertain $\nu\times n$ matrix $\cV[H]$.
By (\ref{dec7eq1}) we have
\begin{equation}\label{dec7eq2}
\|\widehat{\delta}-\delta^*\|_\cB\leq \overline{\Upsilon}[H] \|\Delta^*\|_\cX+\overline{\Upsilon}_0[H].
\end{equation}

Assume now that $E$ is such that $\Upsilon[E]<1$. Then
$$
\|\Delta^*\|_\cX\leq \|\widehat{\Delta}-\Delta^*\|_\cX+\|\widehat{\Delta}\|_\cX\leq \Upsilon[E] \|\Delta^*\|_\cX +[\|\widehat{\Delta}\|_\cX+\Upsilon_0[E]]
$$
whence
\begin{equation}\label{dec07}
\|\Delta^*\|_\cX\leq {\|\widehat{\Delta}\|_\cX+\Upsilon_0[E]\over1-\Upsilon[E]}.
\end{equation} As a result,
\begin{equation}\label{nov7eq3}
\begin{array}{rcll}
\|\widehat{x}-x^*\|_\cX=\|\widehat{\Delta}-\Delta^*\|_\cX
&\leq& {\Upsilon[E]\over 1-\Upsilon[E]}[\|\widehat{\Delta}\|_\cX+\Upsilon_0[E]]+\Upsilon_0[E]&\\
&=& {\Upsilon[E]\over 1-\Upsilon[E]}\left[\|\Pi E^T\overline{q}\|_\cX+\Upsilon_0[E]\right]+\Upsilon_0[E]&(a)\\
\|\widehat{y}-y^*\|_\cB=\|\widehat{\delta}-\delta^*\|_\cB
&\leq& {\overline{\Upsilon}[H]\over 1-\Upsilon[E]}[\|\widehat{\Delta}\|_\cX+\Upsilon_0[E]]+\overline{\Upsilon}_0[H]&\\
&=& {\overline{\Upsilon}[H]\over 1-\Upsilon[E]}\left[\|\Pi E^T\overline{q}\|_\cX+\Upsilon_0[E]\right]+\overline{\Upsilon}_0[H]&(b)
\end{array}
\end{equation}
\par
{\bf Synthesis of linear estimate.} Recall that  the problem of minimizing $\Upsilon[E]$ w.r.t. $E$ is efficiently solvable. If we are lucky to have $\Upsilon_*:=\inf_E\Upsilon[E]<1$,
we can optimize, to some extent, our estimate $H^T\overline{q}$  of $y^*=Bx^*$ in $H$. To this end let us select $E$ which ``nearly minimizes'' the quantity
$$
\Gamma={1\over 1-\Upsilon[E]}\left[\|\Pi E^T\overline{q}\|_\cX+\Upsilon_0[E]\right]
$$
over $E$ under the constraint $\Upsilon(E)<1$; after $E$ is selected, we specify $H$ by minimizing the resulting right hand side of (\ref{nov7eq3}.$b$), that is,
$
\Gamma\overline{\Upsilon}(H) + \overline{\Upsilon}_0[H]
$
in $H$.
\par
``Near-minimization'' of $\Gamma$ over $E$ can be carried out as follows. Let us  select somehow
$\beta<1$ close to 1 (e.g., $\beta=0.9$ or $\beta=0.99$) and set $\Upsilon_i=(1-\beta^i)+\beta^i\Upsilon_*$, $i=0,1,2,...$, so that ${\beta^i\over 1-\Upsilon}\leq {1\over 1-\Upsilon_*}$ is equivalent to $\Upsilon\leq\Upsilon_i$.  We solve one by one feasible convex optimization problems
$$
\Opt_i={1\over 1-\Upsilon_i}\min\limits_E\bigg\{\|\Pi E^T\overline{q}\|_\cX+\Upsilon_0[E]: \Upsilon[E]\leq \Upsilon_i\bigg\}.\leqno{(P_i)}
$$
$i=0,1,...$, run this process until the quantities $\Opt_i$ start to grow, and specify $\Gamma$ as the smallest of the quantities $\Opt_i$ we have generated.
\par
Let us write explicitly the problem ($P_i$) in  the {situation} where
\begin{equation}\label{ssimplecase}
\cX=\Conv\Big\{{\bigcup}_{k\leq K} P_k{\cB_{n_k}}\Big\}
\end{equation}
where $\cB_m$ is  the unit $\|\cdot\|_2$-ball in $\bR^m$ and $P_k\in\bR^{n\times {n_k}}$.
As we know, in this case
\be\begin{array}{c}
\|x\|_\cX=\min\limits_{x_k\in\bR^{n_k},k\leq K}\left\{\sum_k\|x_k\|_2:\sum_kP_kx_k=x\right\},\\
\|x\|_{\cX_*}={\max}_{k\leq K} \|P_k^Tx\|_2,\,\,
\|A\|_{\cX,\cX}=\max_{k\leq K}\|AP_k\|_{\cX,\B_{n_k}}.
\end{array}\ee{simplenorms}\noindent
Exploiting the fact that in our present situation $\cX_*=\{x:\|P_k^Tx\|_2\leq 1,k\leq K\}$ is an ellitope,
($P_i$) may be rewritten as follows (cf. \rf{eqmixture1} in Proposition \ref{propmixture}):
{\footnotesize$$
\begin{array}{rcl}
\|\Pi E^T\overline{q}\|_\cX
&=&\min\limits_{\{x_k,k\leq K\}}\left\{\sum_k\|x_k\|_2:\sum_kP_kx_k=\Pi E^T\overline{q}\right\};\\
\Upsilon_0[E]&=&\min\limits_{\{G_s\in\bR^n,H_s\in\bR,s\leq S\},\atop
\upsilon\in\bR^K,\lambda\in\bR}\bigg\{{1\over 2}[\sum_k\upsilon_k+\lambda]:\\
&&\left\{\begin{array}{l}\upsilon\geq0,\sum_sG_s\preceq \sum_k\upsilon_kP_kP_k^T,\sum_sH_s\leq\lambda \\
\hbox{\scriptsize$\left[\begin{array}{c|c}G_s&\Pi E^T\overline{q}_s\cr\hline
\overline{q}_s^TE\Pi&H_s\cr\end{array}\right]$}\succeq0\\
\end{array}\right.\bigg\}\\
&=&\min\limits_{\{H_s\in\bR,s\leq S\},\atop
\upsilon\in\bR^K,\lambda\in\bR}\bigg\{{1\over 2}[\sum_k\upsilon_k+\lambda]:\\
&&\left\{\begin{array}{l}\upsilon\geq0,\sum_sH_s\leq\lambda\\
\sum_sH_s^{-1}[\Pi E^T\overline{q}_s][\Pi E^T\overline{q}_s]^T\preceq \sum_k\upsilon_kP_kP_k^T\\
\end{array}\right.\bigg\}\\
&=&\min\limits_{\upsilon\in\bR^K,\mu\in\bR^S}\bigg\{{1\over 2}[\sum_k\upsilon_k+\sum_s\mu_s]:\\
&&\upsilon\geq0,\hbox{\scriptsize$\left[\begin{array}{c|c}\sum_k\upsilon_kP_kP_k^T&\Pi E^T[\overline{q}_1,...,\overline{q}_S]\cr\hline
[\overline{q}_1,...,\overline{q}_S]^TE\Pi&\Diag\{\mu\}\cr\end{array}\right]$}\succeq 0\bigg\};\\
\Upsilon[E]&=&\min\limits_{\{G^k_s,H^k_s,\ov G^k,\ov H^k:k\leq K,s\leq S\},\atop
\{\upsilon^k,\lambda^k:k\leq K\}}\bigg\{{1\over 2}\max\limits_{k\leq K}[\sum_{j=1}^K\upsilon^k_j+\lambda^k]:\\
&&\left\{\begin{array}{l}\upsilon^k\geq0,\overline{G}^k+\sum_sG^k_s\leq \sum_{j=1}^K\upsilon^k_jP_jP_j^T,\overline{H}^k+\sum_sH^k_s\preceq \lambda^kI_{n_k}\\
\hbox{\scriptsize$\left[\begin{array}{c|c}\overline{G}^k&\Pi [E^TQ-I_n]\Pi P_k\cr\hline
P_k^T\Pi [Q^TE-I_n]\Pi&\overline{H}^k\cr\end{array}\right]$}\succeq0,\\
\hbox{\scriptsize$\left[\begin{array}{c|c}G^k_s&\Pi E^T\overline{Q}_s\Pi P_k\cr\hline
P_k^T\Pi \overline{Q}_s^TE\Pi&H^k_s\cr\end{array}\right]$}\succeq0,\,k\leq K,s\leq S\\
\end{array}\right.\bigg\};\\
\end{array}
$$}
{\footnotesize$$
\begin{array}{rcl}
\Opt_i&=&{1\over 1-\Upsilon_i}\min\limits_{E,\{x_k,k\leq K\},\upsilon\in\bR^K,\mu\in\bR^S,
\atop {\{G^k_s,H^k_s,\ov G^k,\ov H^k,k\leq K,s\leq S\},\atop
\{\upsilon^k\in\bR^K,\lambda^k,k\leq K\}}}\bigg\{\sum_k\|x_k\|_2+{1\over 2}[\sum_k\upsilon_k+\sum_s\mu_s]:\\
&&\left\{\begin{array}{l}\sum_kP_kx_k=\Pi E^T\overline{q}\\
\upsilon\geq0,\hbox{\scriptsize$\left[\begin{array}{c|c}\sum_k\upsilon_kP_kP_k^T&\Pi E^T[\overline{q}_1,...,\overline{q}_S]\cr\hline
[\overline{q}_1,...,\overline{q}_S]^TE\Pi&\Diag\{\mu\}\cr\end{array}\right]$}\succeq 0\\
\upsilon^k\geq0,\overline{G}^k+\sum_sG^k_s\leq \sum_{j=1}^K\upsilon^k_jP_jP_j^T,\overline{H}^k+\sum_sH^k_s\preceq \lambda^kI_{n_k},k\leq K,\\
\hbox{\scriptsize$\left[\begin{array}{c|c}\overline{G}^k&\Pi [E^TQ-I_n]\Pi P_k\cr\hline
P_k^T\Pi [Q^TE-I_n]\Pi&\overline{H}^k\cr\end{array}\right]$}\succeq0,\\
\hbox{\scriptsize$\left[\begin{array}{c|c}G^k_s&\Pi E^T\overline{Q}_s\Pi P_k\cr\hline
P_k^T\Pi \overline{Q}_s^TE\Pi&H^k_s\cr\end{array}\right]$}\succeq0,\,k\leq K,s\leq S\\
{1\over 2}\max\limits_{k\leq K}[\sum_j\upsilon^k_j+\lambda^k]\leq\Upsilon_i\\
\end{array}\right.\bigg\}
\end{array}
$$}\noindent

\begin{remark}\label{simplenorm} {\rm Rationale behind restricting ourselves to $\cX$ as in \rf{ssimplecase} is as follows. Recall that the norm $\|\cdot\|_\cX$ we consider is  assumed to be both ellitopic and co-ellitopic. There are only two known to us generic situations in which the corresponding unit ball $\cX$ is both ellitopic and co-ellitopic at the same time, and \rf{ssimplecase} is one of them.
The other nice situation, ``symmetric'' to the first, is when $\|\cdot\|_\cX$ is the conjugate of the norm just defined, that is, norm of the form $\max_{k\leq K}\|P_k^Tx\|_2$.
In our context, this second case reduces to the first due to $\|A\|_{\cX,\cX}=\|A^T\|_{\cX_*,\cX_*}$.
}\end{remark}
\par{\sl Numerical illustration} to follow deals with recovery of the parameters of the ``Boeing 747'' model used in Section \ref{peak-to-peak}, which in our present notation reads
\def\LS{{\hbox{\tiny LS}}}
\\[5pt]
\centerline{$u_{t+1}=
\hbox{\tiny$\left[\begin{array}{rrrrrrrr}
0.9957&0.0339&-0.0211&-0.3214&0.0140&0.9886&0.0043&-0.0337\\
    0.0076&0.4699&4.6604&0.0022&-3.4373&1.6648&-0.0079&0.5285\\
    0.0168&-0.0605&0.4038&-0.0029&-0.8219&0.4378&-0.0167&0.0600\\
    0.0091&-0.0370&0.7194&0.9990&-0.4735&0.2491&-0.0091&0.0370\\
\end{array}\right]$}[u_t;r_t]$}\vskip5pt
\noindent
where $u_t\in\bR^4$ are the states, and $r_t\in\bR^4$ are the inputs (``in reality'' the first two entries in $r_t$ are controls, and the last two---external disturbances). We observe $u_t$'s for $0\leq t\leq N=12$ and
$r_t$'s for $0\leq t<N$; in the resulting identification problem,  $m=52$, $n=32$, $S=100$, and $\cL=\bR^n$  (whence $\Pi=I_n$ and  $\bar{x}=0$). Observations of states and inputs are corrupted by ``relative $\epsilon$-noises,''
so that an observable real $r$
and its observation $\overline{r}$ satisfy $|r-\overline{r}|\leq \epsilon\max[|\overline{r}|,1]$. In an experiment, we select a noise level $\epsilon\in(0.001,0.01]$, generate a sample trajectory of the system by selecting at random the initial state and the inputs, then
add to the states and the inputs random $\epsilon$-errors,  and apply to the resulting observations the above robust linear recovery with $B=I_n$ and $\cB=\cX$ being the unit $\|\cdot\|_2$-ball in $\bR^n$
to recover the parameters of the system. We have compared this recovery with the simplest Least Squares recovery $E_\LS^Tq=\argmin_x\|Qx-q\|_2^2$, $E_\LS=Q(Q^TQ)^{-1}$.
\par
\begin{figure}{\tiny
$$
\epsfxsize=350pt\epsfysize=140pt\epsffile{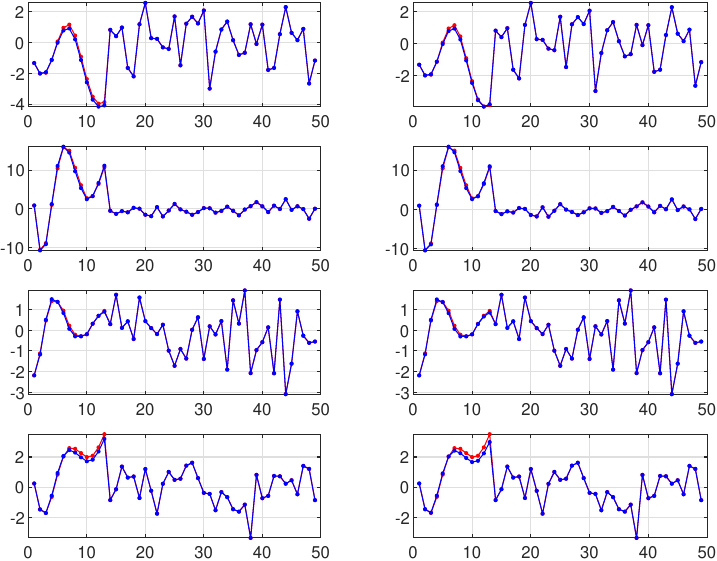}
$$
\caption{\label{figLELS} Experiment \# 10: states of the actual (blue) and the recovered (red) systems vs. time. Left: Least Squares recovery; right: Robust recovery} }
\end{figure}

The results of 10 experiments  are presented in Table \ref{tableLELS}.  In Figure \ref{figLELS}, we present the trajectories of the actual and the recovered (in experiment \# 10, $\epsilon-0.01$) systems on time horizon
$1\leq t\leq 49 $ for random initial state and inputs (different from those used in the experiment).
\vskip-3pt
\begin{table}
{\tiny$$
\begin{array}{|c||c|c|c|c|c|}
\hline
\epsilon&0.001&0.002&0.003&0.004&0.005\\
\hline
\hbox{Least Squares}&0.040/0.498&0.026/0.404&0.032/ Inf&0.065/0.753&0.103/1.889\\
\hline
\hbox{Linear  recovery}&0.043/0.338&0.040/0.320&0.054/0.455&0.058/0.555&0.177/1.447\\
\hline
\epsilon&0.006&0.007&0.008&0.009&0.010\\
\hline
\hbox{Least Squares}&0.077/0.613&0.114/1.224&0.334/5.034&0.110/1.269&0.165/1.744\\
\hline
\hbox{Linear  recovery}&0.100/0.524&0.109/1.035&0.214/3.749&0.126/1.104&0.134/1.513\\
\hline
\end{array}
$$
\caption{\label{tableLELS}  $\|\cdot\|_2$ recovery errors (first numbers in cells), and upper bounds on $\|\cdot\|_2$ recovery errors  as given by (\ref{nov7eq3}) (second numbers in cells)} as functions of noise level $\epsilon$.
}
\end{table}

%\end{center}

\appendix
\section{Proofs}
\subsection{Proof of Theorem \ref{verynewRelaxationTheorem}}
The below proof follows that of Theorem \ref{2020Prop4.6} as given in \cite[Section 4.8.2]{STOPT}, utilizing at some point bilinearity of the quadratic form we want to upper-bound on $\cZ\times\cW$.
\par
Let $q$, $p$ be the dimensions of the embedding spaces of $\cZ$ and $\cW$, and assume w.l.o.g. that $q\leq p$.\footnote{It is immediately seen that the norm bound (\ref{normbound}) is intelligent enough to respect the identity
$\|A\|_{\cB,\cX}=\|A^T\|_{\cX_*,\cB_*}$ where $\cQ_*$  stands for the polar of a set $\cQ$. As a result, to ensure $q\leq p$, we can pass, if necessary, from $\cB,\cX$ and $A$ to $\cX_*,\cB_*$ and $A^T$.}
\par{\bf 1$^o$.} Let
$$
{\mathfrak{T}}=\clsr\{[t;\tau]: \tau>0,t/\tau\in\cT\}\;\mbox{and}\;{\mathfrak{R}}=\clsr\{[r;\theta]:\theta>0,r/\theta\in\cR\}
$$
be the closed conic hulls of $\cT$ and ${\cR}$, so that ${\mathfrak{T}}$ and ${\mathfrak{R}}$ are regular (closed, pointed and convex with nonempty interior) cones such that
$$
\cT=\{t:[t;1]\in{\mathfrak{T}}\},\;\; \cR=\{r:[r;1]\in{\mathfrak{R}}\}.
$$
As is immediately seen, the cones dual to ${\mathfrak{T}}$, ${\mathfrak{R}}$ are
$$
{\mathfrak{T}}_*=\{[g;\tau]: \tau\geq \phi_\cT(-g)\},\,\,{\mathfrak{R}}_*=\{[h;\theta]: \theta\geq \phi_\cR(-h)\}.
$$
In view of these  observations, (\ref{normbound}) is nothing but the conic problem
$$
\Opt(A)=\min\limits_{\lambda,\upsilon,\tau,\theta}\left\{\tau+\theta:\begin{array}{l}\lambda\geq0,\upsilon\geq0,
[-\lambda;\tau]\in{\mathfrak{T}}_*,[-\upsilon;\theta]\in{\mathfrak{R}}_*,\\
\left[\begin{array}{c|c}
\sum_\ell\upsilon_\ell R_\ell&{1\over 2}Q^TAP\cr\hline{1\over 2}P^TA^TQ&\sum_k\lambda_kT_k\cr\end{array}\right]\succeq0\\
\end{array}\right\}.
$$
It is easily seen that this problem is strictly feasible and bounded. By Conic Duality,
{\small\bse
\Opt(A)
&=&{\max}_{r,t,U,V,W}\left\{\Tr(W^TQ^TAP):\begin{array}{l}t\in\cT,r\in\cR\\
\Tr(R_\ell U)\leq r_\ell\,\forall \ell,
\Tr(T_k V)\leq t_k\,\forall k\\
\left[\begin{array}{c|c}U&W\cr\hline W^T&V\cr\end{array}\right]\succeq0\\
\end{array}\right\}\\
&=&{\max}_{r,t,U,V,Y}\left\{\Tr(\vec{}[U^{1/2}YV^{1/2}]^TQ^TAP):\begin{array}{l}r\in\cR,t\in\cT,U\succeq0,V\succeq0,Y^TY\preceq I\cr
\Tr(R_\ell U)\leq r_\ell\,\forall \ell,\Tr(T_kV)\leq t_k\,\forall k\\
\end{array}\right\}\\
&=&{\max}_{r,t,U,V}\left\{\sum_{i=1}^q\sigma_i(U^{1/2}Q^TAPV^{1/2}):\begin{array}{l}
U\succeq0,\Tr(R_\ell U)\leq r_\ell\,\forall \ell,r\in\cR\\V\succeq0,\Tr(T_kV)\leq t_k\,\forall k,t\in\cT\end{array}\right\}
\ese}\noindent
where $\sigma_i(\cdot)$, $i\leq q$,  are the singular values of $q\times p$ matrix (recall that $q\leq p$). At the last two steps of the above derivation, we have used the following well known facts
\begin{itemize}
\item $\hbox{\tiny$\left[\begin{array}{c|c}U&W\cr\hline W^T&V\cr\end{array}\right]$}\succeq0$ if and only if
$U\succeq0,V\succeq0$ and $W=U^{1/2}YV^{1/2}$ with $Y^TY\preceq I$,\\
and
\item the maximum of Frobenius inner products of a given matrix with matrices of spectral norm not exceeding 1 is the nuclear norm of the matrix---the sum of singular values.
\end{itemize}
\par{\bf 2$^o$.} The concluding optimization problem in the above chain clearly is solvable; let $U,V,r,t$ be the optimal solution, and let $\sigma_i=\sigma_i(U^{1/2}Q^TAPV^{1/2})$, and $\sum_{\iota=1}^q\sigma_\iota e_\iota f_\iota^T$ be the singular value decomposition of
$U^{1/2}Q^TAPV^{1/2}$, so that
\begin{equation}\label{sothatmay2021}
\begin{array}{rcl}
\Opt(A)&=&\sum_{\iota=1}^q\sigma_\iota,\\
U^{1/2}Q^TAPV^{1/2}&=&\sum_{\iota=1}^q\sigma_\iota e_\iota f_\iota^T,\\
e_i^Te_j&=&\left\{\begin{array}{ll}1,&i=j\\
0,&i\neq j\cr\end{array}\right.,i,j\leq q\ \&\ f_i^Tf_j=\left\{\begin{array}{ll}1,&i=j\\
0,&i\neq j\cr\end{array}\right.,i,j\leq p.
%\cr
%\Tr(U^{1/2}R_\ell U^{1/2})&\leq&s_\ell,\,\ell\leq L\ \&\ r\in\cR&(d.1)\\
%\Tr(V^{1/2}T_k V^{1/2})&\leq&t_k,\,k\leq K\ \&\ t\in\cT&(d.2)\\
\end{array}
\end{equation}
Let $\epsilon_1,...,\epsilon_p$ be independent random variables taking values $\pm1$ with probabilities 1/2, and let
$$\xi={\sum}_{i=1}^q\epsilon_ie_i,\, \eta={\sum}_{j=1}^p\epsilon_jf_j.$$
 Then in view of (\ref{sothatmay2021}) it holds, identically in $\epsilon_i=\pm1$, $1\leq i\leq p$:
\begin{equation}\label{identically}
\xi^TU^{1/2}Q^TAPV^{1/2}\eta={\sum}_{i,\iota\leq q,j\leq p}
[\epsilon_i\epsilon_j\sigma_\iota e_i^Te_\iota f_\iota^Tf_j]={\sum}_{\iota=1}^q\sigma_\iota=\Opt(A).
\end{equation}
On the other hand, setting $E=[e_1,...,e_q]$, we get an orthonormal $q\times q$ matrix such that $\xi=E\underline{\epsilon}$,
 where $\underline{\epsilon}=[\epsilon_1;...;\epsilon_q]$ is a Rademacher vector (i.e., random vector with independent entries taking values $\pm1$ with probabilities $1/2$), and
$$
\xi^TU^{1/2}R_\ell U^{1/2}\xi=\underline{\epsilon}^T\underbrace{[E^TU^{1/2}R_\ell U^{1/2}E]}_{\overline{R}_\ell}\underline{\epsilon}
$$
By construction, $\overline{R}_\ell\succeq0$
{ and
\[\Tr(\overline{R}_\ell)=\Tr(U^{1/2}R_\ell U^{1/2})=\Tr(R_\ell U)\leq r_\ell.\]
}
For every $\ell$ such that $r_\ell>0$ we have $\Tr(r_\ell^{-1}\overline{R}_\ell)\leq1$. Now let us use the following fact.
%\begin{quote}
{\begin{lemma}\label{2020minilemma} {\rm \cite[Lemma 4.48]{STOPT}} Let $Q$ be positive semidefinite $N\times N$ matrix with trace $\leq1$ and $\zeta$ be $N$-dimensional Rademacher random vector. Then
$$
\bE\left\{\exp\left\{\third\zeta^TQ\zeta\right\}\right\}\leq \sqrt{3}.
$$
\end{lemma}}\noindent
%\end{quote}
By Lemma \ref{2020minilemma}, whenever $r_\ell>0$ we have
$$
\bE\{\exp\{\xi^T[r_\ell^{-1}U^{1/2}R_\ell U^{1/2}]\xi/3\}\}=\bE\{\exp\{\underline{\epsilon}^T[r_\ell^{-1}\overline{R}_\ell]\underline{\epsilon}/3\}\}\leq\sqrt{3}.
$$
As a result, for every $\ell$ such that $r_\ell>0$ we have
$$
\Prob\{\xi^TU^{1/2}R_\ell U^{1/2}\xi>3\ln(4L)r_\ell\}<1/(2L).
$$
The latter relation holds true for those $\ell$ for which $r_\ell=0$ as well, since for these $\ell$ one has $U^{1/2}R_\ell U^{1/2}=0$ because trace of the latter
positive semidefinite matrix is $\leq r_\ell$.  Similar reasoning  with $\overline{\epsilon}=[\epsilon_1;...;\epsilon_p]$ in the role of $\underline{\epsilon}$ and $T_k$, $t_k$ in the roles of $R_\ell$, $r_\ell$ demonstrates that for every $k$ we have
$$
\Prob\{\eta^TV^{1/2}T_k V^{1/2}\eta>3\ln(4K)t_k\}<1/(2K).
$$
 Consequently, invoking (\ref{identically}), we conclude that there exists realization ($\overline{\xi},\overline{\eta}$) of $(\xi,\eta)$ such that
\[\overline{\xi}^TU^{1/2}Q^TAPV^{1/2}\overline{\eta}=\Opt(A),\]
and
\[
\overline{\xi}^TU^{1/2}R_\ell U^{1/2}\overline{\xi}\leq3\ln(4L)r_\ell\,\forall\ell,\;\;\;
 \overline{\eta}^TV^{1/2}T_k V^{1/2}\overline{\eta}\leq3\ln(4K) t_k\, \forall k.\]
Setting $v=QU^{1/2}\overline{\xi}$, $x=PV^{1/2}\overline{\eta}$ and invoking \rf{2020ell2*}, we get $\|x\|_\cX\leq \sqrt{3\ln(4K)}$, $\|v\|_{\cB_*}\leq\sqrt{3\ln(4L)}$, resulting in
$$
\Opt(A)=\overline{\xi}^TU^{1/2}Q^TAPV^{1/2}\overline{\eta}=v^TAx\leq \|x\|_\cX\|\upsilon\|_{\cB_*}\|A\|_{\cB,\cX},
$$ that is,
$$
\Opt(A)\leq3\sqrt{\ln(4K)\ln(4L)}\|A\|_{\cB,\cX}.
$$
as claimed.
\par{\bf 3$^o$.} It remains to consider the case of $K=L=1$. By evident scaling argument, the situation reduces to that where $\cX=\{w:w^TTw\leq1\}$ and $\cB_*=Q\{z:z^TSz\leq 1\}$. In this case,
$$
\begin{array}{rcl}
\|A\|_{\cB,\cX}&=&{\max}_{z:z^TSz\leq1,\atop w:w^TTw\leq1}z^T[Q^TAP]w={\max}_{\zeta:\|\zeta\|_2\leq1\atop \omega:\|\omega\|_2\leq1}\omega^T[S^{-1/2}Q^TAPT^{-1/2}]\zeta\\
&=&\min\limits_\nu\left\{\sqrt{\nu}:
\left[\begin{array}{c|c}\nu I_q&[S^{-1/2}Q^TAPT^{-1/2}]\cr\hline[S^{-1/2}Q^TAPT^{-1/2}]^T&I_p\cr\end{array}\right]\succeq 0\right\}.
\end{array}
$$
On the other hand,
\bse
\lefteqn{\Opt(A)=\min\limits_{\lambda,\upsilon}\left\{\lambda+\upsilon:\left[\begin{array}{c|c}
\upsilon S&{1\over 2}[Q^TAP]\cr\hline
{1\over2}[Q^TAP]^T&\lambda T\cr\end{array}\right]\succeq0\right\}}\\
&=&\min\limits_{\lambda,\upsilon}\left\{\half[\lambda+\upsilon]:\left[\begin{array}{c|c}
\upsilon I_q&[S^{-1/2}Q^TAPT^{-1/2}]\cr\hline
[S^{-1/2}Q^TAPT^{-1/2}]^T&\lambda I_p\cr\end{array}\right]\succeq0\right\}\\
&=& \min\limits_{\lambda\geq0,\upsilon\geq0}\left\{\half
[\lambda+\upsilon]:\left[\begin{array}{c|c}\upsilon\lambda&[S^{-1/2}Q^TAPT^{-1/2}]\cr\hline[S^{-1/2}Q^TAPT^{-1/2}]^T&I_p\cr\end{array}\right]\succeq0\right\}\\
&=&\min\limits_{\lambda\geq0,\upsilon\geq0,\nu}\left\{\half[\lambda+\upsilon]:\upsilon\lambda\geq\nu,\left[\begin{array}{c|c}\nu I_q&[S^{-1/2}Q^TAPT^{-1/2}]\cr\hline[S^{-1/2}Q^TAPT^{-1/2}]^T&I_p\cr\end{array}\right]\succeq0\right\}\\
&=&\min\limits_\nu\left\{\sqrt{\nu}:\left[\begin{array}{c|c}
\nu I_q&[S^{-1/2}Q^TAPT^{-1/2}]\cr\hline
[S^{-1/2}Q^TAPT^{-1/2}]^T&I_p\cr\end{array}\right]\succeq0\right\}=\|A\|_{\cB,\cX}.\qquad\hbox{\qed}
\ese
\subsection{Proof of Proposition \ref{propmatrcube}}
\par{\bf 1$^o$.}
Let $\mathfrak{R},\,\mathfrak{T},\,\mathfrak{R}_*$ and $\mathfrak{T}_*$ be as defined in item 1$^o$ of the proof of Theorem \ref{verynewRelaxationTheorem}. Observe that
{\footnotesize\bse
\lefteqn{\Opt
=\min\limits_{\lambda,\upsilon,G_s,H_s,\alpha,\beta}\left\{\alpha+\beta:\begin{array}{l}\left[\begin{array}{c|c}G_s&{1\over2} Q^TA_sP\cr\hline{1\over2} P^T A_s^TQ&H_s\cr\end{array}\right]\succeq0\,\forall s\leq S,
{[-\upsilon;\alpha]}\in\mathfrak{R}_*,[-\lambda;\beta]\in\mathfrak{T}_*\\
\sum_sG_s\preceq\sum_\ell\upsilon_\ell R_\ell,
\sum_sH_s\preceq \sum_k \lambda_k T_k, \,
\lambda\geq0,\upsilon\geq 0\\
\end{array}\right\}}\nn
&=&{\max}_{Y,X,W_s,r,t}\left\{\sum_s\Tr(W_s^TQ^TA_sP):\begin{array}{l}\left[\begin{array}{c|c} Y&W_s\cr\hline W_s^T&X\cr\end{array}\right]\succeq0\,\forall s\leq S, \,t\in\T,r\in\cR\\
\Tr(YR_\ell)\leq r_\ell,\ell\leq L,\, \Tr(XT_k)\leq t_k,k
\leq K
\end{array}\right\}\nn
&&\mbox{[by conic duality]}\nn
&=& {\max}_{Y,X,r,t}\left\{\sum_s\|\sigma(Y^{1/2}Q^TA_sPX^{1/2})\|_1:\begin{array}{l}
Y\succeq0,X\succeq0,
t\in\T,r\in\cR,\\
\Tr(YR_\ell)\leq r_\ell,\ell\leq L, \Tr(XT_k)\leq t_k,k
\leq K\\
\end{array}\right\}
\ese}\noindent
where $\sigma(A)$ is the singular spectrum of $A$; the last equality in the chain follows from the two simple observations (cf. the proof of Theorem \ref{verynewRelaxationTheorem}):
\begin{itemize}
\item {LMI {\tiny$\left[\begin{array}{c|c}P&Q\cr\hline Q^T&R\cr\end{array}\right]\succeq0$} with $p\times p$ matrix $P$ and $r\times r$ matrix $R$ takes place if and only if $P\succeq0$, $R\succeq0$, and $Q=P^{1/2}YR^{1/2}$ with $p\times r$ matrix $Y$ such that $Y^TY\preceq I_r$}, and
\item {for $p\times r$ matrix $A$, one has ${\max}_Y\{\Tr(Y^TA): Y\in\bR^{p\times r}, Y^TY\preceq I_r\}=\|\sigma(A)\|_1$}
\end{itemize}
With $\cL[B]=\hbox{\tiny$\left[\begin{array}{c|c}&\half B\cr\hline \half B^T&\cr\end{array}\right]$}$, the nonzero eigenvalues of $2\cL[B]$ are exactly plus and minus  nonzero singular values of $B$, and we conclude that
{\footnotesize\be
\Opt_b=\max\limits_{Y,X,r,t}\left\{\sum_s\|\lambda(\cL[Y^{1/2}Q^TA_sPX^{1/2}])\|_1:\begin{array}{l}
Y\in\bS^q_+,X\in\bS^p_+,t\in\T,r\in\cR\\
\Tr(YR_\ell)\leq r_\ell,\ell\leq L,\\
\Tr(XT_k)\leq t_k,k\leq K
\end{array}\right\},
\ee{May2021eq2}}\noindent
where $\lambda(A)$ is the vector of eigenvalues of a symmetric matrix $A$.

Note that $\Opt$ as defined in (\ref{May2021eq2}) clearly is a convex function of $[A_1,...,A_S]$.
\par
Observe that $\|\cA\|_{\cB,\cX}\leq\Opt$. Indeed, the problem specifying $\Opt$ clearly is solvable, and if $\lambda\geq0,\upsilon\geq0,\{G_s,H_s\}$ is its optimal solution, we have for all
 $z\in\cZ$, $w\in\cW$, $\epsilon_s=\pm1:$
\[
\epsilon_sz^TQ^TA_sPw\leq z^TG_sz+w^TH_sw.
\]
Thus,
\bse{\sum}_s\epsilon_sz^TQ^TA_sPw&\leq& z^T\left[{\sum}_\ell\upsilon_\ell R_\ell\right]z+w^T\left[{\sum}_k\lambda_k T_k\right]w\\&\leq&
\ \max_{r\in\cR,t\in\T}\left[\upsilon^Tr+\lambda^Tt\right]\leq \phi_\cR(\upsilon)+\phi_{\T}(\lambda)=\Opt
\ese
for all $w\in\cW,z\in\cZ$, and all $\epsilon_s=\pm1$, implying that $\|\cA\|_{\cB,\cX}\leq\Opt$ (recall that $P\cW=\cX$ and $Q\cZ=\cB_*$).
\par{\bf 2$^o$.}
Now, let $X\succeq0$, $Y\succeq0$, $t,r$ be such that $t\in\T$, $r\in\cR$,  $\Tr(YR_\ell)\leq r_\ell,$ $\ell\leq L$, $\Tr(XT_k)\leq t_k,$ $k
\leq K$, and
$$
\Opt={\sum}_s\|\lambda(\cL[Y^{1/2}Q^TA_sPX^{1/2}])\|_1.
$$
By \cite[Lemma 2.2]{BTNMC} (cf. \cite[Lemma 3.4.3]{LMCO}), if the ranks of all matrices $A_s$  (and thus---matrices $Q^TA_sP$) do not exceed a given $\kappa$, which we assume from now on,
 then for $\omega\sim\cN(0,I_{m+n})$ one has
$$
\bE\left\{|\omega^T\cL[Y^{1/2}Q^TA_sPX^{1/2}]\omega|\right\}\geq \|\lambda(\cL[Y^{1/2}Q^TA_sPX^{1/2}])\|_1/\vartheta(2\kappa),
$$
where $\vartheta(k)$ is defined in \rf{theta}.
It follows that for  {\small$[\eta;\xi]\sim \cN(0,\Diag\{Y,X\})$},
\[
\Opt\leq \vartheta(2\kappa)\bE\left\{{\sum}_s|\omega^T\cL[Y^{1/2}Q^TA_sPX^{1/2}]\omega|\right\}=
\vartheta(2\kappa)\bE\left\{{\sum}_s|\eta^TQ^TA_sP\xi|\right\}.
\]
Now, let $\pi(\cdot)$ be the norm on $\bR^p$ with the unit ball $\cW$, and $\rho(\cdot)$ be the norm on $\bR^q$ with the unit ball $\cZ$.
Taking into account that $\cX=P\cW$ and $\cB_*=Q\cZ$ we conclude that
$$
\forall (\eta\in\bR^q,\xi\in\bR^p): {\sum}_s|\eta^TQ^TA_sP\xi|=\max_{\epsilon_s=\pm1}\eta^TQ^T[{\sum}_s\epsilon_sA_s]P\xi\leq \rho(\eta)\pi(\xi)\|\cA\|_{\cB,\cX},
$$
thus arriving at
\begin{equation}\label{arriveat}
\Opt\leq\vartheta(2\kappa)\|\cA\|_{\cB,\cX}\bE\{\rho(\eta)\pi(\xi)\}=
\vartheta(2\kappa)\|\cA\|_{\cB,\cX}\bE\left\{\pi(\xi)\right\}\bE\left\{\rho(\eta)\right\}.
\end{equation}
\par{\bf 3$^o$.} It remains to invoke
{\begin{lemma}\label{lenlemlem}
Let
$$
\cV=\{v\in\bR^d: \exists r\in\cR: v^TR_jv\leq r_j,1\leq j\leq J\}\subset\bR^d
$$
be a basic ellitope, $W\succeq 0$ be symmetric $d\times d$ matrix such that
$$
\exists r\in\cR: \Tr(WR_j)\leq r_j,j\leq J,
$$
and $\omega\sim\cN(0,W)$. Denoting by $\rho(\cdot)$ the norm on $\bR^d$ with the unit ball $\cV$, we have
\[
\bE\{\rho(\omega)\}\leq \varkappa(J)
\]
where $\varkappa(\cdot)$ is as in \rf{May2021eq2-1}.\end{lemma}}\noindent
The statement of the proposition now follows from \rf{arriveat} by applying Lemma \ref{lenlemlem} to $\cV=\cW$, $W=X$, and to $\cV=\cZ$, $W=Y$.
\par{\bf 4$^o$.} It remains to prove Lemma \ref{lenlemlem}.  Let us start with the case of $J=1$. Setting $\bar{r}=\max\{r:r\in \cR\}$ and $R=R_1/\bar{r}$, we have $\Tr(WR)\leq 1$ and $\rho(u)=\|R^{1/2}u\|_2$. Setting $\bar{W}=R^{1/2}WR^{1/2}$ and $\bar{\omega}=R^{1/2}\omega$, we get $\bar{\omega}\sim\cN(0,\bar{W})$, $\Tr(\bar{W})\leq 1$, and
$$
\bE\{\rho(\omega)\}=\bE\{\|\bar{\omega}\|_2\}\leq \sqrt{\bE\{\bar{\omega}^T\bar{\omega}\}}=\sqrt{\Tr(\bar{W})}\leq1= \varkappa(1).
$$
Now let $J>1$. Observe that if $\Theta\succeq0$ is a $d\times d$ matrix with trace 1, $0\leq t<1/2$, and $\zeta\sim \cN(0,I_d)$ then  by convexity of $\bE\left\{\exp\{t \sum_i\zeta_i^2\lambda_i\}\right\}$ in $\lambda$
\[\bE\big\{\exp\{t \zeta^T\Theta\zeta\}\big\}= \bE\left\{\exp\{t {\sum}_i\zeta_i^2\lambda_i(\Theta)\}\right\}\\
\leq \bE_{\varsigma\sim\cN(0,1)}\{\exp\{t \varsigma^2\}\}=(1-2t)^{-1/2}.
\]
As a result,
\[\forall s\geq0: \;\Prob\left\{\zeta^T\Theta\zeta\geq s^2\}\right\}\leq{\exp\{-t s^2\}\over\sqrt{1-2t}}.
\]
Under the premise of the lemma, let $w\in\W$ be such that $\Tr(WR_j)\leq r_j$ for all $j$. For every $j$ such that $r_j>0$, setting $\Theta_j=W^{1/2}R_jW^{1/2}/r_j$, we get $\Theta_j\succeq0$, $\Tr(\Theta_j)\leq1$, so that by the above for all $s>0$ and $0\leq t<1/2$
$$
\Prob_{\omega\sim\cN(0,W)}\{\omega^TR_j\omega>s^2r_j\}=\Prob_{\zeta\sim\cN(0,I_d)}\{\zeta^T\Theta_j\zeta>s^2\}\leq {\exp\{-t s^2\}\over\sqrt{1-2t}}.
$$
The resulting inequality clearly holds true for $j$ with $r_j=0$ as well. Now, when $\omega$ and $s>0$ are such that $\omega^TR_j\omega\leq s^2r_j$ for all $j$, we have $\rho(\omega)\leq s$. Combining our observations, we get
$$
\Prob_{\omega\sim\cN(0,W)}\{\rho(\omega)>s\}\leq \min\left[1,J{\exp\{-t s^2\}\over\sqrt{1-2t}}\right],
$$
implying that
$$
\bE_{\omega\sim\cN(0,W)}\left\{\rho(\omega)\right\}\leq\int_0^\infty\min\left[1,J{\exp\{-t s^2\}\over\sqrt{1-2t}}\right]ds
$$
Optimizing w.r.t. $t$, we arrive at
$$
\bE_{\omega\sim\cN(0,W)}\left\{\rho(\omega)\right\}\leq {5\over 2}\sqrt{\ln(2J)}=\varkappa(J).\eqno{\hbox{\qed}}$$
%%%%%%%%%%%%%%%%%%%%%%%%%%%%%%%%%%
%%%%%%%%%%%%%%%%%%%%%%%
\subsection{Proof of Proposition \ref{propmixture}}
\par{\bf 0$^o$.} Equalities in \rf{eqmixture2} follow from (\ref{asfollowsa}), (\ref{asfollowsb}).
Consequently, all we need is to prove that for all $i,j$ it holds
\begin{equation}\label{allweneed}
\|\cU_{ij}\|_{\cZ_j^*,\cX_i}\leq\Opt_{ij}[\cU]\leq \max[\varsigma(K_i,L_j)+\varkappa(K_i)\varkappa(L_j)\vartheta(2\kappa)] \|\cU_{ij}\|_{\cZ_j^*,\cX_i}.
\end{equation}
\par{\bf 1$^o$.} Let us fix $i\leq I$, $j\leq J$, and let $\bar{A}_0=Q_j^TA_\n P_i$ and $\bar{A}_s=Q_j^TA_sP_i$, $1\leq s\leq S$. Setting $\overline{\cU}_{ij}=\{\sum_{s=0}^N\epsilon_s\bar{A}_s:\|\epsilon\|_\infty\leq1\}$,
we clearly have $\|\cU_{ij}\|_{\cZ_j^*,\cX_i}=\|\overline{\cU}_{ij}\|_{\cZ_j^*,\cX_i}$. Now, comparing {\rm\rf{May2021eq2-1}} with $\cZ_j^*$ in the role of $\cB$ and $\cX_i$ in the role of $\cX$ with the
definition of $\Opt_{ij}$ in (\ref{eqmixture1}), we see that $\Opt_{ij}$ is nothing but the upper bound, as given by Proposition \ref{propmatrcube}, on $\|\overline{\cU}_{ij}\|_{\cZ_j^*,\cX_i}$,
implying the left inequality in (\ref{allweneed}).
\par{\bf 2$^o$.} Observe that the upper bound on $\alpha:=\|\bar{A}_0\|_{\cZ_j^*,\cX_i}$ as given by Theorem \ref{verynewRelaxationTheorem}, is nothing but
$$
\overline{\alpha}:=\min\limits_{\lambda^\prime,\varkappa^\prime,\atop G,H}\left\{\phi_{\cT_i}(\lambda^\prime)+
\phi_{\cR_j}(\upsilon^\prime):
\begin{array}{l} \lambda^\prime\geq0,\upsilon^\prime\geq0,G\preceq{\sum}_{\ell=1}^{L_j}\upsilon_\ell^\prime  R_{\ell j},\\
H\preceq{\sum}_{k=1}^{K_i}\lambda_k^\prime T_{k i},
\left[\begin{array}{c|c}G&{1\over 2}\bar{A}_0\cr\hline
{1\over 2}\bar{A}_0^T&H\cr\end{array}\right]\succeq0\\
\end{array}\right\},
$$
and by this Theorem,
$$
\alpha\leq\overline{\alpha}\leq\varsigma(K_i,L_j)\alpha.
$$
Next, the upper bound on $\beta:=\|\cA_{ij}\|_{\cZ_j^*,\cX}$, $\cA_{ij}=\{\sum_{s=1}^N\epsilon_s\bar{A}_s:\|\epsilon\|_\infty\leq1\},$ given by Proposition \ref{propmatrcube} is
{\small$$
\overline{\beta}:=\min\limits_{\lambda^{\prime\prime},\upsilon^{\prime\prime},\atop
\{G^s,H^s,s\leq S\}}\left\{\phi_{\cT_i}(\lambda^{\prime\prime})+\phi_{\cR_j}(\upsilon^{\prime\prime}):
\begin{array}{l} \lambda^{\prime\prime}\geq0,\upsilon^{\prime\prime}
\geq0,\sum_sG^s\preceq\sum_{\ell=1}^{L_j}\upsilon_\ell^{\prime\prime}  R_{\ell j},\\
\sum_sH^s\preceq\sum_{k=1}^{K_i}\lambda_k^{\prime\prime} T_{k i},
\left[\begin{array}{c|c}G^s&{1\over 2}\bar{A}_s\cr\hline
{1\over 2}\bar{A}_s^T&H^s\cr\end{array}\right\}\succeq0\\
\end{array}\right\}.
$$}\noindent
and
$$
\overline{\beta}\leq \vartheta(2\kappa)\varkappa(K_i)\varkappa(L_j)\|\cA_{ij}\|_{\cZ_j^*,\cX}
$$
(since the ranks of matrices $\bar{A}_s$, $s\geq1$, do not exceed those of matrices $A_s$). \par
Looking at (\ref{eqmixture1}), we see that if $(\lambda^\prime,\upsilon^\prime,G,H)$,
$(\lambda^{\prime\prime},\upsilon^{\prime\prime},\{G^s,H^s\})$ are feasible solutions to the optimization problems specifying $\overline{\alpha}$ and $\overline{\beta}$, then
$$\lambda^{ij}=\lambda^\prime+\lambda^{\prime\prime},\upsilon^{ij}=\upsilon^\prime+\upsilon^{\prime\prime},G^{ijs}=G^s,H^{ijs}=H^s, \overline{G}^{ij}=G,\overline{H}^{ij}=H
$$
is a feasible solution to the problem specifying $\Opt_{ij}[\cU]$, and the value of the objective of the latter problem at this feasible solution is
$$
\phi_{\cT_i}(\lambda^\prime+\lambda^{\prime\prime})+\phi_{\cR_j}(\upsilon^\prime+\upsilon^{\prime\prime})\leq
\phi_{\cT_i}(\lambda^\prime)+\anc{\psi}{\phi}_{\cT_i}(\lambda^{\prime\prime})+\phi_{\cR_j}(\upsilon^\prime)+\phi_{\cR_j}(\upsilon^{\prime\prime}).
$$
We conclude that
$$\Opt_{ij}[\cU]\leq\overline{\alpha}+\overline{\beta}\leq \varsigma(K_i,L_j)\|\bar{A}_0\|_{\cX_j^*,\cX_i}+\varkappa(K_i)\varkappa(L_j)\vartheta(2\kappa)\|\cA_{ij}\|_{\cX_j^*,\cX_i},
$$
and since by evident reasons one has $\|\cU_{ij}\|_{\cZ_j^*,\cX_i}\geq\max\left[\|\bar{A}_0\|_{\cZ_j^*,\cX_i},\|\cA_{ij}\|_{\cZ_j^*,\cX_i}\right]$, we arrive at the right inequality in (\ref{allweneed}). \qed
%%%%%%%%%%%%%%%%%%%
\newpage
\section{Spectratopic case}
\subsection{Spectratopes} A {\sl basic spectratope} is a {\sl bounded} set $\cW$ represented as
\begin{equation}\label{s2020ell1}
\cW=\{w\in\bR^p:\exists t\in\cT: T_k^2[w]\preceq t_kI_{d_k},1\leq k\leq K\}
\end{equation}
where
$$
T_k[w]=\sum_{i=1}^p w_iT_{ki}
$$
is a linear mapping from $\bR^p$ to $\bS^{d_k}$ (so that $T_{ki}$ are symmetric $d_k\times d_k$ matrices), and $\cT$ is as in the definition of a basic ellitope.
\par
{\sl A spectratope} is a set $\cX$ represented as the linear image of a basic spectratope $\cW$:
$$
\begin{array}{c}
\cX=P\cW=\{x\in\bR^n:\exists w\in\cW:x=Pw\},\\
\cW=\Big\{w\in\bR^p: \exists t\in\cT:\Big[{\sum}_iw_iT_{ki}\Big]^2\preceq t_kI_{d_k},k\leq K\Big\}
\end{array}
$$
{We refer to $D=\sum_{k=1}^K d_k$ as {\em spectratopic size} of $\cW$ and $\cX$.}
\par
Same as ellitopes, spectratopes are convex compact sets symmetric w.r.t. the origin;  a basic spectratope, in addition, has a nonempty interior.
\subsubsection{Examples} First of all, {\sl every ellitope is a spectratope.} Indeed, it suffices to consider the case when the ellitope $\cW$ in question is the basic ellitope (\ref{2020ell1}). In this case, passing to eigenvalue decompositions of matrices $T_k$, we have
$$
T_k=\sum_{i=1}^{\nu_k}e_{ki}e_{ki}^T,\;\nu_k=\rank(T_k),
$$
whence
{\footnotesize$$
\begin{array}{rcl}
\cW&=&\{w:\exists t\in\cT: w^TT_kw\leq t_k,k\leq K\}\\
&=&\{w:\exists \bar{t}
=\{t_{ki},i\leq \nu_k,k\leq K\}\in\overline{\cT}:
[e_{ki}^Tw]^2\preceq t_{ki}I_1,\,1\leq i\leq \nu_k,k\leq K\},\\
\overline{\cT}&:=&\Big\{t_{ki}\geq0:\,\Big[\sum_it_{1i};...;\sum_it_{Ki}\Big]\in\cT\Big\}.
\end{array}
$$}\noindent
An example of a ``genuine'' basic spectratope is the unit $|\cdot|$-ball, $|\cdot|$ being the spectral norm on $\bR^{p\times q}$:
$$
\{w\in\bR^{p\times q}:|w|\leq 1\}=\bigg\{w:\exists t\in\cT[0,1]: T^2[w]:=\left[\begin{array}{c|c}&w\cr\hline w^T&\cr\end{array}\right]^2\preceq tI_{p+q}\bigg\},
$$
Same as ellitopes, spectratopes admit fully algorithmic "calculus," and their family is closed with respect to basic operations preserving convexity
and symmetry w.r.t. the origin, such as
taking finite intersections, linear images, inverse images under linear embedding, direct products, arithmetic summation
 ( see \cite[Section 4.6]{STOPT} for details); what is missing, is taking convex hulls of finite unions.
 \subsubsection{Bounding maximum of  quadratic form over a spectratope}\label{bounding_spectratope}
 Given a linear mapping\\
\centerline{$
 R[w]= {\sum}_{i=1}^\nu w_iR_i: \bR^\nu\to\bS^d
$}
so that $R_i\in\bS^d$, we associate with it  linear mappings
{\small\[\begin{array}{rclrcl}
R^+[W]&=&\sum_{i,j=1}^\nu W_{ij}R_iR_j:\;\bS^\nu\to\bS^d,&
R^{+,*}[\Lambda]&=&\left[\Tr(\Lambda R_iR_j)\right]_{i,j\leq\nu}:\;\bS^d\to\bS^\nu.
\end{array}\]}\noindent
Note that
\be R^+[ww^T]&=&R^2[w]\ee{add1}
and
\be\Tr(R^+[W]\Lambda)&=&\Tr(WR^{+,*}[\Lambda])\,\,\forall (W\in\bS^\nu,\Lambda\in\bS^d).
\ee{add2}
Given a collection $\Lambda=\{\Lambda_k,k\leq K\}$ of symmetric matrices (of, perhaps, different sizes), we set
$$
\lambda[\Lambda]=[\Tr(\Lambda_1);...;\Tr(\Lambda_K)].
$$
Finally, same as above, for a convex compact set $\cT$,
$$
\phi_\cT(\lambda)=\max_{t\in\cT}\lambda^Tt
$$
is the support function of $\cT$.
\par
Given a spectratope\\
\centerline{$
\cX=P\cW,\,\cW=\Big\{w\in\bR^q: \exists t\in\cT: T_k^2[w]:=\left[{\sum}_iw_iT_{ki}\right]^2\preceq t_kI_{d_k},k\leq K\Big\},
$}
an efficiently computable upper bound $\Opt(C)$ on the quantity
$$
\Opt_*(C)=\max_{x\in\cX}x^TCx
$$
can be built as follows. Assume that $\Lambda=\{\Lambda_k\in\S^{d_k}_+,k\leq K\}$ is such that
\be
P^TCP\preceq \sum_k{T}_k^{+,*}[\Lambda_k].
\ee{(a)}
When $x\in \cX$, there exists $w\in\bR^q$ and $t\in\cT$ such that (see (\ref{add1}))\\
\centerline{$
x=Pw\ \&\ T^+_k[ww^T]=T_k^2[w]\preceq t_kI_{d_k},\,k\leq K,
$}
whence
\be
\sum_k\Tr(T^+_k[ww^T]\Lambda_k)\leq \sum_kt_k\Tr(\Lambda_k)\leq \phi_{\cT}(\lambda[\Lambda]).
\ee{(b)}
On the other hand, by (\ref{add2}) we have
\be
\Tr(T_k^+[ww^T]\Lambda_k)=\Tr(T_k^{+,*}[\Lambda_k][ww^T])=w^TT_k^{+,*}[\Lambda_k]w,
\ee{(c)}
so that
$$
x^TCx=w^T[P^TCP]w\underbrace{\leq}_{\hbox{[by \rf{(a)}}]} w^T\left[{\sum}_kT_k^{+,*}[\Lambda_k]\right]w\leq\phi_{\cT}(\lambda[\Lambda])
$$
due to \rf{(c)} and \rf{(b)}.
As a result, the efficiently computable convex function
$$
\Opt(C)=\min\limits_{\Lambda}\left\{\phi_{\cT}(\lambda[\Lambda]):\Lambda=\{\Lambda_k\in\bS^{d_k}_+,k\leq K\},P^TCP\preceq {\sum}_kT_k^{+,*}[\Lambda_k]\right\}
$$
is an upper bound on $\Opt(C)$. It is known (\cite[Proposition 4.8]{STOPT}) that this bound is reasonably tight:
$$
\Opt_*(C)\leq \Opt(C)\leq 2\ln(2D)\Opt_*(C),\;\;D={\sum}_kd_k.
$$
\subsection{Bounding operator norms, spectratopic case}
Similarly to the ellitopic case, our current problem of interest is tight computationally efficient upper-bounding of the norm
$$
\|A\|_{\cB,\cX}=\max_{x\in\cX}\|Ax\|_{\cB}=\max\limits_{[y;x]\in\cB_*\times\cX}[y;x]^T\left[\begin{array}{c|c}&{1\over 2}A\cr\hline {1\over 2}A^T&\cr\end{array}\right][y;x]
$$
in the case when $\cX$ and $\cB_*$ are spectratopes:
\begin{equation}\label{s2020ell2*-**}
\begin{array}{rcl}
\cX&=&P\cW=\{x\in\bR^n:\exists w\in\cW:x=Pw\},\\
&&\cW=\{w\in\bR^p: \exists t\in\cT:T_k^2[w]\preceq t_kI_{d_k},k\leq K\}\\
\cB&=&\{v\in\bR^m:v^Ty\leq1\,\forall y\in\cB_*\},\,\cB_*=Q\cZ=\{y\in\bR^m:\exists z\in\cZ:y=Qz\},\\
&&\cZ=\{z\in\bR^q: \exists r\in\cR: R_\ell^2[z]\preceq s_\ell I_{g_\ell},\ell\leq L\}.\\
\end{array}\end{equation}
%where
%\begin{equation}\label{s2020ell2**}
%\begin{array}{rcl}
%\cB_*&=&\{y\in\bR^m:y^Tv\leq1\,\forall v\in\cB\}=Q\cZ=\{y\in\bR^m:\exists z\in\cZ:y=Qz\},\\
%\cZ&=&\{z\in\bR^q: \exists r\in\cR: R_\ell^2[z]\preceq s_\ell I_{g_\ell},\ell\leq L\}.
%\end{array}
%\end{equation}
In this case the efficiently computable upper bound on $\|A\|_{\cB,\cX}$ and its tightness are given by the following result (which is an improvement of the just cited result from \cite{STOPT}):
{\begin{theorem}\label{sverynewRelaxationTheorem} In the case of {\rm \rf{s2020ell2*-**}} the efficiently computable convex function of $A$ given by
{\small\begin{equation}\label{snormbound}
\Opt(A)=\min\limits_{\Lambda,\Upsilon}\left\{\phi_{\cT}(\lambda[\Lambda])+\phi_{{\cR}}(\lambda[\Upsilon]):
\begin{array}{l}\Lambda=\{\Lambda_k
\in
\bS^{d_k}_+,k\leq K\}, \Upsilon=\{\Upsilon_\ell\in\bS^{g_\ell}_+,\ell\leq L\}\\
\left[\begin{array}{c|c}\sum_\ell R_\ell^{+,*}[\Upsilon_\ell]&{1\over 2}Q^TAP\cr \hline{1\over 2}P^TA^TQ&\sum_kT_k^{+,*}[\Lambda_k]\cr\end{array}\right]\succeq0\\
\end{array}\right\}
\end{equation}}\noindent
is a reasonably tight upper bound on $\|A\|_{\cB,\cX}$:
\begin{equation}\label{sbetterbound}
\begin{array}{c}
\|A\|_{\cB,\cX}\leq\Opt(A)\leq \anc{2\sqrt{\ln(5G)\ln(5D)}}{\overline{\varsigma}\left(\sum_{k=1}^Kd_k\right)}
{\overline{\varsigma}\left(\sum_{\ell=1}^Lg_\ell\right)}\|A\|_{\cB,\cX},\\
\overline{\varsigma}(M)=\sqrt{2\ln(5M)}\\
\end{array}
\end{equation}
\end{theorem}}\noindent
{\bf Proof.}  {\bf 1$^0$} The left inequality in (\ref{sbetterbound}) is evident. Let us prove the right inequality. Let $q$, $p$ be the dimensions of the embedding spaces of $\cZ$ and $\cW$, and assume that $q\leq p$, which is w.l.o.g. for the same reasons as in the ellitopic case.
Same as in the latter case, (\ref{snormbound}) is nothing but the conic problem
$$
\Opt(A)=\min\limits_{\Lambda,\Upsilon,\tau,\theta}\left\{\tau+\theta:\begin{array}{l}\Lambda=\{\Lambda_k\in\bS^{d_k}_+,k\leq K\},
{[-\lambda[\Lambda];\tau]}\in{\mathfrak{T}}_*\\
\Upsilon=\{\Upsilon_\ell\in\bS^{g_\ell}_+,\ell\leq L\},[-\lambda[\Upsilon];\theta]\in{\mathfrak{R}}_*\\
\left[\begin{array}{c|c}
\sum_\ell R_\ell^{+,*}[\Upsilon_\ell]&{1\over 2}Q^TAP\cr\hline{1\over 2}P^TA^TQ&\sum_kT_k^{+,*}[\Lambda_k]\cr\end{array}\right]\succeq0\\
\end{array}\right\}
$$
with the same cones $\mathfrak{T}$, $\mathfrak{R}$ and their duals ${\mathfrak{T}}_*$, ${\mathfrak{R}}_*$ as in the ellitopic case.
Same as in that case, the latter problem is strictly feasible and bounded, and by Conic Duality one has
{\footnotesize$$
\begin{array}{rcl}
&&\Opt(A)\\
&=&\max\limits_{r,t,U,V,W}\left\{\Tr(W^TQ^TAP):\begin{array}{l}
{R}^+_\ell[U]\leq r_\ell I_{g_\ell}\,\forall \ell,
{T}^+_k[V]\leq t_kI_{d_k}\,\forall k\\t\in\cT,r\in\cR,\,
\left[\begin{array}{c|c}U&W\cr\hline W^T&V\cr\end{array}\right]\succeq0\\
\end{array}\right\}\\
&=&\max\limits_{r,t,U,V,Y}\left\{\Tr([U^{1/2}YV^{1/2}]^TQ^TAP):\begin{array}{l}r\in\cR,t\in\cT,U\succeq0,V\succeq0,Y^TY\preceq I\cr
{R}^+_\ell[U]\leq r_\ell I_{g_\ell}\,\forall \ell,{T}^+_k[V]\leq t_kI_{d_k}\,\forall k\\
\end{array}\right\}\\
&=&\max\limits_{r,t,U,V}\left\{\sum_{i=1}^q\sigma_i(U^{1/2}Q^TAPV^{1/2}):
\begin{array}{l}{R}^+_\ell[U]\leq r_\ell I_{g_\ell}\,\forall \ell,\;{T}^+[V]\leq t_kI_{d_k}\forall k\\
U\succeq0,\,V\succeq0,\,r\in\cR,\,t\in\cT\end{array}\right\}\\
\end{array}
$$}\noindent
(cf. item 1$^o$ in the ``ellitopic proof'').
\par{\bf 2$^o$.} The concluding optimization problem in the above chain clearly is solvable; let $U,V,r,t$ be the optimal solution,  and $\sum_{\iota=1}^q\sigma_\iota e_\iota f_\iota^T$ be the singular value decomposition of
$U^{1/2}Q^TAPV^{1/2}$, so that
\begin{equation}\label{ssothatmay2021}
\begin{array}{rcl}
\Opt(A)&=&\sum_{\iota=1}^q\sigma_\iota\\
U^{1/2}Q^TAPV^{1/2}&=&\sum_{\iota=1}^q\sigma_\iota e_\iota f_\iota^T\\
e_i^Te_j&=&\left\{\begin{array}{ll}1,&i=j\\
0,&i\neq j\cr\end{array}\right.,i,j\leq q\ \&\ f_i^Tf_j=\left\{\begin{array}{ll}1,&i=j\\
0,&i\neq j\cr\end{array}\right.,i,j\leq p.
\end{array}
\end{equation}
Let $\epsilon_1,...,\epsilon_p$ be independent random variables taking values $\pm1$ with probabilities 1/2, and let
$$\xi={\sum}_{i=1}^q\epsilon_ie_i,\;\;\; \eta={\sum}_{j=1}^p\epsilon_jf_j.$$
 Then, in view of (\ref{ssothatmay2021}) it holds, identically in $\epsilon_i=\pm1$, $1\leq i\leq p$:
\begin{equation}\label{sidentically}
\xi^TU^{1/2}Q^TAPV^{1/2}\eta={\sum}_{i,\iota\leq q,j\leq p}
[\epsilon_i\epsilon_j\sigma_\iota e_i^Te_\iota f_\iota^Tf_j]={\sum}_{\iota=1}^q\sigma_\iota=\Opt(A).
\end{equation}
On the other hand, setting $E=[e_1,...,e_q]$, we get an orthonormal $q\times q$ matrix such that $\xi=E\underline{\epsilon}$, where $\underline{\epsilon}=[\epsilon_1;...;\epsilon_q]$ is a Rademacher vector. Now let
$\widehat{\xi}=U^{1/2}\xi=U^{1/2}E\underline{\epsilon}$. Observe that for every $\ell\leq L$ and for properly selected matrices $\overline{R}_{\ell i}\in\bS^{g_\ell}$ we have
$$
R_\ell[U^{1/2}Ey]={\sum}_{i=1}^q\overline{R}_{\ell i}y_i,\,\forall y\in\bR^q.
$$
We have
$$\bE\{\widehat{\xi}\widehat{\xi}^T\}=\bE\{U^{1/2}E\underline{\epsilon}\underline{\epsilon}^TE^TU^{1/2}\}=U^{1/2}E\bE\{\underline{\epsilon}
\underline{\epsilon}^T\}E^TU^{1/2}=U^{1/2}EE^TU^{1/2}=U,$$
whence
$$
\bE\left\{R_\ell^2[\widehat{\xi}]\right\}=R^+_\ell[\bE\{\widehat{\xi}\widehat{\xi}^T\}]=R^+_\ell[U]\preceq r_\ell I_{g_\ell}\;\forall \ell\leq L
$$
(we have used \rf{add1}).
On the other hand,
$$
R_\ell[\widehat{\xi}]=R_\ell[U^{1/2}E\underline{\epsilon}]={\sum}_{i=1}^q\overline{R}_{\ell i}\epsilon_i,
$$
so that $\sum_{i=1}^q\overline{R}_{\ell i}^2=\bE\left\{R_\ell^2[\widehat{\xi}]\right\}$, and
we end up with
$$
{\sum}_{i=1}^q\overline{R}_{\ell i}^2\preceq r_\ell I_{g_\ell}.
$$
Applying the noncommutative Khintchine inequality\footnote{Noncommutative Khintchine Inequality due to Lust-Piquard, Pisier, and Buchholz, see \cite[Theorem 4.6.1]{Tropp112}, states that if  $Q_i\in \bS^n$, $1\leq i\leq I$, and $\xi_i$, $i=1,...,I$, are independent Rademacher  or $\cN(0,1)$ random variables, then for all $t\geq 0$ one has
$$ \Prob\left\{\left|{\sum}_{i=1}^I\xi_i Q_i\right|\geq t\right\}\leq
2n\exp\left\{-{t^2\over 2\left|{\sum}_{i=1}^IQ_i^2\right|}\right\}$$
where $|\cdot|$ is the  spectral norm.} we conclude that
{\small$$
\forall s>0: \Prob\{R_\ell^2[\widehat{\xi}]\preceq s^2r_\ell I_{g_\ell}\}=1-\Prob\Big\{\Big|{\sum}_{i}\overline{R}_{\ell i}\epsilon_i\Big|>s\sqrt{r_\ell}\Big\}\geq 1-2d_\ell\exp\big\{-\half s^2\big\}.
$$}\noindent
As a result, when setting $D=\sum_kd_k$ and $s=\sqrt{2\ln(5D)}$ we get
$$
\Prob\{R_\ell^2[\widehat{\xi}]\preceq 2\ln(5D)r_\ell I_{d_\ell},\ell\leq L\} > 1/2,
$$
and
$$
\Prob\big\{\|QU^{1/2}\xi\|_{\cB_*} \leq \sqrt{2\ln(5D)}\big\}\geq \Prob\{R_\ell^2[\widehat{\xi}]\preceq 2\ln(5D)r_\ell I_{d_\ell},\ell\leq L\}>1/2.
$$
By similar reasoning,
$$
\Prob\{\|PV^{1/2}\eta\|_\cX \leq \sqrt{2\ln(5G)}\}>1/2,\,\,G=\sum_\ell g_\ell.
$$
As a result, there exists realization $(\bar{\xi},\bar{\eta})$ of $(\xi,\eta)$  such that
$$
\|QU^{1/2}\bar{\xi}\|_{\cB_*} \leq \sqrt{2\ln(5D)}\; \&\; \|PV^{1/2}\bar{\eta}\|_{\cX}\leq \sqrt{2\ln(5G)}.
$$
On the other hand, invoking (\ref{sidentically}),
$$
\Opt(A)=\bar{\xi}^TU^{1/2}Q^TAPV^{1/2}\bar{\eta}\leq \|QU^{1/2}\bar{\xi}\|_{\cB*}\|PV^{1/2}\bar{\eta}\|_\cX\|A\|_{\cB,\cX}.
$$
Combining our observations, we conclude that
$$
\Opt(A)\leq 2\sqrt{\ln(5D)\ln(5G)}\|A\|_{\cB,\cX}.\eqno{\hbox{\qed}}
$$
%%%%%%%%%%%%%%%%%%%%%%%%%%%%%%
\subsection{Bounding robust norms of uncertain matrices, spectratopic case}\label{uncerftainmatrspectr}
Let spectratopes $\cX\subset\bR^n$, $\cB_*\subset\bR^m$ with nonempty interiors and the polar $\cB$ of $\cB_*$  be given by (\ref{s2020ell2*-**}). Our goal is to conceive a computationally efficient upper-bounding of the robust norm\\
\centerline{$
\|\cA\|_{\cB,\cX}=\max\limits_{A\in\cA}\|A\|_{\cB,\cX}
$}
of uncertain matrix\\
\centerline{$
\cA=\Big\{{\sum}_s \epsilon_sA_s:\|\epsilon\|_\infty\leq1\Big\}\subset\bR^{m\times n}.
$}
\subsubsection{Processing the problem}\label{Beginning}  Acting exactly as in the ellitopic case, with the results of Section \ref{bounding_spectratope} in the role of their ``ellitopic counterparts'' from Section \ref{bounding_ellitope},
we conclude that the efficiently computable quantity
{\footnotesize\begin{equation}\label{sssMay2021eq2}
\Opt:=\min\limits_{\Lambda,\Upsilon,\atop \{G_s,H_s\}}\left\{\phi_{\cR}(\lambda[\Upsilon])+\phi_{\T}(\lambda[\Lambda]):
\begin{array}{l}\left[\begin{array}{c|c}G_s&{1\over2}Q^TA_sP\cr\hline{1\over2} P^TA_s^TQ&H_s\cr\end{array}\right]\succeq0,s\leq S\\
\Upsilon=\{\Upsilon_\ell\in\bS^{g_\ell}_+,\ell\leq L\}, \,\sum_sG_s\preceq\sum_\ell R_\ell^{+,*}[\Upsilon_\ell]\\
\Lambda=\{\Lambda_k\in\bS^{d_k}_+,k\leq K\},\,\sum_sH_s\preceq \sum_k T_k^{+,*}[\Lambda_k]
\end{array}\right\}
\end{equation}}\noindent
---the ``spectratopic analog'' of (\ref{May2021eq2-1})---is an upper bound on $\|\cA\|_{\cB,\cX}$ such that for properly selected matrices $X\in\bS^p_+$, $Y\in\bS^q_+$ and $r\in\cR$, $t\in\cT$ one has\\
\centerline{$
R^+_\ell[Y]\preceq r_\ell I_{g_\ell},\,\ell\leq L,\ \& \ T^+_k[X]\preceq t_k I_{d_k},\,k\leq K,
$}
and {for the norms $\pi(\cdot)$ and  $\rho(\cdot)$ with unit balls $\cW$ and $\cZ$, respectively, and } $[\eta;\xi]\sim\cN(0,\mathrm{Diag}\{Y,X\})$,
\begin{equation}\label{sssarriveat}
\Opt\leq\vartheta(2\kappa)\|\cA\|_{\cB,\cX}\bE\{\rho(\eta)\pi(\xi)\}=
\vartheta(2\kappa)\|\cA\|_{\cB,\cX}\bE\left\{\pi(\xi)\right\}\bE\left\{\rho(\eta)\right\}
\end{equation}
where $\kappa$ is the maximum of ranks of $A_s$ and  $\vartheta(\cdot)$ is given by (\ref{theta}) (cf.(\ref{arriveat})).
\par
We have the following spectratopic analog of Lemma \ref{lenlemlem}.
{\begin{lemma}\label{sslenlemlem}
Let\\
\centerline{$
\cV=\{v\in\bR^d: \exists r\in\cR: R_j^2[v]\preceq r_jI_{\nu_j},1\leq j\leq J\}\subset\bR^d
$}
be a basic spectratope, $W\succeq 0$ be symmetric $d\times d$ matrix such that\\
\centerline{$
\exists r\in\cR: R^+_j[W]\preceq  r_j I_{\nu_j},j\leq J,
$}
and $\omega\sim\cN(0,W)$. Denoting by $\gamma(\cdot)$ the norm on $\bR^d$ with the unit ball $\cV$, we have
\begin{equation}\label{ssMay2021upsilon}
\bE\{\rho(\omega)\}\leq \overline{\varkappa}\left({\sum}_j\nu_j\right),\;
\overline{\varkappa}(F)=2\sqrt{2\ln(2F)}.\\
\end{equation}
\end{lemma}}\noindent
{\bf Proof.}  Let
$\zeta\sim\cN(0,I_d)$. When setting
$$
\overline{R}_j[z]=R_j[W^{1/2}z]={\sum}_{i=1}^d\overline{R}_{ji}z_j,\eqno{[\overline{R}_{ji}\in\bS^{\nu_j}]}, \;j\leq J,
$$
we have
$$
{\sum}_i\overline{R}_{ji}^2=\bE\{\overline{R}_j^2[\zeta]\}=\bE\{R_j^2[W^{1/2}\zeta]\}
=\bE\{R^+_j[W^{1/2}\zeta\zeta^TW^{1/2}]\}=R^+_j[W]\preceq r_jI_{\nu_j}.
$$
Hence for every $s>0$
{\small$$
\begin{array}{rcl}
\Prob\left\{R_j^2[\omega]\preceq s^2r_jI_{\nu_j}\right\}&=&
\Prob\left\{\overline{R}_j^2[\zeta]\preceq s^2r_jI_{\nu_j}\right\}
=1-\Prob\left\{\big|\sum_i\zeta_i\overline{R}_{ji}\big|>s\sqrt{r_j}\right\}\\
&&\multicolumn{1}{r}{\hbox{[as above, $|\cdot|$ is spectral norm]}}\\
&\geq &1-2\nu_j\exp\{-s^2/2\},\\
\end{array}
$$}\noindent
with the concluding $\geq$ given by $\sum_i\overline{R}_{ji}^2\preceq r_jI_{\nu_j}$ combined with the noncommutative Khintchine inequality. As a result,
$$
\Prob\{\gamma(\omega)>s\}\leq 1-\Prob\left\{\exists j: R_j^2[\omega]\preceq s^2r_jI_{\nu_j}\right\}\leq \Big[{\sum}_j2\nu_j\Big]\exp\{-s^2/2\}.
$$
Therefore, when setting $F=\sum_j\nu_j$ we obtain
$$
\bE\{\gamma(\omega)\}\leq \int_0^\infty\min\left[1,2F\exp\{-\gamma^2/2\}\right]d\gamma\leq 2\sqrt{2\ln(2F)}.\eqno{\hbox{\qed}}
$$
Applying the lemma to $\cV=\cW$, $W=X$, and to $\cV=\cZ$, $W=Y$, we get from (\ref{sssarriveat}) the following analog of
Proposition \ref{propmatrcube}:
{\begin{proposition}\label{sspropmatrcube}
In the situation described in the beginning of this section, assuming that ranks of all $A_s$ are $\leq \kappa$, the efficiently computable quantity $\Opt$ as given by {\rm (\ref{sssMay2021eq2})} is a reasonably tight upper bound
on the robust norm $\|\cA\|_{\cB,\cX}$ of uncertain matrix $\cA$, specifically,
\begin{equation}\label{sswehavethat}
\|A\|_{\cB,\cX}\leq\Opt\leq \overline{\varkappa}({\sum}_kd_k)\overline{\varkappa}({\sum}_\ell g_\ell)\vartheta(2\kappa)\|A\|_{\cB,\cX}
\end{equation}
where $\overline{\varkappa}(\cdot)$ is given by  {\rm (\ref{ssMay2021upsilon})} and $\vartheta(\cdot)$, given by {\rm (\ref{theta})}, satisfies
\[
\vartheta(1)=1,\,\vartheta(2)={\pi\over 2},\,\vartheta(4)=2,\,\vartheta(k)\leq \pi\sqrt{k}/2.
\]\end{proposition}}
\subsubsection{Putting things together} Results of  Proposition \ref{sspropmatrcube} (and as a byproduct -- of Theorem \ref{sverynewRelaxationTheorem}) can be extended, in exactly the same fashion
as in the ellitopic case, to the situation where  $\cX$ and the polar $\cB_*$ of $\cB$ are convex hulls of finite unions of spectratopes rahter than plain spectratopes, and the uncertain matrix
in question is not centered, resulting in the following spectratopic analogy of Proposition \ref{propmixture}:

{\begin{theorem}\label{spropmixture} Let $\cU=\{A_\n+\sum_{s=1}^S\epsilon_sA_s:\|\epsilon\|_\infty\leq1\}$ be an uncertain $m\times n$ matrix, $\cX\subset\bR^n$, $\cB,\cB_*\subset\bR^m$ be given by
$$
\begin{array}{rcl}
\cX&=&\Conv\{{\bigcup}_{i=1}^IP_i\cX_i\}=\left\{x=\sum_{i=1}^I \lambda_iP_ix_i:x_i\in\cX_i,\lambda_i\geq0,\sum_i\lambda_i=1\right\}\\
\cB&=&\{v\in\bR^m:\max\limits_{y\in\cB_*} v^Ty\leq 1\},\,\cB_*=\Conv\{{\bigcup}_{J=1}^JQ_j\cZ_j\}\\
&=&\Big\{y={\sum}_{j=1}^J\mu_jQ_jz_j,z_j\in\cZ_j,\mu_j\geq0,{\sum}_j\mu_j=1\Big\}
\end{array}
$$ with basic spectratopes
{\footnotesize$$
\begin{array}{rcl}
\cX_i&=&\left\{x_i\in\bR^{\nu_i}:\exists t^i\in\cT^i: T^2_{ki}[x_i]^2\preceq t^i_kI_{d_{ki}},1\leq k\leq K_i\right\},\,T_{ki}[x]=\sum_{\iota=1}^{\nu_i} x_\iota T_{ki\iota},\,i\leq I\\
\cZ_j&=&\left\{z_j\in\bR^{\mu_j}:\exists r^j\in\cR^j: R^2_{\ell j}[z_j]\preceq r^j_\ell I_{g_{\ell j}},1\leq \ell\leq L_,\right\},R_{\ell j}[z]=
\sum_{\iota=1}^{\mu_j} z_\iota R_{\ell j\iota}\,j\leq J\\
\end{array}
$$}\noindent
Then
the efficiently computable quantity
\[
\Opt[\cU]=
{\max}_{i\leq I,j\leq J}\Opt_{ij}[\cU],
\]
where
{\small
$$
\begin{array}{rcl}
\Opt_{ij}[\cU]&=&\min\limits_{\Lambda^{ij},\Upsilon^{ij},G^{ijs},H^{ijs}\atop
{\overline{G}^{ij},\overline{H}^{ij}\atop
1\leq i\leq I,1\leq j\leq J,1\leq s\leq S}}\Big\{
\phi_{\cT^i}(\lambda[\Upsilon^{ij}])+{\phi}_{\cR^j}(\lambda[\Lambda^{ij}]):\\
&&\left.\begin{array}{l}\Lambda^{ij}=\{\Lambda^{ij}_k\succeq0,k\leq K_i\},\,\Upsilon^{ij}=\{\Upsilon^{ij}_\ell\succeq0,\ell\leq L_j\}\\
{\sum}_{s=1}^SH^{ijs}+\overline{H}^{ij}\preceq \sum_{k=1}^{K_i}T^{+,*}_{ki}[\Lambda^{ij}_k]\\
{\sum}_{s=1}^SG^{ijs}+\overline{G}^{ij}\preceq {\sum}_{\ell=1}^{L_j}R^{+,*}_{\ell j}[\Upsilon^{ij}_\ell],\\
\left[\begin{array}{c|c}G^{ijs}&{1\over 2}[Q_j^TA_sP_i]\cr\hline
{1\over 2}[Q_j^TA_sP_i]^T&H^{ijs}\cr\end{array}\right]\succeq0,s\leq S\\ \left[\begin{array}{c|c}\overline{G}^{ij}&{1\over 2}[Q_j^TA_\n P_i]\cr\hline
{1\over 2}[Q_j^TA_\n P_i]^T&\overline{H}^{ij}\cr\end{array}\right]\succeq0\end{array}\right\}, i\leq I,j\leq J
\Big\}
\end{array}
$$}\noindent
is an efficiently computable {\em convex} in $(A_\n,A_1,...,A_S)$ upper bound on $\|\cU\|_{\cB,\cX}$. This upper bound is reasonably tight,
specifically, setting
$$
\cU_{ij}=Q_j^TA_\n P_i+\Big\{{\sum}_{s=1}^S\epsilon_s[Q_j^TA_sP_i]:\|\epsilon\|_\infty\leq1\Big\},
$$
we have
$$
\begin{array}{c}
\|\cU_{ij}\|_{\cZ_j^*,\cX_i}\leq\Opt_{ij}[\cU]\leq [\overline{\varsigma}\left(D_i\right)\overline{\varsigma}
\left(G_j\right)+\overline{\varkappa}\left(D_i\right)\overline{\varkappa}\left(G_j\right)\vartheta(2\kappa)] \|\cU_{ij}\|_{\cZ_j^*,\cX_i},
\\
D_i=\sum_{k=1}^{K_i}d_{ki},\,\,G_j=\sum_{\ell=1}^{L_j}g_{\ell j}\\
\end{array}
$$
and
$$%\begin{equation}\label{seqmixture2}
\begin{array}{rcl}
\|\cU\|_{\cB,\cX}&=&\max\limits_{i\leq I, j\leq J}\|\cU_{ij}\|_{\cZ_j^*,\cX_i}\leq \Opt[\cU]=\max\limits_{i\leq I, j\leq J}\Opt_{ij}[\cU]\\
&\leq&
\left[\max\limits_{i\leq I, j\leq J}[\overline{\varsigma}\left(D_i\right)\overline{\varsigma}
\left(G_j\right)+\overline{\varkappa}\left(D_i\right)\overline{\varkappa}\left(G_j\right)\vartheta(2\kappa)]\right]\|\cU\|_{\cB,\cX}\\
\end{array}
$$
where $\kappa$ is the maximum of ranks of $A_s$, $1\leq s\leq S$, $\overline{\varsigma}(\cdot)$ and $\overline{\varkappa}(\cdot)$ are as defined in (\ref{sbetterbound}) and
(\ref{ssMay2021upsilon}), and
$\vartheta(\cdot)$  is defined by (\ref{theta}) and satisfies (\ref{4.qwer8}).
\end{theorem}}\noindent
%%%%%%%%%%%%%%%%%%%%%%%%%%%%%%%%%%%%%%%%%%%%%%%%%%%%%%%%%%%%%%%%%%%%%%%%%%%%%%%%%%%%%%%%%%%%%%%%%%%%%%%%%
%%%%%%%%%%%%%%%%%%%%%%%%%%%%%%%%%%%%%%%%%%%%%%%%%%%%%%%%%%%%%%%%%%%%%%%%%%%%%%%%%%%%%%%%%%%%%%%%%%%%%%%%%
%%%%%%%%%%%%%%%%%%%%%%%%%%%%%%%%%%%%%%%%%%%%%%%%%%%%%%%%%%%%%%%%%%%%%%%%%%%%%%%%%%%%%%%%%%%%%%%%%%%%%%%%%
\hide
{ Complete derivation of (\ref{sssarriveat}):\\
With the same
 $\mathfrak{T}$, $\mathfrak{R}$, $\phi_\cT$, $\phi_\cR$ as in item 1$^o$ of the proof of Theorem \ref{verynewRelaxationTheorem}, we have
\begin{equation}\label{ss2020ell2*}
\begin{array}{rcl}
\cX&=&P\cW=\{x\in\bR^n:\exists w\in\cW:x=Pw\},\\
\cW&=&\{w\in\bR^p: \exists t\in\cT:T_k^2[w]\preceq t_kI_{d_k},k\leq K\}\\
\cB&=&\{v\in\bR^m:v^Ty\leq1\,\forall y\in\cB_*\},\\
\cB_*&=&\{y\in\bR^m:y^Tv\leq1\,\forall v\in\cB\}=Q\cZ=\{y\in\bR^m:\exists z\in\cZ:y=Qz\},\\
\cZ&=&\{z\in\bR^q: \exists r\in\cR: R_\ell^2[z]\preceq s_\ell I_{g_\ell},\ell\leq L\},
\\
\end{array}
\end{equation}
{\scriptsize\begin{equation}\label{ssMay2021eq2}
\begin{array}{l}
\Opt:=\min\limits_{\Lambda,\Upsilon,G_s,H_s}\left\{\phi_{\cR}(\lambda[\Lambda])+\phi_{\T}(\lambda[\Upsilon]):
\left\{\begin{array}{l}\left[\begin{array}{c|c}G_s&{1\over2}Q^TA_sP\cr\hline{1\over2} P^TA_s^TQ&H_s\cr\end{array}\right]\succeq0,s\leq S\\
\Lambda=\{\Lambda_\ell\in\S^{g_\ell}_+,\ell\leq L\},\Upsilon=\{\Upsilon_k\in\S^{d_k},k\leq K\}\\
\sum_sG_s\preceq\sum_\ell\cR_\ell^*[\Lambda_\ell],
\sum_sH_s\preceq \sum_k \cT_k^*[\Upsilon_k]\\
\end{array}\right.\right\}\\
=\min\limits_{\Lambda,\Upsilon,G_s,H_s,\varrho,\varsigma}\left\{\varrho+\varsigma:\left\{\begin{array}{l}\left[\begin{array}{c|c}G_s&{1\over2} Q^TA_sP\cr\hline{1\over2} P^T A_s^TQ&H_s\cr\end{array}\right]\succeq0,s\leq S,
{[-\lambda;\varrho]}\in\bR_*,[-\mu;\varsigma]\in\bT_*\\
\Lambda=\{\Lambda_\ell\in\S^{g_\ell}_+,\ell\leq L\},\Upsilon=\{\Upsilon_k\in\S^{d_k},k\leq K\}\\
\sum_sG_s\preceq\sum_\ell\cR_\ell^*[\Lambda_\ell],
\sum_sH_s\preceq \sum_k \cT_k^*[\Upsilon_k]\\
\end{array}\right.\right\}\\
=\max\limits_{Y,X,W_s,s,t}\left\{\sum_s\Tr(W_s^TQ^TA_sP):\left\{\begin{array}{l}t\in\T,r\in\cR\\
\cR_\ell[Y]\preceq s_\ell I_{g_\ell},\ell\leq L,\, \cT_k[X]\preceq t_k I_{d_k},k
\leq K\\
\left[\begin{array}{c|c} Y&W_s\cr\hline W_s^T&X\cr\end{array}\right]\succeq0\,\forall s\leq S\\
\end{array}\right.\right\}\\
\multicolumn{1}{r}{\hbox{[conic duality]}}\\
=\color{magenta} \max\limits_{Y,X,s,t}\left\{\sum_s\|\sigma(Y^{1/2}Q^TA_sPX^{1/2})\|_1:\left\{\begin{array}{l}
Y\succeq0,X\succeq0,
t\in\T,r\in\cR,\\
\cR_\ell[Y]\preceq s_\ell I_{g_\ell},\ell\leq L,\, \cT_k[X]\preceq t_k I_{d_k},k
\leq K\\
\end{array}\right.\right\}\\
=\max\limits_{Y,X,s,t}\left\{\sum_s\|\lambda(\cL[Y^{1/2}Q^TA_sPX^{1/2}])\|_1:\left\{\begin{array}{l}
Y\in\bS^q_+,X\in\bS^p_+,t\in\T,r\in\cR\\
\cR_\ell[Y]\preceq s_\ell I_{g_\ell},\ell\leq L,\, \cT_k[X]\preceq t_k I_{d_k},k
\leq K\\
\end{array}\right.\right\}\\
\end{array}
\end{equation}}\noindent
where $\sigma(A)$ is the vector of singular values of a matrix $A$, $\lambda(A)$ is the vector of eigenvalues of a symmetric matrix $A$, and
$$
\cL[B]=\left[\begin{array}{c|c}&{1\over2} B\cr\hline {1\over2} B^T&\cr\end{array}\right],
$$
\begin{quote}
{\small
In (\ref{ssMay2021eq2}),  $(a)$ follows from the two simple observations (cf. proof of Theorem \ref{verynewRelaxationTheorem}):
\begin{itemize}
\item {LMI {\tiny$\left[\begin{array}{c|c}P&Q\cr\hline Q^T&R\cr\end{array}\right]\succeq0$} with $p\times p$ matrix $P$ and $r\times r$ matrix $R$ takes place if and only if $P\succeq0$, $R\succeq0$, and $Q=P^{1/2}YR^{1/2}$ with $p\times r$ matrix $Y$ such that $Y^TY\preceq I_r$}, and
\item {for $p\times r$ matrix $A$, one has $\max\limits_Y\{\Tr(Y^TA): Y\in\bR^{p\times r}, Y^TY\preceq I_r\}=\|\sigma(A)\|_1$}
\end{itemize}
while $(b)$ stems from the fact that the eigenvalues of $2\cL[B]$ are  positive singular values of $B$, minus these positive singular values, and a number of zeros.}
\par
Note that $\Opt$ as defined in (\ref{ssMay2021eq2}) clearly is a convex function of $[A_1,...,A_S]$.
\end{quote}
Observe that $\|\cA\|_{\cB,\cX}\leq\Opt$. Indeed, the problem specifying $\Opt$ clearly is solvable, and if $\Lambda,\Upsilon,G_s,H_s$ is its optimal solution, we have for all
 $z\in\cZ$, $w\in\cW$, $\epsilon_s=\pm1:$
{\scriptsize$$
\begin{array}{l}
\epsilon_sz^TQ^TA_sPw\leq z^TG_sz+w^TH_sw\Rightarrow \sum_s\epsilon_sz^TQ^TA_sPw\leq z^T\left[\sum_\ell\cR^*_\ell[\Lambda_\ell]\right]z+w^T\left[\sum_k\cT_k^*[\Upsilon_k]\right]w\\
=\sum_\ell\Tr\left(\cR^*_\ell[\Lambda_\ell][zz^T]\right)z+\sum_k\Tr\left(\cT_k^*[\Upsilon_k][ww^T]\right)=\sum_\ell\Tr\left(\Lambda_\ell\cR_\ell[zz^T]\right)+
\sum_k\Tr\left(\Upsilon_k\cT_k[ww^T]\right)\\
=\sum_\ell\Tr\left(\Lambda_\ell R_\ell^2[z]\right)+
\sum_k\Tr\left(\Upsilon_k T_k^2[w]\right)\\
\leq \phi_{\cR}(\lambda[\Lambda])+\phi_{\cT}(\lambda[\Upsilon])
\hbox{\ [recall that $R_\ell^2[z]\preceq s_\ell I_{g_\ell}$, $\ell\leq L$, and $T_k^2[w]\preceq t_k I_{d_k}$, $k\leq K$,}\\ \hbox{ for some
$r\in\cR$ and $t\in\cT$, while $\Lambda_\ell\succeq0$ and $\Upsilon_k\succeq0]$}\\
\Rightarrow
\sum_s\epsilon_sz^TQ^TA_sPw\leq \max_{r\in\cR,t\in\T}\left[\lambda^Ts+\mu^Tt\right]\leq \phi_\cR(\upsilon)+\phi_{\T}(t)=\Opt.
\end{array}
$$
This relation holds true for all $w\in\cW,z\in\cZ$ and all $\epsilon_s=\pm1$, implying that $\|\cA\|_{\cB,\cX}\leq\Opt$ (recall that $P\cW=\cX,Q\cZ=\cB_*$).
\par
Now let $X\succeq0$, $Y\succeq0$, $t,s$ be such that $t\in\T$, $r\in\cR$,  $\cR_\ell[Y]\preceq s_\ell I_{g_\ell},$ $\ell\leq L$, $\cT_k[X]\preceq t_k I_{d_k},$ $k
\leq K$, and
$$
\Opt=\sum_s\|\lambda(\cL[Y^{1/2}Q^TA_sPX^{1/2}])\|_1.
$$
By \cite[Lemma 2.2]{BTNMC} (a.k.a. \cite[Lemma 3.4.3]{LMCO}), if the ranks of all matrices $A_s$ do not exceed a given $\kappa$, which we assume from now on, then for $\omega\sim\cN(0,I_{m+n})$ one has
$$
\bE\left\{|\omega^T\cL[Y^{1/2}A_sX^{1/2}]\omega|\right\}\geq \|\lambda(\cL[Y^{1/2}A_sX^{1/2}])\|_1/\vartheta(2\kappa),
$$
where $\vartheta(k)$ is the universal function given by
\begin{equation}\label{sstheta}
\vartheta(k)={1\over \min\left\{\displaystyle{\int}
|\alpha_1u_1^2+...+\alpha_ku_k^2|p_k(u)du\big|\, \alpha\in{\Bbb
R}^k,\|\alpha\|_1=1\right\}},
\end{equation}
where $p_k(\cdot)$ is the density of $\cN(0,I_k)$ and satisfying (\ref{4.qwer8}).   It follows that
{\scriptsize$$
\begin{array}{rcl}
\Opt&\leq& \vartheta(2\kappa)\bE_{\omega\sim\cN(0,I_{m+n})}\left\{\sum_s|\omega^T\cL[Y^{1/2}Q^TA_sPX^{1/2}]\omega|\right\}\\
&=&
\vartheta(2\kappa)\bE_{[\eta;\xi]\sim\cN(0,\hbox{\scriptsize Diag}\{Y,X\})}\left\{\sum_s|\eta^TQ^TA_sP\xi|\right\}.\\
\end{array}
$$}\noindent
Now let $\pi(\cdot)$ be the norm on $\bR^p$ with the unit ball $\cW$, and $\theta(\cdot)$ be the norm on $\bR^q$ with the unit ball $\cZ$; then, taking into account that $\cX=P\cW$, $\cB_*=Q\cZ$,
$$
\forall (\eta\in\bR^p,\xi\in\bR^q): \sum_s|\eta^TQ^TA_sP\xi|=\max_{\epsilon_s=\pm1}\eta^TQ^T[\sum_s\epsilon_sA_s]P\xi\leq \theta(\eta)\pi(\xi)\|\cA\|_{\cB,\cX},
$$
and we arrive at the relation
{\scriptsize$$
\begin{array}{rcl}
\Opt&\leq&\vartheta(2\kappa)\|\cA\|_{\cB,\cX}\bE_{[\eta;\xi]\sim\cN(0,\hbox{\scriptsize Diag}\{Y,X\})}\{\theta(\eta)\pi(\xi)\}\\
&=&
\vartheta(2\kappa)\|\cA\|_{\cB,\cX}\bE_{\xi\sim\cN(0,X)}\left\{\pi(\xi)\right\}\bE_{\eta\sim\cN(0,Y)}\left\{\theta(\eta)\right\}\\
\end{array}
$$}
which is (\ref{sssarriveat}).}
}%hide

\end{document}